\documentclass[11pt,twoside]{article}

\usepackage{amsbsy,amsfonts,amsmath,amssymb,eucal,mathrsfs}
\usepackage[all]{xy}
\usepackage{pstricks}
\usepackage{pst-plot}

\renewcommand{\ge}{\geqslant}
\renewcommand{\le}{\leqslant}

\renewcommand{\tilde}{\widetilde}

\renewcommand{\bar}{\overline}

\newcommand{\Z}{\mathbb Z}

\newcommand{\C}{\mathbb{C}}

\newcommand{\F}{\mathbb{F}}
\newcommand{\G}{\mathbb{G}}

\newcommand{\com}{\mathtt{Com}}

\newcommand{\p}{{\mathbb P}}

\renewcommand{\iff}{if and only if }

\newcommand{\trace}{\mbox{tr}}
\newcommand{\scal}[1]{\langle #1 \rangle}
\newcommand{\im}{\mathtt{Im}}
\newcommand{\lpara}{
\ \vspace{3pt}

\noindent}

\newcommand{\sectionplus}[1]{\section{#1} \vspace{-5mm} \indent}
\newcommand{\subsectionplus}[1]{\subsection{#1} \vspace{-5mm} \indent}
\newcommand{\subsubsectionplus}[1]{\subsubsection{#1} \vspace{-5mm} \indent}

\newcommand{\liste}{\

\begin{itemize}}
\newcommand{\codim}{\mbox{codim}}

\newtheorem{definition}{Definition}[section]
\newenvironment{defi}{\begin{definition} \rm}{\end{definition}}

\newtheorem{prop}[definition]{Proposition}

\newtheorem{lemm}[definition]{Lemma}
\newtheorem{fact}[definition]{Fact}
\newtheorem{coro}[definition]{Corollary}
\newtheorem{theo}[definition]{Theorem}
\newtheorem{notation}[definition]{Notation}

\newtheorem{construction}[definition]{Construction}

\newtheorem{remark}[definition]{Remark}
\newenvironment{rema}{\begin{remark} \rm}{\end{remark}}
\newtheorem{remarks}[definition]{Remarks}

\newtheorem{example}[definition]{Example}

\newtheorem{examples}[definition]{Examples}

\newtheorem{nothing}[definition]{$\!\!$}

\newenvironment{proo}{{\flushleft \bf Proof :}}{\hfill $\square$ \vspace{2mm}}

\newtheorem{definition*}{Definition}[section]
\newenvironment{defi*}{\begin{definition*} \rm}{\end{definition*}}
\newtheorem{definitions*}[definition*]{Definitions}
\newenvironment{defis*}{\begin{definitions*} \rm}{\end{definitions*}}
\newtheorem{prop*}[definition*]{Proposition}
\newtheorem{lemm*}[definition*]{Lemma}
\newtheorem{coro*}[definition*]{Corollary}
\newtheorem{theo*}[definition*]{Theorem}
\newtheorem{remark*}[definition*]{Remark}
\newenvironment{rema*}{\begin{remark*} \rm}{\end{remark*}}
\newtheorem{remarks*}[definition*]{Remarks}
\newenvironment{remas*}{\begin{remarks*} \rm}{\end{remarks*}}
\newtheorem{example*}[definition*]{Example}
\newenvironment{exam*}{\begin{example*} \rm}{\end{example*}}
\newtheorem{examples*}[definition*]{Examples}
\newenvironment{exams*}{\begin{examples*} \rm}{\end{examples*}}
\newtheorem{nothing*}[definition*]{$\!\!$}
\newenvironment{noth*}{\begin{nothing*} \rm}{\end{nothing*}}

\DeclareMathOperator{\OO}{O}

\newfont{\joli}{cmfi10 scaled \magstep2}

\DeclareMathAlphabet{\mathpzc}{OT1}{pzc}{m}{it}

  \def\cH{{\cal H}}

\def\cP{{\cal P}}   
   \def\cZ{{\cal Z}}

\def \b {{\beta}}
 \def \ga {{\gamma}}
 \def \wt {{\widetilde{w}}}

\def \uh {{\widehat{u}}}
\def \vh {{\widehat{v}}}
\def \wh {{\widehat{w}}}

 \def\s{\sigma}\def\OO{\mathbb O}\def\PP{\mathbb P}\def\QQ{\mathbb Q}
 \def\CC{\mathbb C} \def\ZZ{\mathbb Z}\def\JO{{\mathcal J}_3(\OO)}
 
\begin{document}

 \title{Quantum cohomology of minuscule homogeneous spaces}
 \author{P.E. Chaput, L. Manivel, N. Perrin}
 \maketitle

 {\def\thefootnote{\relax}
 \footnote{ \hspace{-6.8mm}
 }
 }

 \begin{abstract}
We study the quantum cohomology of (co)minuscule homogeneous varieties
under a unified perspective. We show that three points Gromov-Witten 
invariants can always be interpreted as classical intersection numbers
on auxiliary homogeneous varieties. 
Our main combinatorial tools are certain
quivers, in terms of which we obtain a quantum Chevalley formula and a
higher quantum Poincar{\'e} duality. 
In particular we compute the quantum 
cohomology of the two exceptional minuscule homogeneous varieties. 
 \end{abstract}

 \sectionplus{Introduction}
 
The quantum cohomology of complex homogeneous spaces has been studied 
by many people since the fundamental works of Witten and
Kontsevich-Manin -- see \cite{fulton0} for a survey. The (small) 
quantum cohomology ring of Grassmannians 
was investigated by Bertram with the help of Grothendieck's Quot
schemes \cite{bertram}, a method which can be applied in other
situations but can be quite technical. 
An important breakthrough, which greatly simplified
the computations of the quantum products in many cases, was 
the observation by Buch that the Gromov-Witten invariants of certain
homogeneous spaces are controlled by classical intersection numbers,
but on certain auxiliary homogeneous varieties \cite{Buch}. 

In this paper we make use of this idea in the general context of 
minuscule and cominuscule homogeneous varieties (see the next section 
for the definitions). We give a unified 
treatment of the quantum cohomology ring of varieties including 
ordinary and Lagrangian Grassmannians, spinor varieties, quadrics and
also the two exceptional Hermitian symmetric spaces -- the Cayley
plane $E_6/P_1$ and the Freudenthal variety $E_7/P_7$ (in all the 
paper we use the notations of \cite{bou} for root systems).
One of the conclusions of our study is the following presentation
of the quantum cohomology algebras.

\begin{theo}
Let $X=G/P$ be a minuscule homogeneous variety. There exists
a minimal homogeneous presentation of its integer cohomology of the form
$$H^*(X)=\ZZ [H,I_{p_1+1},\ldots ,I_{p_n+1}]/(R_{q_1+1},\ldots ,R_{q_r+1}),$$
and one can choose the relation of maximal degree $R_{q_r+1}$ so that 
the quantum cohomology ring of $X$ is
$$QH^*(X)=\ZZ[q,H,I_{p_1+1},\ldots ,I_{p_n+1}]
/(R_{q_1+1},\ldots ,R_{q_r+1}+q).$$
\end{theo}
Note that $q_r=h$ is the Coxeter number of $G$. 

In fact this was already known for the classical minuscule 
homogeneous
varieties (see \cite{fulton,KT1,KT2,siebtian}), so our contribution to this
statement only concerns the exceptional cases, which are treated in 
section \ref{section_quantum_exceptional} 
(see Theorems \ref{e6_quantique} and \ref{e7_quantique}). It is easy
to deduce from these theorems the quantum product of any two Schubert
classes in the two exceptional Hermitian symmetric spaces. This relies
on the computations of the classical cohomology rings done in
\cite{iliev} for the Cayley plane, and very recently in \cite{nikolenko}
for the Freudenthal variety. 

Our treatment of these
homogeneous spaces relies on combinatorial tools that we 
develop in the general context of (co)minuscule varieties. 
Namely, we use the combinatorics of certain quivers, first introduced
in \cite{carquois}, to give a convenient interpretation of Poincar{\'e}
duality on a  (co)minuscule $X=G/P$. We deduce a nice combinatorial 
version of the quantum Chevalley formula (Proposition \ref{chevalley}).
 
As we mentionned, this relies on the interpretation of degree one
Gromov-Witten invariants as classical intersection numbers on the 
Fano variety of lines on $X$, which remains
$G$-homogeneous. Generalizing the case by case analysis of \cite{BKT}
for classical groups, we extend 
this interpretation to degree $d$ invariants in Corollary
\ref{gw-d}: Gromov-Witten invariants of degree $d$ are classical 
intersection numbers on certain auxiliary $G$-homogeneous varieties $F_d$. 

In fact this statement gives a very special role to a certain Schubert 
subvariety $Y_d$ of $X$. In particular, we are able to define, in 
terms of this variety, a duality on a certain
family of Schubert classes. In degree zero this is just the usual 
Poincar{\'e} duality, which we thus extend to a ``higher quantum Poincar{\'e}
duality'', see Proposition \ref{q-poincare}. 
We also obtain a combinatorial
characterization, again in terms of quivers, of the minimal power 
of $q$ that appears in the quantum product of two Schubert classes:
see Corollary \ref{smallest}. 
This extends a result of Buch \cite{Buch}, which itself relied on a 
general study of this question by Fulton and Woodward \cite{FW}.

\tableofcontents

\sectionplus{Minuscule and cominuscule varieties} 
\label{section_minuscule}

We begin by reminding what are the (co)minuscule homogeneous
varieties, and how their Chow ring can be computed, in particular 
in the two exceptional cases. Then we introduce quivers in relation 
with Poincar{\'e} duality. These tools will be useful for a unified
approach of quantum cohomology. 

\subsectionplus{First properties}

Let $G$ be a simple complex algebraic group, and $B$ a Borel subgroup. 
Let $V$ be an irreducible representation of $G$, with highest weight 
$\omega$.
Recall that $\omega$ is {\it minuscule} if $|\langle\omega,
\check{\alpha}\rangle | \le 1$ for any root 
$\alpha$ (see \cite{bou}, VI,1,exercice 24). 
This implies that $\omega$ is a fundamental weight, and all
the weights of $V$ are in the orbit of $\omega$ under the Weyl group $W$. 
Correspondingly, the closed $G$-orbit $G/P\subset\PP V $ is a 
{\it minuscule homogeneous variety}. Here and in the sequel $P$ 
is a maximal parabolic subgroup of $G$.

A related notion is the following: a fundamental weight $\omega$
is {\it cominuscule} if $\langle\omega,\check{\alpha_0}\rangle 
=1$, where $\alpha_0$ denotes the highest root. The two definitions
characterize the same fundamental weights if $G$ is simply laced. 
Equivalently, the {\it cominuscule homogeneous variety} 
$G/P \subset\p V$ is a $G$-Hermitian symmetric space in its minimal 
homogeneous embedding. (Beware that these are called {\it minuscule 
homogeneous varieties} in \cite{landsberg}.)

The list of cominuscule homogeneous varieties is 
the following. There are infinite series of classical homogeneous 
spaces plus two exceptional cases.

$$\begin{array}{ccccc}
type & variety & diagram  & dimension & \hspace*{5mm}index \hspace*{5mm}\\
 A_{n-1} &\hspace{5mm} \G(k,n)\hspace{5mm} &\setlength{\unitlength}{2.5mm}
\begin{picture}(15,3)(-2,0)
\put(0,0){$\circ$}
\multiput(2,0)(2,0){5}{$\circ$}
\multiput(.73,.4)(2,0){5}{\line(1,0){1.34}}
\put(4,0){$\bullet$}
\end{picture} & k(n-k) & n \\
C_n & \hspace{5mm}\G_{\omega}(n,2n)\hspace{5mm} &
\setlength{\unitlength}{2.5mm}
\begin{picture}(15,3)(-2,0)
\put(0,0){$\circ$}
\multiput(2,0)(2,0){5}{$\circ$}
\multiput(.73,.4)(2,0){4}{\line(1,0){1.34}}
\multiput(8.73,.2)(0,.4){2}{\line(1,0){1.34}}
\put(10,0){$\bullet$}
\end{picture} & \frac{n(n+1)}{2} & n+1 \\
D_n & \hspace{5mm}\G_Q(n,2n)\hspace{5mm} &
\setlength{\unitlength}{2.5mm}
\begin{picture}(15,3)(-2,0)
\put(2,0){$\circ$}
\multiput(4,0)(2,0){4}{$\circ$}
\multiput(2.73,.4)(2,0){4}{\line(1,0){1.34}}
\put(0,1.2){$\bullet$}
\put(0,-1.1){$\circ$}
\put(.6,1.5){\line(5,-3){1.5}}
\put(.6,-.64){\line(5,3){1.5}}
\end{picture}
\vspace{.2cm}
 & \frac{n(n-1)}{2} & 2n-2 \\
B_n,D_n & \hspace{5mm}\QQ^m\hspace{5mm} &
 & m & m \\
E_6 &\hspace{5mm} \OO\PP^2 \hspace{5mm}& 
\setlength{\unitlength}{2.5mm}
\begin{picture}(15,3)(-3,-.5)
\put(0,0){$\circ$}
\multiput(2,0)(2,0){4}{$\circ$}
\multiput(.73,.4)(2,0){4}{\line(1,0){1.34}}
\put(0,0){$\bullet$}
\put(4,-2){$\circ$}
\put(4.42,-1.28){\line(0,1){1.36}}
\end{picture}
& 16 & 12 \\
E_7 & \hspace{5mm}E_7/P_7\hspace{5mm} & 
\setlength{\unitlength}{2.5mm}
\begin{picture}(15,3)(-2,0)
\put(0,0){$\circ$}
\multiput(2,0)(2,0){5}{$\circ$}
\multiput(.73,.4)(2,0){5}{\line(1,0){1.34}}
\put(10,0){$\bullet$}
\put(4,-2){$\circ$}
\put(4.42,-1.28){\line(0,1){1.36}}
\end{picture}
& 27 & 18 
\end{array}
$$
\bigskip

\bigskip
Recall that a simply-connected simple complex algebraic group $G$
is determined by its Dynkin diagram, and that the conjugacy classes
of its maximal parabolic subgroups $P$ correspond to the nodes of this
diagram. In the above array, the marked diagrams therefore represent
such conjugacy classes, or equivalently the homogeneous varieties 
$X=G/P$. These are projective varieties, described in the second column,
where the notation $\G(k,n)$ (resp. $\G_\omega(n,2n),\G_Q(n,2n)$)
stands for the Grassmannian of $k$-planes in a fixed $n$-dimensional
space (resp. the Grassmannian of maximal isotropic subspaces in a fixed
symplectic or quadratic vector space of dimension $2n$). To be more 
precise, in the orthogonal case $\G_Q(n,2n)$ will be only one of the 
two (isomorphic) connected component of this Grassmannian; 
its minimal embedding is in the projectivization of a half-spin
representation, and for this reason we sometimes call it a spinor
variety; another consequence is that its index is twice as large as
what one could expect. An
$m$-dimensional quadric is denoted $\QQ^m$ and, finally, the Cayley
plane $\OO \PP^2=E_6/P_1$ and the Freudenthal variety $E_7/P_7$ 
are varieties described in subsections
\ref{subsubsection_cayley} and \ref{subsubsection_freudenthal}.

\smallskip
Beware that a  quadric $\QQ^m$ is minuscule if and only if its dimension $m$
is even. 
Moreover there are two series of $G$-homogeneous varieties which are
cominuscule but not minuscule, but actually their automorphism groups
$H$ are bigger than $G$ and they are minuscule when considered as 
$H$-varieties: namely, the projective spaces of odd dimensions acted
upon by the symplectic groups, and the maximal orthogonal
Grassmannians for orthogonal groups in odd dimensions.

\smallskip      
For a maximal parabolic subgroup $P$ of $G$, the 
Picard group of $G/P$ is free of rank one : $Pic(G/P)=\ZZ H$, where
the very ample generator $H$ defines the minimal homogeneous 
embedding  $G/P\subset \PP V$. 

The index $c_1(G/P)$ is defined by the relation $-K_{G/P}=c_1(G/P)H$. 
This integer is of special importance with respect to quantum
cohomology since it defines the 
degree of the quantum parameter $q$. 
A combinatorial recipe that allows to compute the index of any 
rational homogeneous space can be found in \cite{snow}. Geometrically, 
the index is related with the dimension of the Fano variety $F$
of projective lines on $G/P\subset \PP V$ \cite[p.50]{fulton}:
$$\dim(F)=\dim(G/P)+c_1(G/P)-3.$$
Equivalently, the Fano variety $F_o$ of lines through the base point
$o\in G/P$ has dimension 
$$\dim(F_o)=c_1(G/P)-2.$$
It was proved in \cite[Theorem 4.8]{landsberg}
that when $P$ is defined by a long simple 
root, $F_o$ is homogeneous under the semi-simple part $S$ of $P$. Moreover, the 
weighted Dynkin diagram of $F_o$ can be obtained from that of $G/P$
by suppressing the marked node and marking the nodes that were
connected to it. For example, if $G/P=E_7/P_7$ then $F_o$ is a copy 
of the Cayley plane $E_6/P_6\simeq E_6/P_1$, whose dimension is $16$
-- hence $c_1(E_7/P_7)=18$. 

We also observe  that when $G/P$ is minuscule, the relation
$$c_1(G/P)=e_{max}(G)+1=h$$
does hold, where  $e_{max}(G)$ denotes the maximal exponent of the
Weyl group of $G$ and $h$ is the Coxeter number 
\cite[V, 6, definition 2]{bou}. 

\subsectionplus{Chow rings}

In this section we recall some fundamental facts about the Chow ring 
of a complex rational homogeneous space $X=G/P$. 
This graded Chow ring with coefficients in a ring
$k$ will be denoted $A^*(G/P)_k$, and we simply write $A^*(G/P)$ for
$A^*(G/P)_\Z$ (which coincides with the usual cohomology ring).

First, we recall the {\it Borel presentation} of this Chow ring with
rational coefficients.
Let $W$ (resp. $W_P$) be the Weyl group of $G$ (resp. of $P$). 
Let $\cP$ denote the weight lattice of $G$. The Weyl group $W$
acts on $\cP$. We have
$$A^*(G/P)_{\QQ}\simeq \QQ [\cP]^{W_P}/\QQ [\cP]^W_+, $$
where $\QQ [\cP]^{W_P}$ denotes the ring of 
$W_P$-invariants polynomials on the 
weight lattice, and $\QQ [\cP]^W_+$ is the 
ideal of $\QQ [\cP]^{W_P}$ generated 
by $W$-invariants without constant term 
(see \cite{borel}, Proposition 27.3 or 
\cite{bgg}, Theorem 5.5). 

Recall that the full invariant algebra $\QQ [\cP]^W$ is a polynomial
algebra $\QQ[F_{e_1+1},\ldots ,F_{e_{max}+1}]$, where $e_1,\ldots
,e_{max}$ is the set $E(G)$ of exponents of $G$.  
If $d_1,\ldots ,d_{max}$ denote the
exponents of $S$, we get that $\QQ [\cP]^{W_P}=\QQ[H,I_{d_1+1},\ldots
,I_{d_{max}+1}]$, where $H$ represents the fundamental weight
$\omega_P$ defining
$P$; we denote it this way since geometrically, it corresponds to the 
hyperplane class. Now each $W$-invariant $F_{e_i+1}$ must be
interpreted as a polynomial relation between the $W_P$-invariants
$H,I_{d_1+1},\ldots ,I_{d_{max}+1}$. In particular, if $e_i$ is also 
an exponent of the semi-simple part $S$ of $P$, 
this relation allows to eliminate $I_{e_i+1}$. 
We thus get the presentation, by generators and relations,
\begin{eqnarray}\label{presentation}
A^*(G/P)_{\QQ}\simeq \QQ[H,I_{p_1+1},\ldots ,I_{p_n+1}]/
(R_{q_1+1},\ldots ,R_{q_r+1}),
\end{eqnarray}
where $\{p_1,\ldots ,p_n\}=E(S)-E(G)$ and $\{q_1,\ldots
,q_r\}=E(G)-E(S)$. Note that $q_r=e_{max}$. 

\smallskip
Over the integers, the Chow ring is a free $\ZZ$-module admitting for
basis the classes of the Schubert varieties, the $B$-orbit 
closures in $X$. The Schubert subvarieties of $X$ are parametrised by 
the quotient $W/W_P$. Recall that in any class, 
there exists a unique element of minimal length. 
We denote this set of representatives of  $W/W_P$ by
$W_X$. For $w\in W_X$, let $X(w)$ be the
corresponding Schubert subvariety of $X$. In particular, we let
$w_X$ the unique element of maximal length in $W_X$, such that $X(w_X)=X$.
When $P$ is maximal there is also a unique element of length 
$\ell(w_X)-1$, defining the hyperplane class $H$ with respect to the
minimal homogeneous embedding of $G/P$. All the inclusion relations 
are given by the restriction to $W_X$ of the Bruhat order on $W$. 
The {\it Hasse diagram} of $G/P$ is the graph whose vertices are the 
Schubert classes, and whose edges encode  the inclusion relations in 
codimension one. 

To compute the product of two Schubert classes in the Schubert basis,
we have the following tools:
\begin{enumerate}
\item {\it Poincar{\'e} duality} is known to define an involution of
  $W_X$, given by $w\mapsto w_0ww_0w_X$ (see e.g. \cite{kock}). 
\item {\it the Chevalley formula} allows to multiply any Schubert
  class by the hyperplane class:
$$[X(w)]\cdot H=\sum_{w\rightarrow
  v=ws_{\alpha}}\langle\omega_P,\check{\alpha}\rangle 
[X(v)],$$
where $w\rightarrow v$ must be an arrow in the Hasse diagram
\cite[V, coro 3.2]{hiller}. 
\end{enumerate}
Note that in the minuscule cases the integers
$\langle\omega_P,\check{\alpha}\rangle$
are always equal to one. In particular the degree of each Schubert 
variety can be computed as a number of paths in the Hasse diagram.

\subsectionplus{The exceptional minuscule varieties} 

As we have seen, there exists two exceptional minuscule homogeneous 
spaces, the {\it Cayley plane}  $E_6/P_1\simeq E_6/P_6$, and the 
{\it Freudenthal variety} $E_7/P_7$. We briefly recall how they can
be constructed and a few properties that will be useful in the
sequel. 

\subsubsectionplus{The Cayley plane} 

\label{subsubsection_cayley}

Let ${\bf O}$ denote the normed algebra of (real) octonions, 
and let  $\OO$ be its complexification. The space
$$\JO = \Bigg\{ 
\begin{pmatrix} 
c_1            & x_3            & \bar{x}_2  \\ 
\bar{x}_3 & c_2            &  x_1            \\ 
x_2            & \bar{x}_1 &  c_3            \\ 
\end{pmatrix}
: c_i \in {\bf C}, x_i \in \OO \Bigg\} \cong {\bf C}^{27}
$$
of $\OO$-Hermitian matrices of order $3$, is the exceptional 
simple complex Jordan algebra.

The subgroup of $GL(\JO)$ consisting of automorphisms 
preserving a certain cubic form called determinant is the
simply-connected group of type $E_6$. 
The Jordan algebra $\JO$ and its
dual are the minimal representations of this group. 

The Cayley plane can be defined as the closed $E_6$-orbit in
$\PP\JO$. It is usually denoted $\OO\PP^2$ and indeed it can be 
interpreted as a projective plane over the octonions. In particular
it is covered by a family of eight dimensional quadrics, which are
interpreted as $\OO$-lines.

In fact there exists only two other $E_6$-orbits in $\PP\JO$. The open 
orbit is the complement of the determinantal cubic hypersurface $Det$. 
The remaining one is the complement of $\OO\PP^2$ in $Det$, which can 
be interpreted as the projectivization of the set of rank two matrices
in $\JO$. Observe that geometrically, $Det$ is just the secant variety
of the Cayley plane. For a general point $z\in Det$, the closure of
the union of the 
secant lines to $\OO\PP^2$ passing through $z$ is a linear space
$\Sigma_z$, and $Q_z := \Sigma_z \cap X$ is an $\OO$-line
(see \cite{lv,zak} for more on this).  We deduce the following statement:

\begin{lemm}\label{e6p1coniques} 
Let $q$ be a conic on $\OO\PP^2$, whose supporting plane is not contained 
in $\OO\PP^2$. Then $q$ is contained in a unique $\OO$-line. 
\end{lemm}

\begin{proo} 
Let $P_q$ denote the supporting plane of $q$. Choose any 
$z\in P_q-q$. Then it follows from the definition of $Q_z$
that $q\subset Q_z$. Moreover, two different $\OO$-lines meet
along a linear subspace of $\OO\PP^2$
\cite[proposition 3.2]{zak}. Since  $P_q$ is not contained in $\OO\PP^2$, 
there is no other $\OO$-line containing $q$.  
\end{proo}

The Chow ring of the Cayley plane was computed in \cite{iliev}. 
The Hasse diagram of the Schubert varieties is the following,
classes are indexed by their degree, that is, the codimension of 
the corresponding Schubert variety: 

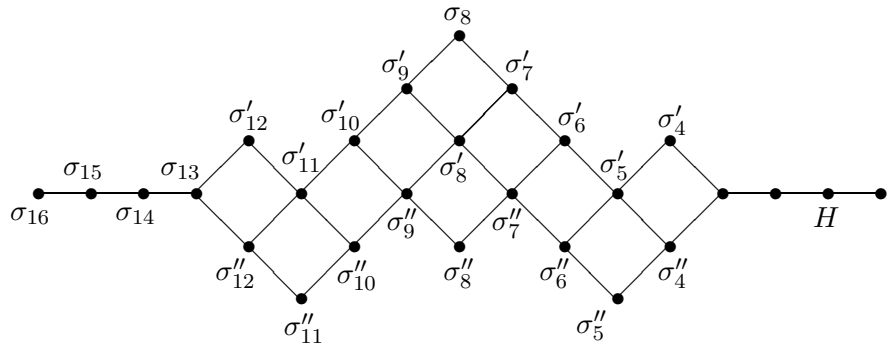
\begin{center}
\setlength{\unitlength}{3.5mm}
\begin{picture}(40,13)(-9,2)

\multiput(-6.2,7.7)(2,0){4}{$\bullet$}
\multiput(1.8,9.7)(4,0){5}{$\bullet$}
\multiput(1.8,5.7)(4,0){5}{$\bullet$}
\multiput(3.8,7.7)(4,0){5}{$\bullet$}
\multiput(3.8,3.7)(12,0){2}{$\bullet$}
\multiput(7.8,11.7)(4,0){2}{$\bullet$}
\multiput(21.8,7.7)(2,0){3}{$\bullet$}
\put(9.8,13.7){$\bullet$}
\put(.8,4.7){$\sigma_{12}''$}
\put(1.3,10.7){$\sigma_{12}'$}
\put(-1.3,8.7){$\sigma_{13}$}
\put(-3,7){$\sigma_{14}$}
\put(-5,8.7){$\sigma_{15}$}
\put(-7,7){$\sigma_{16}$}
\put(9.45,4.7){$\s_8''$}
\put(9.3,8.7){$\s_8'$}
\put(9.45,14.6){$\s_8$}
\put(23.5,6.7){$H$}
\put(17.6,4.7){$\s_4''$}
\put(17.6,10.5){$\s_4'$}
\put(14.5,2.7){$\s_5''$}
\put(15.4,8.8){$\s_5'$}
\put(13.1,4.7){$\s_6''$}
\put(13.8,10.6){$\s_6'$}
\put(11.3,6.6){$\s_7''$}
\put(11.8,12.6){$\s_7'$}
\put(7.3,6.6){$\s_9''$}
\put(7.1,12.6){$\s_9'$}
\put(5.4,4.6){$\s_{10}''$}
\put(4.8,10.7){$\s_{10}'$}
\put(3.4,2.6){$\s_{11}''$}
\put(3.3,9.1){$\s_{11}'$}

\put(10,10){\line(1,1){2}}

\put(10,14){\line(1,-1){8}}
\put(20,8){\line(1,0){6}}
\put(-6,8){\line(1,0){2}}
\put(18,10){\line(1,-1){2}}
\put(14,6){\line(1,1){4}}
\put(16,4){\line(1,1){4}}
\put(12,8){\line(1,1){2}}

\put(-4,8){\line(1,0){4}}
\put(0,8){\line(1,1){2}}
\put(2,6){\line(1,1){8}}
\put(4,4){\line(1,1){8}}
\put(10,6){\line(1,1){2}}

\put(0,8){\line(1,-1){4}}
\put(2,10){\line(1,-1){4}}
\put(6,10){\line(1,-1){4}}
\put(8,12){\line(1,-1){8}}

\end{picture}
\end{center}

Among the Schubert classes some are particularly remarkable. The class
of an $\OO$-line is $\s_8$. Moreover the Cayley plane has two families 
of maximal linear spaces, some of dimension five whose class is
$\s''_{11}$, and some of dimension four  whose class is $\s_{12}'$. 
Note also that Poincar{\'e}
duality is given by the obvious  symmetry of the Hasse diagram. 

Let $\s'_4$ be the Schubert class which is Poincar{\'e} dual to
$\s_{12}'$. The Hasse diagram shows that the Chow ring 
$A^*(\OO\PP^2)$ is generated by $H$, 
$\s'_4$ and $\s_8$ (proposition 
\ref{e6_classique} will give a stronger result), 
and the multiplication table of its Schubert cells is completely 
determined by the Chevalley formula, Poincar{\'e} duality and the two
formulas
$$(\s'_4)^2=\s_8+\s'_8+\s''_8, \qquad \s'_4\s_8=\s_{12}'.$$

An easy consequence is a presentation of $A^*(\OO \PP^2)$ as a quotient
of a polynomial algebra over $\Z$:

\begin{prop}
\label{e6_classique}
Let ${\cal H} = \ZZ[h,s]/(3hs^2-6h^5s+2h^9,s^3-12h^8s+5h^{12})$.
Mapping $h$ to $H$ and $s$ to $\sigma'_4$
yields an isomorphism of graded algebras $${\cal H} \simeq
A^*(\OO\PP^2).$$
\end{prop}

\begin{proo}
Let $\sigma=\sigma'_4$. 
Using the Chevalley formula, we get successively that
$$
\begin{array}{rcl}
\sigma''_4 & = & H^4-\sigma\\
\sigma'_5  & = & H\sigma\\
\sigma''_5 & = & -2H\sigma + H^5\\
\sigma''_6 & = & -2H^2\sigma + H^6\\
\sigma'_6  & = & 3 H^2 \sigma -H^6\\
\sigma''_7 & = & -2 H^3\sigma + H^7\\
\sigma'_7  & = & 5H^3\sigma -2H^7\\
\sigma_8   & = & \sigma^2 + 2H^4\sigma -H^8\\
\sigma'_8  & = & -\sigma^2 + 3 H^4 \sigma -H^8\\
\sigma''_8 & = & \sigma^2 - 5 H^4\sigma + 2H^8
\end{array}
$$
The last three identities have been found taking into account the fact
$\sigma^2 = \sigma_8 + \sigma'_8 + \sigma''_8$ \cite[(16), p.11]{iliev}.
Computing a hyperplane section of these cells, we get the relation
$(\sigma_8-\sigma'_8+\sigma''_8)H = 0$, namely
\begin{equation}
\label{classique9}
3H\sigma^2 - 6 H^5\sigma + 2 H^9 = 0.
\end{equation}
In the following, we compute over the rational numbers and
use this relation to get rid of the terms involving
$\sigma^2$. We get
$$
\begin{array}{rcl}
\sigma'_9 & = & 4H^5\sigma - 5/3 H^9\\
\sigma''_9 & = & -3H^5\sigma + 4/3 H^9\\
\sigma'_{10} & = & 4H^6\sigma -5/3 H^{10}\\
\sigma''_{10} & = & -7H^6\sigma + 3H^{10}\\
\sigma'_{11} & = & 4H^7\sigma - 5/3 H^{11}\\
\sigma''_{11} & = & -11 H^7\sigma +14/3 H^{11}\\
\sigma''_{12} & = & -11H^8\sigma + 14/3 H^{12}\\
\sigma'_{12} & = & 15 H^8\sigma - 19/3 H^{12}
\end{array}
$$
Using the relation
$\sigma'_4.\sigma_8 = \sigma'_{12}$ \cite[proposition 5.2]{iliev}, we therefore
get the relation in degree 12:
\begin{equation}
\label{classique12}
\sigma^3-12H^8\sigma + 5H^{12} = 0.
\end{equation}
This implies $26H^9\sigma = 11 H^{13}$, which we use to eliminate $\sigma$
and get $H^i=78\sigma_i$ for $13 \leq i \leq 16$.

Therefore, we see that there is indeed a morphism of algebras
$f:{\cal H} \rightarrow A^*(\OO\PP^2)$ mapping $h$ to $H$ and 
$s$ to $\sigma'_4$.
This map is surjective up to degree 8, by the previous
computation of the Schubert classes. In the remaining degrees, one can
similarly express the Schubert cells as integer polynomials in $H$ and
$\sigma'_4$, proving that $f$ is
surjective. We claim that $\cH$ is a free $\Z$-module of rank
27. Therefore, $f$ must also be injective and it is an isomorphism.

\smallskip
Let us check that $\cH$ is indeed a free $\Z$-module of rank 27. Let
$\Z[h,s]_d \subset \Z[h,s]$ (resp. $\cH_d \subset \cH$) denote the
degree-$d$ part (where of course $h$ has degree 1 and $s$ has degree
4). Since there are no relations in degree $d$ for $0 \leq d \leq 8$,
$\cH_d$ is a free $\Z$-module. For $9 \leq d \leq 12$, since
$\cH_d = ( \Z.h^{d-8}s^2 \oplus \Z.h^{d-4}s \oplus \Z.h^d) /
(3h^{d-8}s^2-6h^{d-4}s+2h^d)$ and $3,6,2$ are coprime, $\cH_d$ is
free.

Let $13 \leq d \leq 16$.
Eliminating $h^{d-12}s^3$, it comes that 
$\cH_d = ( \Z.h^{d-8}s^2 \oplus \Z.h^{d-4}s \oplus \Z.h^d) /
(3h^{d-8}s^2 - 6h^{d-4}s + 2h^d, 26 h^{d-4}s - 11 h^d)$.
Since the $(2\times 2)$-minors of the matrix
$$
\left (
\begin{array}{ccc}
3 & -6 & 2\\
0 & 26 & -11
\end{array}
\right )
$$
are 78, 33, 118, and are therefore coprime, the $\Z$-module
$\cH_d$ is again free.

Finally, in degree 17, the relations are
$$
\left \{
\begin{array}{cccccl}
&& 3h^9s^2 & -6h^{13}s & +2 h^{17} & =\ \ 0\\
& 3h^5s^3 & -6 h^9s^2 & +2 h^{13}s & & =\ \ 0\\
3hs^4 & -6h^5s^3 & +2 h^9s^2 &&& =\ \ 0\\
& h^5s^3 & -12 h^{13}s & +5h^{17} & & =\ \ 0\\
hs^4 & -12 h^9s^2 & +5 h^{13}s & & & =\ \ 0\ \ \ \ .
\end{array} 
\right .
$$
Since the determinant of the matrix
$$
\left (
\begin{array}{ccccc}
0 & 0 & 3 & -6 & 2\\
0 & 3 & -6 & 2 & 0\\
3 & -6 & 2 & 0 & 0\\
0 & 1 & 0 & -12 & 5\\
1 & 0 & -12 & 5 & 0\\
\end{array}
\right )
$$
is 1, $\cH_{17} = 0$ and we are done.
\end{proo}

\subsubsectionplus{The Freudenthal variety} 

\label{subsubsection_freudenthal}

The other exceptional minuscule homogeneous variety can 
be interpreted as the 
{\it twisted cubic over the exceptional Jordan algebra}. 
Consider the {\it Zorn algebra} 
$$\cZ_2(\OO)=\CC\oplus\JO\oplus\JO\oplus\CC.$$
One can prove that there exists an
action of $E_7$ on that $56$-dimensional vector space 
(see \cite{freudenthal}). 
Then the closed $E_7$-orbit inside $\PP \cZ_2(\OO)$ is the {\it 
Freudenthal variety} $E_7/P_7$. It was studied extensively by 
Freudenthal, and more recently in \cite{kaji} through a slightly 
different point of view. It can be constructed explicitely as the 
closure of the set of elements of the form $[1,X,\com(X),\det(X)]$ 
in $\PP \cZ_2(\OO)$, where $X$ belongs to $\JO$ and its comatrix
$\com(X)$ is defined by the usual formula for order three matrices,
so that $X\com(X)=\det(X)I$.  

For future use we notice the following two geometric properties of the 
low-degree rational curves in the Freudenthal variety:

\begin{lemm}\label{e7p7coniques}
A general conic  on the Freudenthal variety
is contained in a unique maximal quadric. 
\end{lemm}

\begin{proo} Let $q$ be such a conic. We may suppose that $q$ passes 
through the point $[1,0,0,0]$, and we let $[0,Y,0,0]$ be the tangent 
direction to $q$ at $[1,0,0,0]$. A general point in $q$ 
is of the form $[1,X,\com(X),\det(X)]$. 

Now the supporting plane of $q$ is the closure of the set of points of
the form $[1,sX+tY,s\com(X),s\det(X)]$. To belong to the Freudenthal
variety, such a point must verify  the condition
$(sX+tY)s\com(X)=s\det(X)I$, hence  $tY\com(X)=(1-s)\det(X)I$ if $s\ne
0$. But $q$ being smooth cannot verify any linear condition, so we
must have $\det(X)=0$ and  $Y\com(X)=0$. This means that $X$
has rank at most two, in fact exactly two by the genericity of $q$.
So up to the action of $E_6$ we may suppose that 
$$X=\begin{pmatrix} 1&0 &0\\0&1 &0\\0&0 &0\end{pmatrix},
\qquad \com(X)=\begin{pmatrix} 0&0 &0\\0&0 &0\\0&0 &1\end{pmatrix},
\qquad Y=\begin{pmatrix} *&* &0\\ *&* &0\\0&0 &0\end{pmatrix}.$$
Note that $Z=\com(X)$ is a rank one matrix, hence defines a point of
$\OO\PP^2$, and since $YZ=0$, $Y$ has to be contained in the orthogonal 
(with respect to the quadratic form $M\mapsto\trace(M^2)$) $\Sigma_Z$
to the 
affine tangent space of $\OO\PP^2$ at $[Z]$.  Finally, $q$ must be
contained in 
$\CC\oplus \Sigma_Z\oplus \CC Z\simeq \CC^{12}$, whose intersection
with the Freudenthal variety is a maximal ten dimensional quadric
uniquely defined by $q$. This proves the claim. 
\end{proo}

\begin{lemm}\label{e7p7cubiques}
Through three general points of the Freudenthal variety, there 
is a unique twisted cubic curve. 
\end{lemm}

\begin{proo}
Let $p_0$, $p_{\infty}$ be two of the three points. Because of the
Bruhat decomposition, we may suppose
that $p_0$ defines a highest weight line, and $p_{\infty}$ a lowest 
weight line  in $\cZ_2(\OO)$. Then their common stabilizer in 
$E_7$ contains a copy of $E_6$, and the restriction of the
$E_7$-module $\cZ_2(\OO)$ to $E_6$ decomposes as 
$\CC\oplus\JO\oplus\JO\oplus\CC$. Otherwise said, we may suppose
that in this decomposition, $p_0=[1,0,0,0]$ and $p_{\infty} =[0,0,0,1]$. 

Now our third generic point is of the form 
$[1,M,\com(M),\det(M)]$ with $\det(M) \not = 0$.
The twisted cubic
$C_0 := \{ [t^3,t^2uM,tu^2\com(M),u^3\det(M)] : [t,u] \in \p^1 \}$
obviously passes through $p_0,p_1,p_{\infty}$; let us prove that it is
the only such curve.  

For simplicity we denote the Freudenthal variety by $X$. 
Let $C$ be a rational curve of degree three 
passing through $p_0,p_1,p_{\infty}$ and let us proove it must be $C_0$. 
First we notice that $C$ must be irreducible, since otherwise there
would be  a line or a conic
through $p_0$ and $p_{\infty}$, yielding a contradiction with
$\widehat{T_{p_0}X} \cap \widehat{T_{p_{\infty}}X} = \{0\}$. 

So let $T_i := \widehat {T_{p_i}C} \subset \cZ_2(\OO)$ for
$i \in \{0,\infty\}$. Since $C$ is a twisted cubic, the linear span $S$
of $\widehat C$ in $\cZ_2(\OO)$ satisfies $S = T_0 \oplus T_\infty$. So
$p_1=[1,M,\com(M),\det(M)] \in T_0 \oplus T_\infty$. Now, $T_0$ is included in
$(*,*,0,0)$ and contains $(1,0,0,0)$, and similarly for $T_\infty$. It
follows that $T_0 = \{ (\lambda,\mu M,0,0) : \lambda,\mu \in \C \}$
and $T_\infty = \{ (0,0,\lambda \com(M),\mu) : \lambda,\mu \in \C \}$.
Therefore, $S = \{ (\lambda,\mu M,\nu \com(M),\omega) :
\lambda,\mu,\nu,\omega \in \C \}$, and since $X \cap \p S = C_0$, we
get $C = C_0$.\end{proo}

The Chow ring of the Freudenthal variety was recently computed in 
\cite{nikolenko} with the help of a computer. The Hasse diagram 
of Schubert classes is the following:

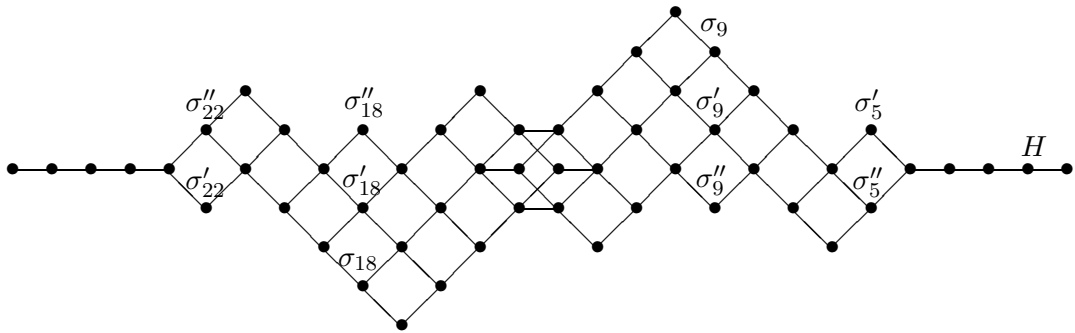
\begin{center}
\setlength{\unitlength}{2.6mm}
\begin{picture}(40,16)(-1,0)

\put(-8,8){\line(1,0){8}}
\put(0,8){\line(1,1){4}}
\put(0,8){\line(1,-1){2}}
\put(2,10){\line(1,-1){10}}
\put(2,6){\line(1,1){4}}
\put(6,6){\line(1,1){4}}
\put(4,12){\line(1,-1){10}}
\put(10,10){\line(1,-1){6}}
\put(16,12){\line(1,-1){4}}
\put(14,10){\line(1,-1){4}}
\put(8,4){\line(1,1){8}}
\put(10,2){\line(1,1){8}}
\put(12,0){\line(1,1){6}}
\put(18,6){\line(1,1){2}}
\put(16,8){\line(1,0){2}}
\put(18,6){\line(1,0){2}}
\put(18,10){\line(1,0){2}}
\put(20,8){\line(1,0){2}}
\put(18,8){\line(1,1){8}}
\put(18,8){\line(1,-1){4}}
\put(20,6){\line(1,1){8}}
\put(22,4){\line(1,1){8}}
\put(20,10){\line(1,-1){4}}
\put(22,12){\line(1,-1){6}}
\put(24,14){\line(1,-1){10}}
\put(26,16){\line(1,-1){10}}
\put(28,6){\line(1,1){4}}
\put(32,6){\line(1,1){4}}
\put(34,4){\line(1,1){4}}
\put(36,10){\line(1,-1){2}}
\put(38,8){\line(1,0){8}}
\multiput(-8.2,7.7)(2,0){5}{$\bullet$}
\multiput(1.7,9.7)(4,0){5}{$\bullet$}
\multiput(1.7,5.7)(4,0){5}{$\bullet$}
\multiput(3.7,7.7)(4,0){5}{$\bullet$}
\multiput(7.7,3.7)(4,0){3}{$\bullet$}
\multiput(3.7,11.7)(12,0){2}{$\bullet$}
\multiput(9.7,1.7)(4,0){2}{$\bullet$}
\put(11.7,-.3){$\bullet$}
\multiput(17.7,7.7)(4,0){5}{$\bullet$}
\multiput(19.7,9.7)(4,0){5}{$\bullet$}
\multiput(19.7,5.7)(4,0){5}{$\bullet$}
\multiput(21.7,11.7)(4,0){3}{$\bullet$}
\multiput(23.7,13.7)(4,0){2}{$\bullet$}
\multiput(25.7,15.7)(4,0){1}{$\bullet$}
\multiput(21.7,3.7)(12,0){2}{$\bullet$}
\multiput(37.7,7.7)(2,0){5}{$\bullet$}
\put(43.7,8.7){$H$}
\put(35.2,11){$\sigma'_5$}
\put(35.1,7){$\sigma''_5$}
\put(27.3,15.1){$\sigma_9$}
\put(27.1,11.1){$\sigma'_9$}
\put(27.1,7.1){$\sigma''_9$}
\put(1,7){$\sigma'_{22}$}
\put(1,11){$\sigma''_{22}$}
\put(8.8,3.1){$\sigma_{18}$}
\put(9.1,11.1){$\sigma''_{18}$}
\put(9,7.1){$\sigma'_{18}$}
\end{picture}
\end{center}

We did not indicate all the Schubert cells in this diagram, but we use
the same convention as for the Cayley plane: when there are several
Schubert cells of the same codimension $c$, we denote them 
$\sigma''_c,\sigma'_c,\sigma_c$, starting with the lowest cell
on the diagram up to degree $13$, and with the highest one from
degree $14$. Poincar{\'e} duality is given by the obvious central 
symmetry of this diagram. 

As we will see, the Chow ring is generated by the hyperplane section $H$
and two Schubert classes of degree $5$ and $9$ (in accordance with
(\ref{presentation})). We choose the Schubert 
classes $\s'_5$ and $\s_9$. The results of Nikolenko and Semenov 
can be stated as follows:
\begin{eqnarray} \nonumber
(\s'_5)^2 &=& 2\s'_9H, \\ \nonumber
\s'_5\s''_5 &=& (\s_9+\s'_9+2\s''_9)H, \\ \nonumber
(\s''_5)^2 &=& (3\s'_9+\s''_9)H, \\ \nonumber
\s_9\s'_5 &=& \s_{14}+2\s'_{14}+2\s''_{14}, \\ \nonumber
\s'_9\s'_5 &=& 3\s_{14}+4\s'_{14}+4\s''_{14}, \\ \nonumber
\s''_9\s'_5 &=& 2\s_{14}+3\s'_{14}+2\s''_{14}, \\ \nonumber
\s_9\s''_5 &=& \s_{14}+3\s'_{14}+3\s''_{14}, \\ \nonumber
\s'_9\s''_5 &=& 4\s_{14}+6\s'_{14}+5\s''_{14}, \\ \nonumber
\s''_9\s''_5 &=& 3\s_{14}+3\s'_{14}+3\s''_{14}, \\ \nonumber
(\s_9)^2 &=& 2\s_{18}+2\s'_{18}, \\ \nonumber
(\s'_9)^2 &=& 4\s_{18}+10\s'_{18}+6\s''_{18}, \\ \nonumber
(\s''_9)^2 &=& 2\s_{18}+4\s'_{18}+2\s''_{18}, \\ \nonumber
\s_9\s'_9 &=& 2\s_{18}+4\s'_{18}+3\s''_{18}, \\ \nonumber
\s'_9\s''_9 &=& 3\s_{18}+6\s'_{18}+4\s''_{18}, \\ \nonumber
\s_9\s''_9 &=& 3\s'_{18}+2\s''_{18}.
\end{eqnarray}

These formulas, plus the Chevalley formula and Poincar{\'e} duality, 
completely determine the multiplication table of Schubert classes in 
$E_7/P_7$. As for the case of the Cayley plane, we deduce a
presentation of $A^*(E_7/P_7)$.

\begin{theo}
\label{e7_classique}
Let ${\cal H} = \ZZ[h,s,t]/(s^2 - 10sh^5 + 2th + 4h^{10},
2st - 12sh^9 + 2th^5 + 5h^{14}, t^2 + 922 sh^{13} - 198 th^9 - 385 h^{18})$.
Mapping $h$ to $H$, $s$ to $\sigma'_5$ and $t$ to $\sigma_9$
yields an isomorphism of graded algebras $${\cal H} \simeq
A^*(E_7/P_7).$$
\end{theo}

\begin{proo}
The fact that the displayed relations are relations in the Chow ring
follows from the previous formulas and an expression of the Schubert
cells as polynomials in the generators similar to that we did for the
proof of proposition \ref{e6_classique}.

We therefore have a morphism
$f:\cH \rightarrow A^*(E_7/P_7)$ mapping
$h$ to $H$, $s$ to $\sigma'_5$ and $t$ to $\sigma_9$; 
the fact that $f$ is surjective can be read on the Hasse
diagram as for proposition \ref{e6_classique}, except in degree 14.
But in this degree, we note that
$\sigma_{14} + \sigma'_{14} = \sigma_{13}.H,
\sigma_{14} + \sigma''_{14} = \sigma'_{13}.H$
and $\sigma'_{14} +\sigma''_{14} = \sigma''_{13}.H$ belong to the
image of $f$. This does not imply that $A^{14}(E_7/P_7) \subset \im(f)$, 
but the surjectivity of $f$ follows from
the equality $\sigma'_9.\sigma'_5 = \sigma_{14} + 2 \sigma'_{14}
+ 2 \sigma''_{14}$.

The injectivity of $f$ will again follow from the fact that $\cH$ is a
free module of rank 56. This is a lengthy computation; we only give some
indications to the reader for the relevant degrees. Our strategy is that
we use the first and the third relation to get rid of monomials involving
$s^2$ or $t^2$, and then use Gauss elimination to prove that $\cH$ is
free. We denote $\cH_i$ the component of $\cH$ of degree $i$.
So for example, in degree 14, $\cH_{14}$ is generated as a module by
$h^{14},h^9s,h^5t,st$, which satisfy the relation
$2st-12h^9s+2h^5t+5h^{14} = 0$. Since $2,-12,2,5$ are comprime, $\cH_{14}$
is a free $\Z$-module. For $\cH_{19}$, the second relation, multiplied by
$h^5$ and $s$, gives
$$
\left \{
\begin{array}{c}
2h^5st-12h^{14}s+2h^{10}t+5h^{19} = 0\\
22h^5st+3573h^{14}s-776h^{10}t-1492h^{19} = 0.
\end{array}
\right .
$$
After Gauss elimination, we find that this is equivalent to
$$
\left \{
\begin{array}{l}
3094h^5st-39h^{14}s-8962h^{10}t = 0\\
1238h^5st-18h^{14}s-358h^{10}t+h^{19} = 0,
\end{array}
\right .
$$
therefore $\cH_{19}$ is again a free module.
Similaly, in degree 23, the relations are~:
$$
\left \{
\begin{array}{l}
1238h^9st-18h^{18}s-358h^{14}t+h^{23} = 0\\
-1312h^9st+52h^{18}s -h^{14}t = 0 \\
5586h^9st-221h^{18}s = 0.
\end{array}
\right .
$$
Finally, in degree 28, we multiply the second relation by
$h^{14},h^9s,h^5t,st$, and get~:
$$
\left \{
\begin{array}{l}
2h^{14}st-12h^{23}s+2h^{19}t+5h^{28} = 0 \\
22h^{14}st-776h^{19}t+3573h^{23}s-1492h^{28} = 0 \\
384h^{14}st+4089h^{19}t-19514h^{23}s+8146h^{28} = 0 \\
7929h^{14}st-114572h^{19}t+521102h^{23}s-217624h^{28} = 0.
\end{array}
\right .
$$
Since the determinant of this system is 1, $\cH_{28}$ is the trivial
$\Z$-module.
\end{proo}

\subsectionplus{Quivers and Poincar{\'e} duality}

\label{subsection_quivers}

The archetypal minuscule homogeneous variety is the Grassmannian,
whose Schubert varieties are indexed by partitions whose Ferrers
diagram are contained in a fixed rectangle. In particular, Poincar{\'e}
duality for Schubert classes is easily visualized: it 
associates to such a partition its complementary partition 
inside the rectangle. 
 
In this section, we argue that there exists a general 
very convenient way to visualize
Poincar{\'e} duality for any minuscule or cominuscule homogeneous variety
$X=G/P$. This was first observed in \cite{carquois}. The main idea 
is to associate to $X$ a quiver $Q_X$, which when $X$ is a
Grassmannian will be the rectangle we just mentionned. 

\smallskip
We start with $X$ any rational homogeneous space and a reduced
expression for $w_X$, the longest element 
in $W_X$ -- say $w_X=s_{\b_1}\cdots s_{\b_N}$ where $N=\dim(X)$. 
 An important point here is that in 
the (co)minuscule case,
this reduced expression is unique up to commutation relations (see
\cite{stembridge}).

 \begin{definition}\
\label{definition_carquois}
   \begin{itemize}
     \item
       For $\beta$ a simple root, let $m_X(\beta)$ be the number of occurences
       of $\beta$ in the reduced expression
       $w_X=s_{\b_1}\cdots s_{\b_N}$
       ($m_X(\beta) = \# \{ j : \beta_j = \beta \}$).
     \item
       For $(\beta,i)$ such that $1 \leq i \leq m_X(\beta)$, let
       $r(\beta,i)$ denote the integer $j$ such that $\beta_j = \beta$
       and $\# \{ k \leq j : \beta_k = \beta \} = i$. If $i>m_X(\beta)$,
       set $r(\beta,i) = + \infty$. Set also $r(\beta,0)=0$.
     \item
       Let $Q_X$ be the quiver whose set of vertices is the set of pairs
       $(\beta,i)$, where $\beta$ is a simple root and 
       $1\leq i \leq m_X(\beta)$
       and whose arrows are given as follows. There is an arrow
       from $(\beta,i)$ to $(\gamma,j)$
       if $\scal{\gamma^\vee,\beta}\neq 0$ and 
$r(\gamma,j-1) < r(\beta,i) < r(\gamma,j) < r(\beta,i+1)$.
   \end{itemize}
 \end{definition}
Beware that this definition is slightly different from that given in 
\cite{carquois}, but the corresponding quivers can easily be recovered
one from the other. 
For $X$ (co)minuscule, since the commutation relations do not change 
the quiver (to check this, it is enough to check that the quiver does not
change when one commutes two commuting reflexions), 
$Q_X$ is {\it uniquely defined}. 

\smallskip
Another important property is that it is {\it symmetric}. Indeed,
recall that the Poincar{\'e} duality on $X$ is defined by the 
involution 
$$w\mapsto w_0ww_0w_X, \quad w\in W_X.$$ 
In particular, since $w_X$ defines the fundamental class in $X$,
its Poincar{\'e} dual is the class of a point, so
$w_0w_X$ must be an involution.
Thus $w_X^{-1}=w_0w_Xw_0^{-1}$.
Now, if $w_X=s_{\b_1}\cdots s_{\b_N}$ is any reduced expression, it
follows that $w_X=s_{w_0(\beta_N)} \cdots s_{w_0(\beta_1)}$ is also a
reduced expression. Since
$w_0(\beta_j) = \iota(\beta_j)$, where $\iota$ is the Weyl involution on the
Dynkin diagram, it follows
that the involution 
$(\beta,k) \mapsto i_X(\beta,k) = (\iota(\beta),m_X(\beta)+1-k)$ 
induces an arrow-reversing automorphism of the quiver $Q_X$.

\smallskip
More generally, any $w\in W_X$ is given by a subexpression 
of some reduced expression of $w_X$, and again its reduced 
decomposition is unique up to commutation relations. We can 
therefore define by the same recipe a unique quiver $Q_w$,
which is a subquiver of $Q_{w_X}=Q_X$. 
We think of these quivers 
as combinatorial tools generalizing the (strict) 
partitions parametrizing Schubert subvarieties in
(isotropic)  Grassmannians. 

 \begin{example}
 \label{exemple}
 (\i) For a Grassmannian $X=\G(p,n)$, the quiver
 $Q_X$ is the $p\times(n-p)$-rectangle. The set $W_X$ can be
 identified with the set of partitions contained in this rectangle, 
 and the quiver $Q_w$, for any $w\in W_X$, is the
 complement in the rectangle of the partition defining $X(w)$ (see
 \cite{carquois} and examples \ref{exemple-grass1} and
 \ref{exemple-grass2} below).

 (\i\i)  Let $X=E_6/P_1$ be the Cayley plane.  A symmetric reduced
   expression for $w_X$ is given by 
$$w_X= s_{\alpha_6}s_{\alpha_5}s_{\alpha_4}
   s_{\alpha_2} s_{\alpha_3} s_{\alpha_4} s_{\alpha_5}
   s_{\alpha_1} s_{\alpha_6} s_{\alpha_3}  s_{\alpha_4}  s_{\alpha_5}
   s_{\alpha_2}
   s_{\alpha_4}s_{\alpha_3} s_{\alpha_1}.$$ 
The quiver $Q_X$ is the following (all arrows are going down).

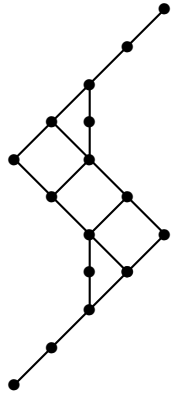
\psset{xunit=0.5cm}
\psset{yunit=0.5cm}
\centerline{
\begin{pspicture*}(-4.3,-4.3)(0.3,7.3)
\psline(-0,0)(-1,1)
\psline(-0,6)(-3,3)
\psline(-2,4)(-2,3)
\psline(-2,2)(-1,1)
\psline(-2,2)(-3,3)
\psline(-1,-1)(-4,2)
\psline(-2,2)(-2,3)
\psline(-3,3)(-4,2)
\psline(-2,2)(-3,1)
\psline(-1,1)(-2,0)
\psline(-0,0)(-4,-4)
\psline(-2,0)(-2,-2)
\put(-0.1,-0.1){$\bullet$}
\put(-0.6,0.4){$\bullet$}
\put(-1.1,0.9){$\bullet$}
\put(-1.6,1.4){$\bullet$}
\put(-1.1,-0.1){$\bullet$}
\put(-1.6,0.4){$\bullet$}
\put(-2.1,0.9){$\bullet$}
\put(-0.6,-0.6){$\bullet$}
\put(-1.1,1.4){$\bullet$}
\put(-1.1,-0.6){$\bullet$}
\put(-2.1,-2.1){$\bullet$}
\put(-1.6,-1.6){$\bullet$}
\put(-1.1,-1.1){$\bullet$}
\put(-0.6,-0.6){$\bullet$}
\put(-0.1,2.9){$\bullet$}
\put(-0.6,2.4){$\bullet$}
\put(-1.1,1.9){$\bullet$}
\end{pspicture*}}
 \end{example}

(\i\i\i) Let $X=E_7/P_7$ be the Freudenthal variety. Then $Q_X$ is the 
following (again with all arrows going down).

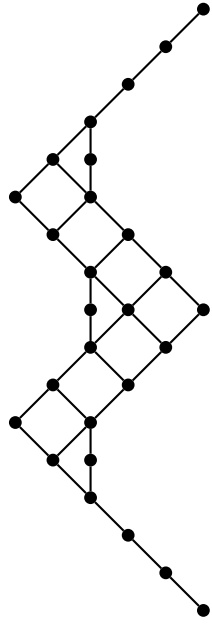
\psset{unit=5mm}

\newcommand{\cercleplein}[2]{
\pscircle*(#1,#2){.17}   }

\newcommand{\cerclevide}[2]{
\pscircle(#1,#2){.17}   }

\newcommand{\cerclepleinrouge}[2]{
\pscircle*[linecolor=red](#1,#2){.17}    }

\newcommand{\cercleviderouge}[2]{
\pscircle[linecolor=red](#1,#2){.17}   }

\newcommand{\cerclepleinbleu}[2]{
\pscircle*[linecolor=blue](#1,#2){.17}   }

\newcommand{\cerclevidebleu}[2]{
\pscircle[linecolor=blue](#1,#2){.17}   }

\centerline{
\begin{pspicture*}(-.5,-.5)(6.5,18.5)
\psline[linecolor=black](4.88,17.88)(4.12,17.12)
\psline[linecolor=black](3.88,16.88)(3.12,16.12)
\psline[linecolor=black](2.88,15.88)(2.12,15.12)
\psline[linecolor=black](1.88,14.88)(1.12,14.12)
\psline[linecolor=black](2,14.83)(2,14.17)
\psline[linecolor=black](0.88,13.88)(0.12,13.12)
\psline[linecolor=black](1.12,13.88)(1.88,13.12)
\psline[linecolor=black](2,13.83)(2,13.17)
\psline[linecolor=black](0.12,12.88)(0.88,12.12)
\psline[linecolor=black](1.88,12.88)(1.12,12.12)
\psline[linecolor=black](2.12,12.88)(2.88,12.12)
\psline[linecolor=black](1.12,11.88)(1.88,11.12)
\psline[linecolor=black](2.88,11.88)(2.12,11.12)
\psline[linecolor=black](3.12,11.88)(3.88,11.12)
\psline[linecolor=black](2,10.83)(2,10.17)
\psline[linecolor=black](2.12,10.88)(2.88,10.12)
\psline[linecolor=black](3.88,10.88)(3.12,10.12)
\psline[linecolor=black](4.12,10.88)(4.88,10.12)
\psline[linecolor=black](2,9.83)(2,9.17)
\psline[linecolor=black](2.88,9.88)(2.12,9.12)
\psline[linecolor=black](3.12,9.88)(3.88,9.12)
\psline[linecolor=black](4.88,9.88)(4.12,9.12)
\psline[linecolor=black](1.88,8.88)(1.12,8.12)
\psline[linecolor=black](2.12,8.88)(2.88,8.12)
\psline[linecolor=black](3.88,8.88)(3.12,8.12)
\psline[linecolor=black](0.88,7.88)(0.12,7.12)
\psline[linecolor=black](1.12,7.88)(1.88,7.12)
\psline[linecolor=black](2.88,7.88)(2.12,7.12)
\psline[linecolor=black](0.12,6.88)(0.88,6.12)
\psline[linecolor=black](1.88,6.88)(1.12,6.12)
\psline[linecolor=black](2,6.83)(2,6.17)
\psline[linecolor=black](1.12,5.88)(1.88,5.12)
\psline[linecolor=black](2,5.83)(2,5.17)
\psline[linecolor=black](2.12,4.88)(2.88,4.12)
\psline[linecolor=black](3.12,3.88)(3.88,3.12)
\psline[linecolor=black](4.12,2.88)(4.88,2.12)
\pscircle*(0,13){.17}
\pscircle*(0,7){.17}
\pscircle*(1,14){.17}
\pscircle*(1,12){.17}
\pscircle*(1,8){.17}
\pscircle*(1,6){.17}
\pscircle*(2,15){.17}
\pscircle*(2,14){.17}
\pscircle*(2,13){.17}
\pscircle*(2,11){.17}
\pscircle*(2,10){.17}
\pscircle*(2,9){.17}
\pscircle*(2,7){.17}
\pscircle*(2,6){.17}
\pscircle*(2,5){.17}
\pscircle*(3,16){.17}
\pscircle*(3,12){.17}
\pscircle*(3,10){.17}
\pscircle*(3,8){.17}
\pscircle*(3,4){.17}
\pscircle*(4,17){.17}
\pscircle*(4,11){.17}
\pscircle*(4,9){.17}
\pscircle*(4,3){.17}
\pscircle*(5,18){.17}
\pscircle*(5,10){.17}
\pscircle*(5,2){.17}
\end{pspicture*}}

\psset{unit=1cm}
\smallskip
For the quivers of the remaining classical (co)minuscule varieties 
we refer to \cite{carquois}. 

To avoid confusion will always draw our quivers vertically as 
above, contrary to Hasse diagrams which we draw horizontally. 
Nevertheless there is a strong connexion between quivers and 
Hasse diagrams, at least in the (co)minuscule cases. since we will
not need it in hte sequel we just state the following result:

\begin{prop}\label{quiver_Hasse}
Let $X=G/P$ be a cominuscule variety with base point $o$. 
Then the quiver $Q_X$ coincides with the Hasse diagram of the 
Fano variety $F_o\subset\PP T_oX$ of lines through $o$. 
\end{prop}

For example, if $X=\G(k,n)$, the Fano variety 
of lines through $o$ is $F_o=\PP^{k-1}\times \PP^{n-k-1}$. 
These two projective spaces can 
be interpreted as the two sides of the rectangular quiver $Q_X$.   
\medskip

For any (co)minuscule homogeneous variety $X$, 
the quiver $Q_{X}$ has a natural partial order given by $i\preccurlyeq j$ if
there exists an oriented path from $j$ to $i$ (see \cite{carquois} for
more details). This induces a partial order on each of the 
subquivers $Q_w$. 

\begin{defi}
\label{def_pic}
The {\it peaks} of $Q_w$ are its maximal
elements. We denote by $p(Q_w)$ the set of peaks of $Q_w$.
\end{defi}

The fact that, for (co)minuscule homogeneous spaces, the Bruhat order
is generated by simple reflections (see for example \cite{LMS})
implies that the subquivers $Q_w$ of $Q_X$ corresponding to the
Schubert subvarieties $X(w)$ are obtained from $Q_{X}$ inductively by
removing peaks. In other words, they are the Schubert subquivers of
$Q_X$, according to the following definition: 

\begin{defi}
\label{sous-carquois}
A {\it Schubert subquiver} of $Q_X$ is a full subquiver of $Q_X$ 
whose set of vertices is an order ideal.
\end{defi}

Note that $w$ itself can easily be recovered from $Q_w$: 
simply  remove from the reduced expression of $w_X$ the
reflexions corresponding to the vertices removed from $Q_X$ to
obtain $Q_w$. 

In the more general case where $X$ is any rational homogeneous space
and we have chosen a reduced expression $\wt_X$ of $w_X$ giving a
quiver $Q_{\wt_X}$ depending on $\wt_X$, we can still define the
Schubert subquivers of $Q_{\wt_X}$. As in the (co)minuscule case,
these subquivers correspond to Schubert subvarieties in $X$. We have
the following fact: 

\begin{fact}
  \label{cas-gen}
Let $X$ be any homogeneous variety and $Q_u$ and $Q_v$ two
Schubert subquivers of $Q_{\tilde w_X}$
corresponding to the Schubert varieties
$X(u)$ and $X(v)$.

(\i) If $Q_u\subset Q_v$ then $X(u)\subset X(v)$.

(\i\i) If $X$ is (co)minuscule, the converse is true.
\end{fact}

\begin{proo}
The first part is clear. The second one comes from the fact that any
Schubert subvariety is obtained inductively by removing peaks.
\end{proo}

\smallskip
Now we come to Poincar{\'e} duality for $X$ (co)minuscule. 
The involution $i_X$ on $Q_X$ induces an involution on
the set of subquivers attached to the Schubert classes.
Indeed, we can let 
$$Q_w\mapsto Q_{i_X(w)}=i_X(Q_X-Q_w).$$ 
This is well defined since $i_X$ completely reverses the partial 
order on $Q_X$: thus the set of vertices of $Q_X-Q_w$ is mapped by $i_X$
to an order ideal. 

 \begin{prop}
\label{dualite}
The Schubert classes $[X(w)]$ and
   $[X(i_X(w))]$ are Poincar{\'e} dual in $A^*(X)$.
 \end{prop}

 \begin{proo} 
Let $w\in W_X$. There exists a  
symmetric reduced expression $w_X=s_{\b_1}\cdots s_{\b_N}$ 
such that $s_{\b_{k+1}}\cdots s_{\b_N}$ is a reduced expression 
for $w$. Since $i(\beta)=-w_0(\beta)$, the element in $W_X$ defining 
the Schubert class which is Poincar{\'e} dual to $\sigma(w)$ is
  $$w^*=s_{i(\b_{k+1})}\cdots s_{i(\b_N)}w_X
=s_{\b_{N-k+1}}\cdots s_{\b_N}.$$
This is nothing else than $i_X(w)$.
 \end{proo}

 \begin{example}
\label{exemplebis}
 (\i) For a Grassmannian, we recover the fact that the
   Poincar{\'e} duality is given by the complementarity of partitions 
in the corresponding rectangle. The previous proposition is a 
generalization of this fact.

 (\i\i) Let $X$ be the Cayley plane. Consider the
   Schubert classes  $\sigma(v)=\sigma'_{12}$ and $\sigma(w)=\sigma''_{12}$. 
We have reduced expressions $v=s_{\alpha_5}
   s_{\alpha_4} s_{\alpha_3} s_{\alpha_1}$ and $w=s_{\alpha_2}
   s_{\alpha_4} s_{\alpha_3} s_{\alpha_1}$.
The Poincar{\'e} duals are 
$w^*=i_X(w)=s_{\alpha_3} s_{\alpha_4}
s_{\alpha_5} s_{\alpha_1} s_{\alpha_6} s_{\alpha_3}  s_{\alpha_4}
s_{\alpha_2} s_{\alpha_5} s_{\alpha_4} s_{\alpha_3} s_{\alpha_1}$, and 
$v^*=i_X(v)= s_{\alpha_2} s_{\alpha_4} s_{\alpha_5} s_{\alpha_1}
s_{\alpha_6} s_{\alpha_3} s_{\alpha_4} s_{\alpha_2} s_{\alpha_5}
s_{\alpha_4}s_{\alpha_3} s_{\alpha_1}.$
The corresponding quivers are given by the
following pictures, where on the left (resp. right) in black we have 
the quivers 
$Q_v$ and $Q_w$ (resp. $Q_{v^*}$ and $Q_{w^*}$) and in red their
complements in $Q_{X}$. This means that an arrow between two vertices
of the subquiver is drawn in black, and the other arrows are drawn in
red. The same convention will be used in all the article.
\vskip 0.5 cm

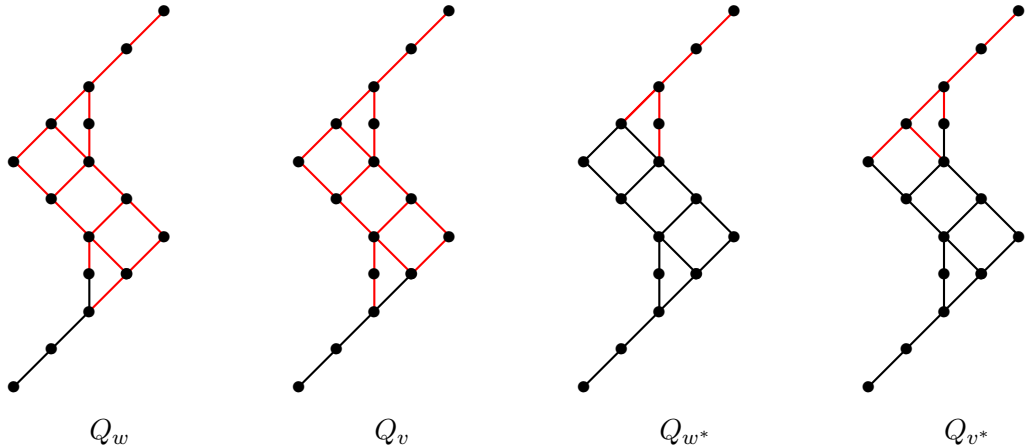
\psset{xunit=0.5cm}
\psset{yunit=0.5cm}
\centerline{\begin{pspicture*}(-5,-5.8)(0.3,6.3)
\psline[linecolor=red](-0,0)(-1,1)
\psline[linecolor=red](-0,6)(-3,3)
\psline[linecolor=red](-2,4)(-2,3)
\psline[linecolor=red](-2,2)(-1,1)
\psline[linecolor=red](-2,2)(-3,3)
\psline[linecolor=red](-1,-1)(-4,2)
\psline[linecolor=red](-2,2)(-2,3)
\psline[linecolor=red](-3,3)(-4,2)
\psline[linecolor=red](-2,2)(-3,1)
\psline[linecolor=red](-1,1)(-2,0)
\psline[linecolor=red](-0,0)(-2,-2)
\psline(-2,-2)(-4,-4)
\psline[linecolor=red](-2,-1)(-2,0)
\psline(-2,-1)(-2,-2)
\put(-0.1,-0.1){$\bullet$}
\put(-0.6,0.4){$\bullet$}
\put(-1.1,0.9){$\bullet$}
\put(-1.6,1.4){$\bullet$}
\put(-1.1,-0.1){$\bullet$}
\put(-1.6,0.4){$\bullet$}
\put(-2.1,0.9){$\bullet$}
\put(-0.6,-0.6){$\bullet$}
\put(-1.1,1.4){$\bullet$}
\put(-1.1,-0.6){$\bullet$}
\put(-2.1,-2.1){$\bullet$}
\put(-1.6,-1.6){$\bullet$}
\put(-1.1,-1.1){$\bullet$}
\put(-0.6,-0.6){$\bullet$}
\put(-0.1,2.9){$\bullet$}
\put(-0.6,2.4){$\bullet$}
\put(-1.1,1.9){$\bullet$}
\put(-1,-2.7){$Q_w$}
\end{pspicture*}
\begin{pspicture*}(-7,-5.8)(0.3,6.3)
\psline[linecolor=red](-0,0)(-1,1)
\psline[linecolor=red](-0,6)(-3,3)
\psline[linecolor=red](-2,4)(-2,3)
\psline[linecolor=red](-2,2)(-1,1)
\psline[linecolor=red](-2,2)(-3,3)
\psline[linecolor=red](-1,-1)(-4,2)
\psline[linecolor=red](-2,2)(-2,3)
\psline[linecolor=red](-3,3)(-4,2)
\psline[linecolor=red](-2,2)(-3,1)
\psline[linecolor=red](-1,1)(-2,0)
\psline[linecolor=red](-0,0)(-1,-1)
\psline(-1,-1)(-4,-4)
\psline[linecolor=red](-2,0)(-2,-2)
\put(-0.1,-0.1){$\bullet$}
\put(-0.6,0.4){$\bullet$}
\put(-1.1,0.9){$\bullet$}
\put(-1.6,1.4){$\bullet$}
\put(-1.1,-0.1){$\bullet$}
\put(-1.6,0.4){$\bullet$}
\put(-2.1,0.9){$\bullet$}
\put(-0.6,-0.6){$\bullet$}
\put(-1.1,1.4){$\bullet$}
\put(-1.1,-0.6){$\bullet$}
\put(-2.1,-2.1){$\bullet$}
\put(-1.6,-1.6){$\bullet$}
\put(-1.1,-1.1){$\bullet$}
\put(-0.6,-0.6){$\bullet$}
\put(-0.1,2.9){$\bullet$}
\put(-0.6,2.4){$\bullet$}
\put(-1.1,1.9){$\bullet$}
\put(-1,-2.7){$Q_v$}
\end{pspicture*}
\begin{pspicture*}(-7,-5.8)(0.3,6.3)
\psline[linecolor=red](0,6)(-3,3)
\psline(-0,0)(-1,1)
\psline[linecolor=red](-2,4)(-3,3)
\psline[linecolor=red](-2,4)(-2,3)
\psline(-2,2)(-1,1)
\psline(-2,2)(-3,3)
\psline(-1,-1)(-4,2)
\psline[linecolor=red](-2,2)(-2,3)
\psline(-3,3)(-4,2)
\psline(-2,2)(-3,1)
\psline(-1,1)(-2,0)
\psline(-0,0)(-4,-4)
\psline(-2,0)(-2,-2)
\put(-0.1,-0.1){$\bullet$}
\put(-0.6,0.4){$\bullet$}
\put(-1.1,0.9){$\bullet$}
\put(-1.6,1.4){$\bullet$}
\put(-1.1,-0.1){$\bullet$}
\put(-1.6,0.4){$\bullet$}
\put(-2.1,0.9){$\bullet$}
\put(-0.6,-0.6){$\bullet$}
\put(-1.1,1.4){$\bullet$}
\put(-1.1,-0.6){$\bullet$}
\put(-2.1,-2.1){$\bullet$}
\put(-1.6,-1.6){$\bullet$}
\put(-1.1,-1.1){$\bullet$}
\put(-0.6,-0.6){$\bullet$}
\put(-0.1,2.9){$\bullet$}
\put(-0.6,2.4){$\bullet$}
\put(-1.1,1.9){$\bullet$}
\put(-1,-2.7){$Q_{w^*}$}
\end{pspicture*}
\begin{pspicture*}(-7,-5.8)(0.3,6.3)
\psline[linecolor=red](0,6)(-4,2)
\psline(-0,0)(-1,1)
\psline[linecolor=red](-2,4)(-2,3)
\psline(-2,2)(-1,1)
\psline[linecolor=red](-2,2)(-3,3)
\psline(-1,-1)(-4,2)
\psline(-2,2)(-2,3)
\psline(-2,2)(-3,1)
\psline(-1,1)(-2,0)
\psline(-0,0)(-4,-4)
\psline(-2,0)(-2,-2)
\put(-0.1,-0.1){$\bullet$}
\put(-0.6,0.4){$\bullet$}
\put(-1.1,0.9){$\bullet$}
\put(-1.6,1.4){$\bullet$}
\put(-1.1,-0.1){$\bullet$}
\put(-1.6,0.4){$\bullet$}
\put(-2.1,0.9){$\bullet$}
\put(-0.6,-0.6){$\bullet$}
\put(-1.1,1.4){$\bullet$}
\put(-1.1,-0.6){$\bullet$}
\put(-2.1,-2.1){$\bullet$}
\put(-1.6,-1.6){$\bullet$}
\put(-1.1,-1.1){$\bullet$}
\put(-0.6,-0.6){$\bullet$}
\put(-0.1,2.9){$\bullet$}
\put(-0.6,2.4){$\bullet$}
\put(-1.1,1.9){$\bullet$}
\put(-1,-2.7){$Q_{v^*}$}
\end{pspicture*}}
 \end{example}

\sectionplus{From classical to quantum invariants} 

In this section, we give a unified presentation of some results of 
A. Buch, A. Kresch and H. Tamvakis \cite{BKT}, according to which 
the Gromov-Witten invariants of (co)minuscule
homogeneous varieties under some group $G$, 
can be computed as classical invariants on some 
other $G$-homogeneous varieties. In particular we extend these
results to the exceptional cases (corollary \ref{gw-d}).

\subsectionplus{Gromov-Witten invariants}

First we briefly recall the definition of the quantum Chow
ring of a rational  homogeneous variety  $X=G/P$, where $P$ is a
maximal parabolic subgroup of $G$. We refer to \cite{fulton}
for more details. 

The small (in the following, this adjective will be skipped)
quantum Chow ring $QA^*(G/P)$ is the free 
$\ZZ[q]$-module over the Schubert classes in $A^*(G/P)$, with a
product given by the following formula:
$$[X(u)]*[X(v)]=\sum_{d\ge 0}q^d\sum_{w\in W_X}I_d([X(u)],[X(v)]
,[X(w)])[X(w^*)].$$
The coefficients $I_d([X(u)],[X(v)],[X(w)])$ in this formula are the
degree $d$ {\it Gromov-Witten invariants}, which are defined as intersection
numbers on the moduli space of stable curves. In the homogeneous
setting, such an invariant counts the number of pointed maps 
$f :\PP^1\rightarrow G/P$ of degree $d$, such that $f(0)\in g\cdot
X(u)$, $f(1)\in g'\cdot X(v)$ and $f(\infty)\in g''\cdot X(w)$
for three general elements $g,g',g''$ of $G$ (see \cite{fulton}, Lemma
13). In fact $I_d([X(u)],[X(v)],[X(w)])$ can be non zero only when
$$\codim(X(u))+\codim(X(v))+\codim(X(w))=\dim(X)+dc_1(X).$$
The quantum product makes of $QA^*(X)$ a commutative and associative 
graded ring, $q$ being given the degree $c_1(X)$ and the Schubert 
classes their codimensions, as in the usual Chow ring.

\subsectionplus{Lines and quantum cohomology}

The associativity of the quantum Chow ring $QA^*(G/P)$ is a very
strong property, and a nice consequence is that the quantum product
can often be completely determined by the classical product and a
small list of quantum multiplications. 

We begin with a simple observation which gives almost for free such a 
quantum product, and makes the link with the discussion of the index 
in section \ref{section_minuscule}. 

Recall that we denoted by $F$ the Fano variety of lines in
$X=G/P\subset \PP V$, with respect to its minimal homogeneous embedding. 
Also we denoted by $F_o\subset F$ the closed subscheme of
 lines through the base point $o$. Mapping such a line to its
 tangent direction at $o$ yields a closed embedding of $F_o$ in $\p
 T_oX$. In the
 following proposition, the degree of $F_o$ is understood 
with respect to  this embedding. 

 \begin{prop}
 \label{h}
The quantum power of the hyperplane class $H$ of exponent $c_1(X)$ is
$$H^{*c_1(X)}=H^{c_1(X)}+\deg(F_o)q.$$
 \end{prop}

 \begin{proo}
 Let $L_1,L_2 \subset \p V$ be general linear subspaces of codimension $l_1,l_2$
 with $l_1 + l_2 = \dim(F) - \dim(X) + 3$. We need to prove that 
 $I_1([o],[L_1],[L_2]) = \deg(F_o)$. As we have seen,
since $\codim(o) + \codim(L_1 \cap X) + \codim(L_2 \cap X) = \dim(F) +
 3 = \dim(X) + c_1(X)$, the Gromov-Witten invariant
$I_1([o],[L_1],[L_2])$ counts the number of lines in $X$ through $o$
meeting $L_1$ and $L_2$. 

Let $\widehat L_i \subset V$ denote the affine cone over $L_i$ and let
 $p:\widehat{T_oX} \rightarrow T_oX = \widehat{T_oX}/\ell_o$ be the natural
 projection, where $\ell_o\subset V$ is the line defined by $o$.
 In terms of the projective embedding $F_o \subset \p T_oX$, a line
 meets $L_1$ and $L_2$ \iff the corresponding point lies in
 $L := \p (p(\widehat{L_1} \cap \widehat{T_oX}) \cap 
 p(\widehat{L_2} \cap \widehat{T_oX}))$. Since
 $\dim F_o = \dim F - \dim X + 1 = l_1+l_2-2 = \codim_{\p T_oX} L$,
 there are $\deg(F_o)$ such points.  \end{proo}

In fact this proposition is true in a much larger generality. 
The previous proof adapts almost verbatim to any projective variety 
for which quantum cohomology is defined.

 \subsectionplus{Quivers and Poincar{\'e} duality for the Fano variety of
 lines}
\label{subsection_fano}

 As in the first section, let us denote by $F$ the Fano
 variety of lines in $X=G/P$. Denote by $I$ the incidence variety and by
 $p$ and $q$ the projections from $I$ to $X$ and $F$:
 $$\xymatrix{I\ar[r]^p\ar[d]^q&X\\
 F&}$$

If $\ell$ is a point in $F$, we denote by $L$ the corresponding line
 in $X$. We have $L=p(q^{-1}(\ell))$. The varieties $F$ and $I$ are
 homogeneous: if $X=G/P$ where $P$ is
 associated to a simple root $\beta$,
then $F=G/Q$ where $Q$
 is associated to the simple roots linked to $\beta$ in
 the Dynkin diagram (see \cite{landsberg}, Theorem 4.8; but beware that this
 is not true for all rational homogeneous varieties). Moreover 
 $I=G/R$ where $R=P\cap Q$.

 Remark that the fiber $Z$ of $p$ is isomorphic to a product of
(co)minuscule homogeneous varieties under a subgroup of
 $G$. In particular there is a uniquely defined quiver $Q_Z=Q_{w_Z}$ 
for $Z$,  constructed from any reduced
 expression of $w_Z$. 

We have $w_I=w_Xw_{Z}$. This can be seen as follows. By the Bruhat 
decomposition, $Bw_ZR$ is open in $P$. By the projection map to
$X=G/P$, the subset $Bw_XBw_ZR/R$ of $I=G/R$ is thus mapped to 
$Bw_XP/P$, the open cell of $X$. The fiber of the base point
is $p^{-1}(w_XP/P)\cap Bw_XBw_ZR/R=(Bw_X\cap w_XP)Bw_ZR/R$, a
dense open subset of $p^{-1}(w_XP/P)=w_XP/R$. Thus $Bw_XBw_ZR/R$
is an open subset of $I$. Its dimension $\ell(w_X)+\ell(w_Z)$ 
is thus  equal to the dimension $\ell(w_I)$ of $I$. In particular, 
$w_Xw_Z$ is reduced and $Bw_XBw_ZR/R$ contains $Bw_Xw_ZR/R$ as an
open subset. But the later must be the open cell in $I$, thus 
$w_Xw_Z=w_I$ as claimed. 

By the same argument, we get that $w_I=s_{\b_1}w_F$,
 where $\b_1$ denotes the unique simple root such that 
$l(s_{\b_1}w)=l(w)-1$.
Now, if we
 choose reduced expressions $s_{\b_1}\cdots s_{\b_N}$ for $w_X$ and
 $s_{\b'_1}\cdots s_{\b'_M}$ for $w_Z$ (with $M=\dim(Z)$), we obtain a
 reduced expression
 $s_{\b_2}\cdots s_{\b_{N}}s_{\b'_1}\cdots s_{\b'_M}$ for $w_F$. 
This reduced expression is uniquely defined modulo commutation
relations, although it is not true that $w_F$ has a unique 
reduced  expression modulo commutation relations. We thus get a 
uniquely defined quiver $Q_F=Q_{w_F}$. Of course we also have a quiver 
$Q_I$, which is deduced from $Q_F$ by attaching a vertex
(corresponding to $\b_1$) to the highest vertex of $Q_F$.

As in the  (co)minuscule case, the quiver $Q_F$ is symmetric:
the fact that $w_0w_X$ and $w_0w_F$ are involutions ensures that 
our reduced expression of $w_F$ can be chosen good, in the sense of
 the following definition:
\begin{defi}
\label{good_expression}
A good reduced expression for $w_F$ is a reduced expression 
$$w_F = s_{\ga_1}\cdots s_{\ga_R} = 
s_{\b_2}\cdots s_{\b_{N}}s_{\b'_1}\cdots s_{\b'_M},$$
where $w_X = s_{\b_1}\cdots s_{\b_{N}}$, $w_Z = s_{\b'_1}\cdots s_{\b'_M}$,
$R=\dim(F)$, and $i(\ga_k)=\ga_{{R+1-k}}$.
\end{defi}

 \begin{example}
   \label{exemple_incidence_e6}
 Let us fix $X$ as in example \ref{exemple} (\i\i). The element $w_Z$
   has the following reduced expression
   $$w_Z= s_{\alpha_2}
   s_{\alpha_4} s_{\alpha_5} s_{\alpha_6} s_{\alpha_3} s_{\alpha_4}
   s_{\alpha_5} s_{\alpha_2} s_{\alpha_4} s_{\alpha_3}$$ 
and we deduce for $w_F$ a good reduced expression 
$$w_F=s_{\alpha_5}s_{\alpha_4}
   s_{\alpha_2} s_{\alpha_3} s_{\alpha_4} s_{\alpha_5}
   s_{\alpha_1} s_{\alpha_3}  s_{\alpha_4} s_{\alpha_2} s_{\alpha_6}
   s_{\alpha_5}
   s_{\alpha_4}s_{\alpha_3} s_{\alpha_1} s_{\alpha_2} s_{\alpha_4}
   s_{\alpha_5} s_{\alpha_6} s_{\alpha_3} s_{\alpha_4} s_{\alpha_5}
   s_{\alpha_2} s_{\alpha_4} s_{\alpha_3}.$$ 
The quivers $Q_{I}$ and
   $Q_{F}$ have the following forms:

\vskip 0.5 cm

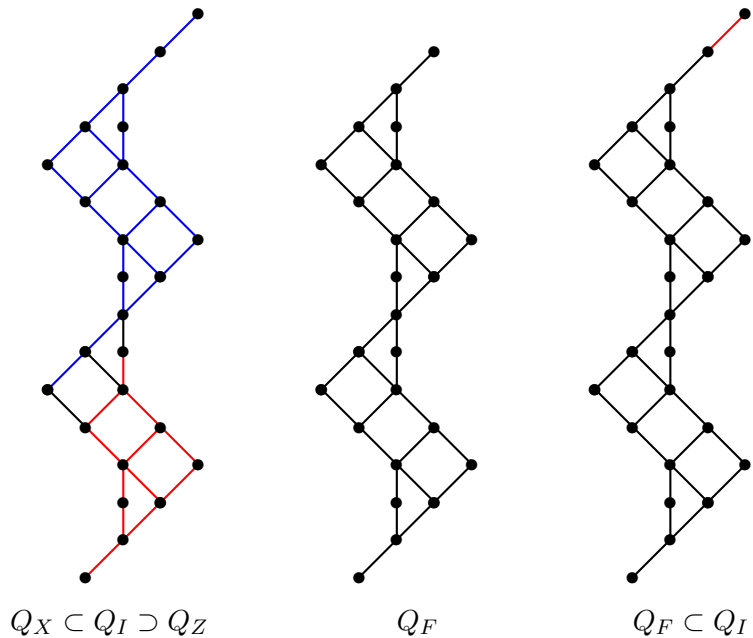
\psset{xunit=0.5cm}
\psset{yunit=0.5cm}
\centerline{\begin{pspicture*}(-5,-11)(2,6.3)
\psline[linecolor=blue](-0,0)(-1,1)
\psline[linecolor=blue](-0,6)(-3,3)
\psline[linecolor=blue](-2,4)(-2,3)
\psline[linecolor=blue](-2,2)(-1,1)
\psline[linecolor=blue](-2,2)(-3,3)
\psline[linecolor=blue](-1,-1)(-4,2)
\psline[linecolor=blue](-2,2)(-2,3)
\psline[linecolor=blue](-3,3)(-4,2)
\psline[linecolor=blue](-2,2)(-3,1)
\psline[linecolor=blue](-1,1)(-2,0)
\psline[linecolor=blue](-0,0)(-4,-4)
\psline[linecolor=blue](-2,0)(-2,-2)
\put(-0.1,-0.1){$\bullet$}
\put(-0.6,0.4){$\bullet$}
\put(-1.1,0.9){$\bullet$}
\put(-1.6,1.4){$\bullet$}
\put(-1.1,-0.1){$\bullet$}
\put(-1.6,0.4){$\bullet$}
\put(-2.1,0.9){$\bullet$}
\put(-0.6,-0.6){$\bullet$}
\put(-1.1,1.4){$\bullet$}
\put(-1.1,-0.6){$\bullet$}
\put(-2.1,-2.1){$\bullet$}
\put(-1.6,-1.6){$\bullet$}
\put(-1.1,-1.1){$\bullet$}
\put(-0.6,-0.6){$\bullet$}
\put(-0.1,2.9){$\bullet$}
\put(-0.6,2.4){$\bullet$}
\put(-1.1,1.9){$\bullet$}
\psline(-2,-2)(-2,-3)
\psline(-3,-3)(-2,-4)
\psline(-4,-4)(-3,-5)
\psline[linecolor=red](-2,-3)(-2,-4)
\psline[linecolor=red](-2,-4)(0,-6)
\psline[linecolor=red](-3,-5)(-1,-7)
\psline[linecolor=red](0,-6)(-3,-9)
\psline[linecolor=red](-2,-6)(-2,-8)
\psline[linecolor=red](-2,-4)(-3,-5)
\psline[linecolor=red](-1,-5)(-2,-6)
\put(-0.1,-3.1){$\bullet$}
\put(-0.6,-2.6){$\bullet$}
\put(-1.1,-2.1){$\bullet$}
\put(-1.6,-1.6){$\bullet$}
\put(-1.1,-3.1){$\bullet$}
\put(-1.6,-2.6){$\bullet$}
\put(-2.1,-2.1){$\bullet$}
\put(-0.6,-3.6){$\bullet$}
\put(-1.1,-1.6){$\bullet$}
\put(-1.1,-3.6){$\bullet$}
\put(-1.6,-4.6){$\bullet$}
\put(-1.1,-4.1){$\bullet$}
\put(-0.6,-3.6){$\bullet$}
\put(-2.5,-5.2){$Q_{X}\subset Q_{I}\supset Q_{Z}$}
\end{pspicture*}
\begin{pspicture*}(-5,-11)(2,6.3)
\psline(-0,0)(-1,1)
\psline(-1,5)(-3,3)
\psline(-2,4)(-2,3)
\psline(-2,2)(-1,1)
\psline(-2,2)(-3,3)
\psline(-1,-1)(-4,2)
\psline(-2,2)(-2,3)
\psline(-3,3)(-4,2)
\psline(-2,2)(-3,1)
\psline(-1,1)(-2,0)
\psline(-0,0)(-4,-4)
\psline(-2,0)(-2,-2)
\put(-0.1,-0.1){$\bullet$}
\put(-0.6,0.4){$\bullet$}
\put(-1.1,0.9){$\bullet$}
\put(-1.6,1.4){$\bullet$}
\put(-1.1,-0.1){$\bullet$}
\put(-1.6,0.4){$\bullet$}
\put(-2.1,0.9){$\bullet$}
\put(-0.6,-0.6){$\bullet$}
\put(-1.1,1.4){$\bullet$}
\put(-1.1,-0.6){$\bullet$}
\put(-2.1,-2.1){$\bullet$}
\put(-1.6,-1.6){$\bullet$}
\put(-1.1,-1.1){$\bullet$}
\put(-0.6,-0.6){$\bullet$}
\put(-0.6,2.4){$\bullet$}
\put(-1.1,1.9){$\bullet$}
\psline(-2,-2)(-2,-4)
\psline(-3,-3)(0,-6)
\psline(-4,-4)(-1,-7)
\psline(0,-6)(-3,-9)
\psline(-2,-6)(-2,-8)
\psline(-2,-4)(-3,-5)
\psline(-1,-5)(-2,-6)
\put(-0.1,-3.1){$\bullet$}
\put(-0.6,-2.6){$\bullet$}
\put(-1.1,-2.1){$\bullet$}
\put(-1.6,-1.6){$\bullet$}
\put(-1.1,-3.1){$\bullet$}
\put(-1.6,-2.6){$\bullet$}
\put(-2.1,-2.1){$\bullet$}
\put(-0.6,-3.6){$\bullet$}
\put(-1.1,-1.6){$\bullet$}
\put(-1.1,-3.6){$\bullet$}
\put(-1.6,-4.6){$\bullet$}
\put(-1.1,-4.1){$\bullet$}
\put(-0.6,-3.6){$\bullet$}
\put(-1,-5.2){$Q_{F}$}
\end{pspicture*}
\begin{pspicture*}(-5,-11)(2,6.3)
\psline[linecolor=red](0,6)(-1,5)
\psline(-0,0)(-1,1)
\psline(-1,5)(-3,3)
\psline(-2,4)(-2,3)
\psline(-2,2)(-1,1)
\psline(-2,2)(-3,3)
\psline(-1,-1)(-4,2)
\psline(-2,2)(-2,3)
\psline(-3,3)(-4,2)
\psline(-2,2)(-3,1)
\psline(-1,1)(-2,0)
\psline(-0,0)(-4,-4)
\psline(-2,0)(-2,-2)
\put(-0.1,-0.1){$\bullet$}
\put(-0.6,0.4){$\bullet$}
\put(-1.1,0.9){$\bullet$}
\put(-1.6,1.4){$\bullet$}
\put(-1.1,-0.1){$\bullet$}
\put(-1.6,0.4){$\bullet$}
\put(-2.1,0.9){$\bullet$}
\put(-0.6,-0.6){$\bullet$}
\put(-1.1,1.4){$\bullet$}
\put(-1.1,-0.6){$\bullet$}
\put(-2.1,-2.1){$\bullet$}
\put(-1.6,-1.6){$\bullet$}
\put(-1.1,-1.1){$\bullet$}
\put(-0.6,-0.6){$\bullet$}
\put(-0.1,2.9){$\bullet$}
\put(-0.6,2.4){$\bullet$}
\put(-1.1,1.9){$\bullet$}
\psline(-2,-2)(-2,-4)
\psline(-3,-3)(0,-6)
\psline(-4,-4)(-1,-7)
\psline(0,-6)(-3,-9)
\psline(-2,-6)(-2,-8)
\psline(-2,-4)(-3,-5)
\psline(-1,-5)(-2,-6)
\put(-0.1,-3.1){$\bullet$}
\put(-0.6,-2.6){$\bullet$}
\put(-1.1,-2.1){$\bullet$}
\put(-1.6,-1.6){$\bullet$}
\put(-1.1,-3.1){$\bullet$}
\put(-1.6,-2.6){$\bullet$}
\put(-2.1,-2.1){$\bullet$}
\put(-0.6,-3.6){$\bullet$}
\put(-1.1,-1.6){$\bullet$}
\put(-1.1,-3.6){$\bullet$}
\put(-1.6,-4.6){$\bullet$}
\put(-1.1,-4.1){$\bullet$}
\put(-0.6,-3.6){$\bullet$}
\put(-1.5,-5.2){$Q_{F}\subset Q_{I}$}
\end{pspicture*}}

\noindent
where on the left picture we have drawn in blue the quiver $Q_{X}$
 embedded in $Q_{I}$ and in red the quiver $Q_{Z}$ embedded in
 $Q_{I}$. In the right one we have drawn in black the quiver
 $Q_{F}$ embedded in $Q_{I}$ and in red its complement. The
 middle picture is $Q_{F}$.
 \end{example}

 As in  subsection \ref{subsection_quivers} 
 the symmetry of the quiver $Q_F$ induces an
 involution $i_F : Q_w\mapsto Q_{i_F(w)}$ on the Schubert subquivers
 (recall definition \ref{sous-carquois} of $Q_{F}$). 
The same proof as that of Proposition 
\ref{dualite} gives the following partial interpretation of 
Poincar{\'e} duality on $F$:

 \begin{prop}
\label{dualite-F}
   Let $F(w)$ be a Schubert subvariety of $F$ such that
   there exists a good reduced expression
 $w_F = s_{\gamma_1} \cdots s_{\gamma_R}$ and an integer $k$ such that
$w=s_{\ga_{k+1}}\cdots s_{\ga_R}$.
Then the classes $[F(w)]$ and
   $[F(i_F(w))]$ are Poincar{\'e} dual.
\end{prop} 

\begin{rema}
(\i) Beware  that not all Schubert varieties $F(w), w\in W_F$ satisfy the
hypothesis of the proposition. This is
because $F$ is not minuscule and  in consequence there may be 
braid relations. However, all Schubert varieties $F(w)$ associated to a
Schubert subquiver $Q_F(w)$ of $Q_F$ satisfy the property.

(\i\i) We will denote by $F(u^\star)$ the Poincar{\'e} dual of $F(u)$.
\end{rema}

Let $X(w)$ be a Schubert variety in $X$. We define the Schubert
variety $F(\wh)=q(p^{-1}(X(w)))$ of $F$. The variety $F(\wh)$ parametrises
the lines in $X$ meeting $X(w)$. For any Schubert variety $X(w)$,
we can choose a reduced expression $w_X = s_{\beta_1} \cdots s_{\beta_N}$
such that $w$ has a reduced expression $s_{\b_{k+1}}\cdots
s_{\b_N}$. Then we get for $\wh$ the reduced expression $\wh=s_{\b_{k+1}}\cdots
s_{\b_N}w_Z= s_{\ga_{k+1}}\cdots s_{\ga_R}$. 
In particular $F(\wh)$ 
satisfies the hypothesis of proposition \ref{dualite-F}.

Consider the quivers $Q_{X}$ and
$Q_{F}$ as subquivers of $Q_{I}$: the first one is
obtained by removing $Q_{Z}$ from the bottom of  $Q_{I}$;  the second
one by removing the top vertex of $Q_{I}$. The quiver of $F(\wh)$
is deduced from the quiver of $X(w)$ by removing from $Q_{F}\subset
Q_{I}$ the vertices of $Q_{X}$ not contained 
in $Q_w$.
In particular, if $X(w)$ is different from $X$, then
$$\codim_F(F(\wh))=\codim_X(X(w))-1.$$ Remark that a Schubert subvariety
$F(v)$, where $v=s_{\ga_{k+1}}\cdots s_{\ga_R}$, is of the form $F(\wh)$
for some $w$, if and only if the quiver $Q_{F(v)}$  contains $Q_{Z}$. 
In that case $Q_w$ is the complement of $Q_{Z}$ in $Q_{F(v)}$.

\begin{example}
Let $X$, $X(w)$ and $X(v)$ as in example \ref{exemplebis} (\i\i). Then
the quivers of $F(\wh)$ and $F(\vh)$ are the following (in black):

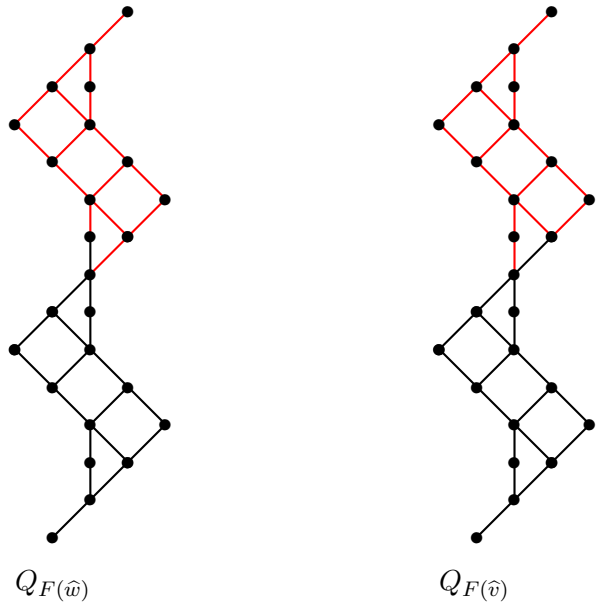
\psset{xunit=0.5cm}
\psset{yunit=0.5cm}
\centerline{\begin{pspicture*}(-5,-11)(2,6.3)
\psline[linecolor=red](-0,0)(-1,1)
\psline[linecolor=red](-2,4)(-2,3)
\psline[linecolor=red](-2,2)(-1,1)
\psline[linecolor=red](-2,2)(-3,3)
\psline[linecolor=red](-1,-1)(-4,2)
\psline[linecolor=red](-2,2)(-2,3)
\psline[linecolor=red](-3,3)(-4,2)
\psline[linecolor=red](-2,2)(-3,1)
\psline[linecolor=red](-1,1)(-2,0)
\psline[linecolor=red](-0,0)(-2,-2)
\psline(-2,-2)(-4,-4)
\psline[linecolor=red](-2,-1)(-2,0)
\psline(-2,-1)(-2,-2)
\psline[linecolor=red](-1,5)(-3,3)
\put(-0.1,-0.1){$\bullet$}
\put(-0.6,0.4){$\bullet$}
\put(-1.1,0.9){$\bullet$}
\put(-1.6,1.4){$\bullet$}
\put(-1.1,-0.1){$\bullet$}
\put(-1.6,0.4){$\bullet$}
\put(-2.1,0.9){$\bullet$}
\put(-0.6,-0.6){$\bullet$}
\put(-1.1,1.4){$\bullet$}
\put(-1.1,-0.6){$\bullet$}
\put(-2.1,-2.1){$\bullet$}
\put(-1.6,-1.6){$\bullet$}
\put(-1.1,-1.1){$\bullet$}
\put(-0.6,-0.6){$\bullet$}
\put(-0.6,2.4){$\bullet$}
\put(-1.1,1.9){$\bullet$}
\psline(-2,-2)(-2,-4)
\psline(-3,-3)(0,-6)
\psline(-4,-4)(-1,-7)
\psline(0,-6)(-3,-9)
\psline(-2,-6)(-2,-8)
\psline(-2,-4)(-3,-5)
\psline(-1,-5)(-2,-6)
\put(-0.1,-3.1){$\bullet$}
\put(-0.6,-2.6){$\bullet$}
\put(-1.1,-2.1){$\bullet$}
\put(-1.6,-1.6){$\bullet$}
\put(-1.1,-3.1){$\bullet$}
\put(-1.6,-2.6){$\bullet$}
\put(-2.1,-2.1){$\bullet$}
\put(-0.6,-3.6){$\bullet$}
\put(-1.1,-1.6){$\bullet$}
\put(-1.1,-3.6){$\bullet$}
\put(-1.6,-4.6){$\bullet$}
\put(-1.1,-4.1){$\bullet$}
\put(-0.6,-3.6){$\bullet$}
\put(-2,-5.2){$Q_{F(\wh)}$}
\end{pspicture*}
\hskip 2 cm
\begin{pspicture*}(-5,-11)(2,6.3)
\psline[linecolor=red](-0,0)(-1,1)
\psline[linecolor=red](-2,4)(-2,3)
\psline[linecolor=red](-2,2)(-1,1)
\psline[linecolor=red](-2,2)(-3,3)
\psline[linecolor=red](-1,-1)(-4,2)
\psline[linecolor=red](-2,2)(-2,3)
\psline[linecolor=red](-3,3)(-4,2)
\psline[linecolor=red](-2,2)(-3,1)
\psline[linecolor=red](-1,1)(-2,0)
\psline[linecolor=red](-0,0)(-1,-1)
\psline(-1,-1)(-4,-4)
\psline[linecolor=red](-2,0)(-2,-2)
\psline[linecolor=red](-1,5)(-3,3)
\put(-0.1,-0.1){$\bullet$}
\put(-0.6,0.4){$\bullet$}
\put(-1.1,0.9){$\bullet$}
\put(-1.6,1.4){$\bullet$}
\put(-1.1,-0.1){$\bullet$}
\put(-1.6,0.4){$\bullet$}
\put(-2.1,0.9){$\bullet$}
\put(-0.6,-0.6){$\bullet$}
\put(-1.1,1.4){$\bullet$}
\put(-1.1,-0.6){$\bullet$}
\put(-2.1,-2.1){$\bullet$}
\put(-1.6,-1.6){$\bullet$}
\put(-1.1,-1.1){$\bullet$}
\put(-0.6,-0.6){$\bullet$}
\put(-0.6,2.4){$\bullet$}
\put(-1.1,1.9){$\bullet$}
\psline(-2,-2)(-2,-4)
\psline(-3,-3)(0,-6)
\psline(-4,-4)(-1,-7)
\psline(0,-6)(-3,-9)
\psline(-2,-6)(-2,-8)
\psline(-2,-4)(-3,-5)
\psline(-1,-5)(-2,-6)
\put(-0.1,-3.1){$\bullet$}
\put(-0.6,-2.6){$\bullet$}
\put(-1.1,-2.1){$\bullet$}
\put(-1.6,-1.6){$\bullet$}
\put(-1.1,-3.1){$\bullet$}
\put(-1.6,-2.6){$\bullet$}
\put(-2.1,-2.1){$\bullet$}
\put(-0.6,-3.6){$\bullet$}
\put(-1.1,-1.6){$\bullet$}
\put(-1.1,-3.6){$\bullet$}
\put(-1.6,-4.6){$\bullet$}
\put(-1.1,-4.1){$\bullet$}
\put(-0.6,-3.6){$\bullet$}
\put(-2,-5.2){$Q_{F(\vh)}$}
\end{pspicture*}}
\end{example}

\subsectionplus{Computing degree one Gromov-Witten invariants}

\label{subsection_computing}

In this subsection, we explain how to apply the technics of A. Buch,
A. Kresch and H. Tamvakis \cite{BKT} to calculate degree one
Gromov-Witten invariants.

\begin{lemm}
\label{intersection}
  Let $X(w)$ and $X(v^*)$ be two Schubert varieties such that
  $X(w)\subset X(v^*)$. Then there exists an element $g\in G$ such that
  the intersection $X(w)\cap g\cdot X(v)$ is a reduced point.
\end{lemm}

\begin{proo}
We proceed by induction on $a=\dim(X(v^*))-\dim(X(w))$. If $a=0$ this
is just Poincar{\'e} duality. Assume that the result holds for $X(w')$
such that $X(w)$ is a divisor in $X(w')$. Then there exists an
element $g\in G$ such that $X(w')\cap g\cdot X(v)$ is a reduced
point $x$. Now we use the following fact on (co)minuscule 
Schubert varieties (cf. \cite{LMS}):

\begin{fact}
  Any divisor $X(w)$ in $X(w')$ is a moving divisor, in particular we
  have 
$$\displaystyle{X(w')=\bigcup_{h\in {\rm Stab}(X(w'))}h\cdot
  X(w).}$$
\end{fact}

We deduce from this fact that there exist an element $h\in{\rm
  Stab}(X(w'))$ such that $x\in h\cdot X(w)$. But now $h\cdot X(w)$
  meets $g\cdot X(v)$ in $x$ at least. Since $h\cdot
  X(w)$ is contained in $X(w')$ and $X(w')\cap g\cdot
  X(v)$ is the reduced point $x$, $h\cdot X(w)$ and $g\cdot
  X(v)$ meet only in $x$, and transversely at that point.
\end{proo}

\begin{rema}
 Applying the lemma to the case where $X(v^*)$ is the codimension one
  Schubert subvariety of $X$, we obtain that if $X(w)$ is different from
  $X$ and $\ell$ is a general point in $F(\wh)$, then $L$ meets $X(w)$
  in a unique point.   
In this paragraph we will only use this version of the lemma. In
section \ref{subsection_higher}, we use it in its general formulation.
\end{rema}

\begin{lemm}
Let $X(u)$, $X(v)$ and $X(w)$ be three proper Schubert subvarieties 
of $X$,  such that 
$$\codim(X(u))+\codim(X(v))+\codim(X(w)) = \dim(X)+c_1(X).$$
Then for $g$, $g'$ and $g''$ three general elements in $G$, the
intersection $g\cdot F(\uh)\cap g'\cdot F(\vh)\cap g''\cdot F(\wh)$ is
a finite set of reduced points. 

Let $\ell$ be a point in this
intersection, then the line $L$ meets each of $g\cdot X(u)$, $g'\cdot X(v)$ and
$g''\cdot X(w)$ in a unique point and these points are in general
position in $L$.
\end{lemm}

\begin{proo}
The codimension condition and the fact that $X(u)$, $X(v)$ and $X(w)$
  are different from $X$ imply that 
$$\codim(F(\uh))+\codim(F(\vh))+\codim(F(\wh)) = \dim(F).$$
In particular, the first part of the proposition follows from Bertini's
theorem (see \cite{kleiman}). Furthermore, by Bertini again, we may
assume that any $\ell$ in the intersection is general in $g\cdot
F(\uh)$, $g'\cdot F(\vh)$ and $g''\cdot F(\wh)$. In particular, $L$
meets each of $g\cdot X(u)$, $g'\cdot X(v)$ and $g''\cdot X(w)$ in a
unique point by the previous lemma. Finally, the stabiliser of $\ell$
acts transitively on $L$ and by modifying $g$, $g'$ and $g''$ by
elements in this stabiliser me may assume that the points are in
general position in $L$.  
\end{proo}

\begin{coro}
\label{XF_degre1}
Let $X(u)$, $X(v)$ and $X(w)$ be three proper Schubert subvarieties of
$X$. Suppose that the sum of their codimensions is $\dim(X)+c_1(X)$.  
Then
$$I_1([X(u)],[X(v)],[X(w)])=I_0([F(\uh)],[F(\vh)],[F(\wh)]).$$
\end{coro}

\begin{proo}
The image of any morphism counting in the invariant
$I_1([X(u)],[X(v)],[X(w)])$ is a line $\ell$ in the intersection
$g\cdot F(\uh)\cap g'\cdot 
F(\vh)\cap g''\cdot F(\wh)$ for general elements $g$, $g'$ and $g''$
in $G$. The preceding lemma implies that there is a finite number
$I_0([F(\uh)],[F(\vh)],[F(\wh)])$ of such lines and that given such a
line, there exists a unique morphism from $\p^1$ to $L$ with 0, 1 and
$\infty$ sent to the intersection of $L$ with $g\cdot X(u)$, $g'\cdot
X(v)$ and $g''\cdot X(w)$. 
\end{proo}

In order to make some computations on $F$ we prove the following lemma.

\begin{lemm}
\label{annul}
  In any homogeneous variety $Y$, 
  the intersection product $[Y(u)]\cdot[Y(v)]$ vanishes if and only if
  $Y(v)$ does not contain $Y(u^\star)$ where $\star$ is the Poincar{\'e}
  duality on $Y$.
\end{lemm}

\begin{proo}
The product $[Y(u)]\cdot[Y(v)]$ does not vanish if and only there
exists a sequence of codimension one Schubert subvarieties 
$H_1,\ldots,H_a$, with
$a=\dim(Y)-\codim(Y(u))-\codim(Y(v))$, such that the
product $[H_1]\cdots [H_a]\cdot [Y(u)]\cdot[Y(v)]$ is
non zero (recall that ${\rm Pic}(Y)$ may be bigger than $\Z$).
But we have
$$[H_1]\cdots [H_a]\cdot[Y(v)]=\sum_wc^{H_1\cdots
  H_a}_w[Y(w)]$$
where the sum runs over all $w$ such that
$\dim(Y(w))=\codim(Y(u))$. The product
$[H_1]\cdots [H_a]\cdot [Y(u)]\cdot[Y(v)]$ is non zero if and
only if $c^{H_1\cdots H_a}_{u^\star}$ is non zero. 

As a consequence the product $[Y(u)]\cdot[Y(v)]$ does not vanish if
and only there exists a sequence of codimension one Schubert varieties
$H_1,\cdots,H_a$ such that $c^{H_1\cdots H_a}_{u^\star}\neq0$. But this is
equivalent to the inclusion $Y(u^\star)\subset Y(v)$.
\end{proo}

Using this lemma on a (co)minuscule Schubert variety $X$ will be easy
because inclusions of Schubert varieties are equivalent to 
inclusions of their quivers (Fact \ref{cas-gen}). This is not the case
for a general rational homogeneous variety and in particular for
$F$. However we have the following lemma:

\begin{lemm}
  \label{cas-gen1}
Let $u$ and $v$ in $W_F$ such that $F(u)$ and $F(v)$ are represented
by Schubert subquivers $Q_F(u)$ and $Q_F(v)$ of $Q_F$.

(\i) If $Q_F(u)\subset Q_F(v)$ then $F(u)\subset F(v)$ (see Fact
\ref{cas-gen} (\i)).

(\i\i) Conversely, if $F(u)\subset F(v)$, then we have the inclusion
$Q_F(u)\cap i_F(Q_Z)\subset Q_F(v)\cap i_F(Q_Z)$.
\end{lemm}

\begin{proo}
For (\i\i), remark that we have the
  equality $w_I=w_Fs_{\iota(\beta_1)}$ so that we get a quiver for $I$ by
  adding at the end of $Q_F$ a vertex corresponding to
  $i(\b_1)$. Adding the same vertex at the end of $i_F(Q_Z)$ gives the
  quiver $Q_X$. More generally, adding the same vertex at the end of
  $Q_F(u)\cap i_F(Q_Z)$
gives the quiver of $p(q^{-1}(F(u))$. In particular, if
  $F(u)\subset F(v)$ we must have $p(q^{-1}(F(u))\subset
  p(q^{-1}(F(v))$ and an inclusion of the corresponding quivers
  because $X$ is minuscule. This
  gives the desired condition.
\end{proo}

\begin{example}\label{sigma8}
  We come back to example \ref{exemple} (\i\i).   Let $X(v)$ be the Schubert
  subvariety corresponding to the class $\sigma_8$. The quiver $Q_v$  
  has the following form:

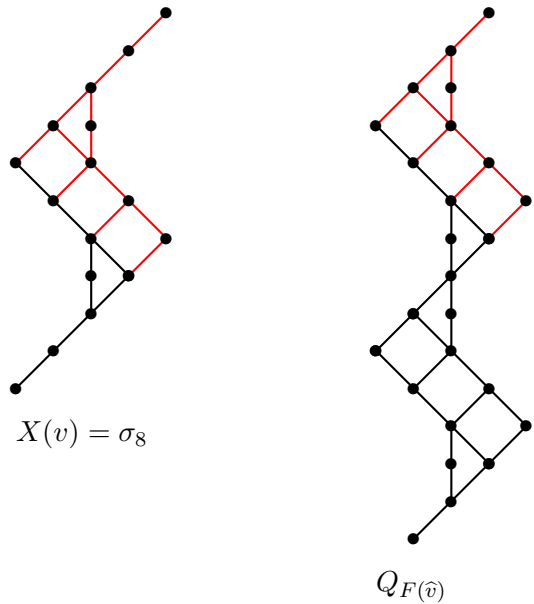
\psset{xunit=0.5cm}
\psset{yunit=0.5cm}
\centerline{\begin{pspicture*}(-5,-10)(0.3,6.3)
\psline[linecolor=red](-0,0)(-1,1)
\psline[linecolor=red](-0,6)(-3,3)
\psline[linecolor=red](-2,4)(-2,3)
\psline[linecolor=red](-2,2)(-1,1)
\psline[linecolor=red](-2,2)(-3,3)
\psline(-1,-1)(-4,2)
\psline[linecolor=red](-2,2)(-2,3)
\psline[linecolor=red](-3,3)(-4,2)
\psline[linecolor=red](-2,2)(-3,1)
\psline[linecolor=red](-1,1)(-2,0)
\psline(-1,-1)(-2,-2)
\psline[linecolor=red](0,0)(-1,-1)
\psline(-2,-2)(-4,-4)
\psline(-2,-1)(-2,0)
\psline(-2,-1)(-2,-2)
\put(-0.1,-0.1){$\bullet$}
\put(-0.6,0.4){$\bullet$}
\put(-1.1,0.9){$\bullet$}
\put(-1.6,1.4){$\bullet$}
\put(-1.1,-0.1){$\bullet$}
\put(-1.6,0.4){$\bullet$}
\put(-2.1,0.9){$\bullet$}
\put(-0.6,-0.6){$\bullet$}
\put(-1.1,1.4){$\bullet$}
\put(-1.1,-0.6){$\bullet$}
\put(-2.1,-2.1){$\bullet$}
\put(-1.6,-1.6){$\bullet$}
\put(-1.1,-1.1){$\bullet$}
\put(-0.6,-0.6){$\bullet$}
\put(-0.1,2.9){$\bullet$}
\put(-0.6,2.4){$\bullet$}
\put(-1.1,1.9){$\bullet$}
\put(-2,-2.7){$X(v)=\sigma_{8}$}
\end{pspicture*}
\hskip 2 cm
\begin{pspicture*}(-5,-11)(2,6.3)
\psline[linecolor=red](-0,0)(-1,1)
\psline[linecolor=red](-2,4)(-2,3)
\psline[linecolor=red](-2,2)(-1,1)
\psline[linecolor=red](-2,2)(-3,3)
\psline(-1,-1)(-4,2)
\psline[linecolor=red](-2,2)(-2,3)
\psline[linecolor=red](-3,3)(-4,2)
\psline[linecolor=red](-2,2)(-3,1)
\psline[linecolor=red](-1,1)(-2,0)
\psline[linecolor=red](-0,0)(-1,-1)
\psline(-1,-1)(-4,-4)
\psline(-2,0)(-2,-2)
\psline[linecolor=red](-1,5)(-3,3)
\put(-0.1,-0.1){$\bullet$}
\put(-0.6,0.4){$\bullet$}
\put(-1.1,0.9){$\bullet$}
\put(-1.6,1.4){$\bullet$}
\put(-1.1,-0.1){$\bullet$}
\put(-1.6,0.4){$\bullet$}
\put(-2.1,0.9){$\bullet$}
\put(-0.6,-0.6){$\bullet$}
\put(-1.1,1.4){$\bullet$}
\put(-1.1,-0.6){$\bullet$}
\put(-2.1,-2.1){$\bullet$}
\put(-1.6,-1.6){$\bullet$}
\put(-1.1,-1.1){$\bullet$}
\put(-0.6,-0.6){$\bullet$}
\put(-0.6,2.4){$\bullet$}
\put(-1.1,1.9){$\bullet$}
\psline(-2,-2)(-2,-4)
\psline(-3,-3)(0,-6)
\psline(-4,-4)(-1,-7)
\psline(0,-6)(-3,-9)
\psline(-2,-6)(-2,-8)
\psline(-2,-4)(-3,-5)
\psline(-1,-5)(-2,-6)
\put(-0.1,-3.1){$\bullet$}
\put(-0.6,-2.6){$\bullet$}
\put(-1.1,-2.1){$\bullet$}
\put(-1.6,-1.6){$\bullet$}
\put(-1.1,-3.1){$\bullet$}
\put(-1.6,-2.6){$\bullet$}
\put(-2.1,-2.1){$\bullet$}
\put(-0.6,-3.6){$\bullet$}
\put(-1.1,-1.6){$\bullet$}
\put(-1.1,-3.6){$\bullet$}
\put(-1.6,-4.6){$\bullet$}
\put(-1.1,-4.1){$\bullet$}
\put(-0.6,-3.6){$\bullet$}
\put(-2,-5.2){$Q_{F(\vh)}$}
\end{pspicture*}}

On the left we have drawn in black the quiver $Q_v$ inside
$Q_{X}$ with its complement in red. On the right we have drawn in
black the quiver of $F(\vh)$ inside $Q_{F}$ with its complement in
red.

Now let $X(u)$ be any Schubert subvariety of codimension four.
For degree reasons, the quantum product
  $X(u)*X(v)$ can only have terms in $q^0$ or $q^1$. Let us
  concentrate on $q$-terms. We need to compute
  $I_0([F(\uh)],[F(\vh)],[F(\wh)])$ for $w$ such that
  $\codim(X(w))=16$,  that is,  $X(w)=\{{\rm pt}\}$. 
By Lemma \ref{cas-gen1}, the Schubert variety $F(\vh^\star)$ is not
contained in $F(\wh)$. Indeed, both quivers are contained in
$\iota_F(Q_Z)$ but we don't have $Q_F(\vh^\star)\subset Q_F(\wh)$ 
(see the following pictures).

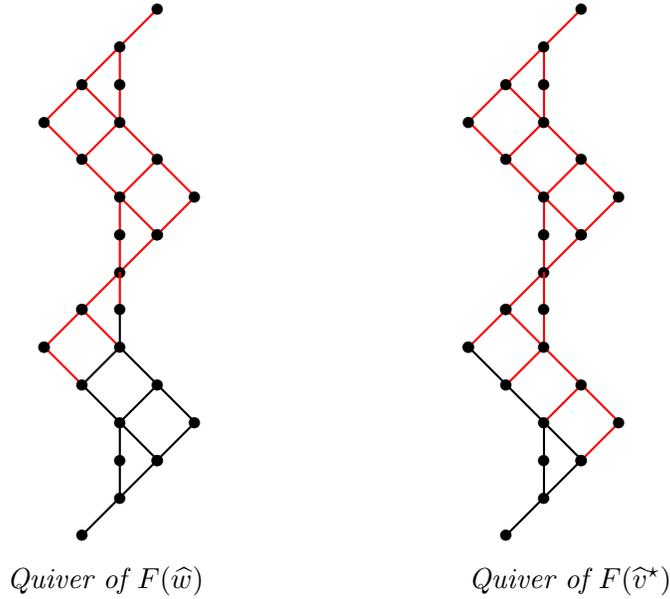
\psset{xunit=0.5cm}
\psset{yunit=0.5cm}
\centerline{
\begin{pspicture*}(-5,-11)(2,6.3)
\psline[linecolor=red](-0,0)(-1,1)
\psline[linecolor=red](-1,5)(-3,3)
\psline[linecolor=red](-2,4)(-2,3)
\psline[linecolor=red](-2,2)(-1,1)
\psline[linecolor=red](-2,2)(-3,3)
\psline[linecolor=red](-1,-1)(-4,2)
\psline[linecolor=red](-2,2)(-2,3)
\psline[linecolor=red](-3,3)(-4,2)
\psline[linecolor=red](-2,2)(-3,1)
\psline[linecolor=red](-1,1)(-2,0)
\psline[linecolor=red](-0,0)(-4,-4)
\psline[linecolor=red](-2,0)(-2,-2)
\put(-0.1,-0.1){$\bullet$}
\put(-0.6,0.4){$\bullet$}
\put(-1.1,0.9){$\bullet$}
\put(-1.6,1.4){$\bullet$}
\put(-1.1,-0.1){$\bullet$}
\put(-1.6,0.4){$\bullet$}
\put(-2.1,0.9){$\bullet$}
\put(-0.6,-0.6){$\bullet$}
\put(-1.1,1.4){$\bullet$}
\put(-1.1,-0.6){$\bullet$}
\put(-2.1,-2.1){$\bullet$}
\put(-1.6,-1.6){$\bullet$}
\put(-1.1,-1.1){$\bullet$}
\put(-0.6,-0.6){$\bullet$}
\put(-0.6,2.4){$\bullet$}
\put(-1.1,1.9){$\bullet$}
\psline[linecolor=red](-2,-2)(-2,-3)
\psline[linecolor=red](-3,-3)(-2,-4)
\psline[linecolor=red](-4,-4)(-3,-5)
\psline(-2,-3)(-2,-4)
\psline(-2,-4)(0,-6)
\psline(-3,-5)(-1,-7)
\psline(0,-6)(-3,-9)
\psline(-2,-6)(-2,-8)
\psline(-2,-4)(-3,-5)
\psline(-1,-5)(-2,-6)
\put(-0.1,-3.1){$\bullet$}
\put(-0.6,-2.6){$\bullet$}
\put(-1.1,-2.1){$\bullet$}
\put(-1.6,-1.6){$\bullet$}
\put(-1.1,-3.1){$\bullet$}
\put(-1.6,-2.6){$\bullet$}
\put(-2.1,-2.1){$\bullet$}
\put(-0.6,-3.6){$\bullet$}
\put(-1.1,-1.6){$\bullet$}
\put(-1.1,-3.6){$\bullet$}
\put(-1.6,-4.6){$\bullet$}
\put(-1.1,-4.1){$\bullet$}
\put(-0.6,-3.6){$\bullet$}
\put(-2.5,-5.2){Quiver of $F(\wh)$}
\end{pspicture*}
\hskip 2 cm
\begin{pspicture*}(-5,-11)(2,6.3)
\psline[linecolor=red](-0,0)(-1,1)
\psline[linecolor=red](-2,4)(-2,3)
\psline[linecolor=red](-2,2)(-1,1)
\psline[linecolor=red](-2,2)(-3,3)
\psline[linecolor=red](-1,-1)(-4,2)
\psline[linecolor=red](-2,2)(-2,3)
\psline[linecolor=red](-3,3)(-4,2)
\psline[linecolor=red](-2,2)(-3,1)
\psline[linecolor=red](-1,1)(-2,0)
\psline[linecolor=red](-0,0)(-1,-1)
\psline[linecolor=red](-1,-1)(-4,-4)
\psline[linecolor=red](-2,0)(-2,-2)
\psline[linecolor=red](-1,5)(-3,3)
\put(-0.1,-0.1){$\bullet$}
\put(-0.6,0.4){$\bullet$}
\put(-1.1,0.9){$\bullet$}
\put(-1.6,1.4){$\bullet$}
\put(-1.1,-0.1){$\bullet$}
\put(-1.6,0.4){$\bullet$}
\put(-2.1,0.9){$\bullet$}
\put(-0.6,-0.6){$\bullet$}
\put(-1.1,1.4){$\bullet$}
\put(-1.1,-0.6){$\bullet$}
\put(-2.1,-2.1){$\bullet$}
\put(-1.6,-1.6){$\bullet$}
\put(-1.1,-1.1){$\bullet$}
\put(-0.6,-0.6){$\bullet$}
\put(-0.6,2.4){$\bullet$}
\put(-1.1,1.9){$\bullet$}
\psline[linecolor=red](-2,-3)(-2,-4)
\psline[linecolor=red](-2,-2)(-2,-3)
\psline[linecolor=red](-3,-3)(-2,-4)
\psline(-4,-4)(-3,-5)
\psline[linecolor=red](-2,-4)(0,-6)
\psline(-3,-5)(-1,-7)
\psline(-1,-7)(-3,-9)
\psline[linecolor=red](0,-6)(-1,-7)
\psline(-2,-6)(-2,-8)
\psline[linecolor=red](-2,-4)(-3,-5)
\psline[linecolor=red](-1,-5)(-2,-6)
\put(-0.1,-3.1){$\bullet$}
\put(-0.6,-2.6){$\bullet$}
\put(-1.1,-2.1){$\bullet$}
\put(-1.6,-1.6){$\bullet$}
\put(-1.1,-3.1){$\bullet$}
\put(-1.6,-2.6){$\bullet$}
\put(-2.1,-2.1){$\bullet$}
\put(-0.6,-3.6){$\bullet$}
\put(-1.1,-1.6){$\bullet$}
\put(-1.1,-3.6){$\bullet$}
\put(-1.6,-4.6){$\bullet$}
\put(-1.1,-4.1){$\bullet$}
\put(-0.6,-3.6){$\bullet$}
\put(-2,-5.2){Quiver of $F(\vh^\star)$}
\end{pspicture*}}

This implies thanks to lemma \ref{annul} that
  $[F(\vh)]\cdot[F(\wh)]=0$ and in particular
  $I_0([F(\uh)],[F(\vh)],[F(\wh)])=0$.
We conclude that for any codimension four class $\tau$ in
  $A^*(X)$, we have:
$$\tau*\sigma_8=\tau\cdot\sigma_8.$$
\end{example}

\subsectionplus{Some geometry of rational curves}

\label{geometry}

We will give in the next section a way of computing higher 
Gromov-Witten invariants similar to what we did for degree one.
The starting point is to find a variety playing the role of $F$ for 
rational curves of degree $d>1$. We have seen that the fact that $F$ 
is homogeneous plays a crucial role in the proofs. However, the
variety of degree $d$ rational curves on $X$ is not homogeneous for 
$d\geq2$. We will introduce in this section a homogeneous
variety $F_d$ which will play the role of $F$ for $d\geq 2$.

F.L. Zak suggested to study the following integer $d(x,y)$.

\begin{defi}\
\label{y_d}
\begin{itemize}
\item
Let $x,y \in X$. We denote $d(x,y)$
the least integer $\delta$ such that there exists a degree $\delta$
union of rational curves through $x$ and $y$. This gives $X$ the structure
of a metric space.
\item
For $x,y \in X$, we define $Y(x,y)$ as the union of all degree $d(x,y)$
unions of rational curves through $x$ and $y$. Let $Y_d$ denote the
abstract variety $Y(x,y)$, for any couple $(x,y)$ such that $d(x,y)=d$.
\item
Let $\alpha$ denote the root defining $P$ and let 
$d_{max}(X)$ denote the number of occurences of $s_\alpha$ in a reduced
expression of $w_X$.
\item
If $d \in [0,d_{max}]$, we denote $F_d$ the set of all $Y(x,y)$'s,
for $x,y \in X$ such that $d(x,y)=d$.
\end{itemize}
\end{defi}

The varieties $F_d$ and $Y_d$ are well-defined in view of the
following proposition:

\begin{prop}
\label{proposition_yd}
\label{couples}
The group
$G$ acts transitively on the set of couples of points $(x,y)$
with $d(x,y)=d$, and $d(x,y)$ takes exactly all 
the values between $0$ and $d_{max}$.
If $\omega \in F_d$, then the stabilizer of 
$\omega$ is a parabolic subgroup of $G$,
thus giving $F_d$ the structure of a projective variety. Moreover,
$Y_d$ and $F_d$ are as in the following table~:
$$\begin{array}{ccccc}
 X & d_{max} & d & F_d  & Y_d \\
 & & & \\
\G(p,n) & \min(p,n-p) & & \F(p-d,p+d;n) &
\G(d,2d)\\
\G_{\omega}(n,2n) & n & & \G_{\omega}(n-d,2n) & \G_{\omega}(d,2d)\\
\G_{Q}(n,2n) & \frac{n}{2} & & \G_{Q}(n-2d,2n) &  \G_{Q}(2d,4d) \\ 
\QQ^n & 2 & d=2 & \{{\rm pt}\} & \QQ^n \\
E_6/P_1 & 2 & d=2 & E_6/P_6 & \QQ^8\\
E_7/P_7 & 3 & d=2 & E_7/P_1 & \QQ^{10}\\
E_7/P_7 & 3 & d=3 & \{{\rm pt}\} & E_7/P_7
\end{array}$$
\end{prop}
\begin{proo}
Let us denote temporarily $F_d',Y_d'$ the varieties in definition
\ref{y_d}, and let $F_d,Y_d$ denote the homogeneous varieties in the above
array. The result of the proposition, that $F_d=F'_d$ and $Y_d=Y'_d$, will
follow from facts about $F_d,Y_d$.
If $X=G/P$, then the variety $F_d$ is a homogeneous variety
under $G$ of the form $F_d=G/Q$ with $Q$ a parabolic subgroup (if
$F_d=\{{\rm pt}\}$ then $Q=G$). In particular we have an incidence
variety $I_d$ and morphisms $p_d:I_d\to X$ and $q_d:I_d\to F_d$ giving
rise to the diagram:
$$\xymatrix{I_d\ar[r]^{p_d}\ar[d]^{q_d}&X\\
 F_d&}$$

If $\omega \in F_d$, we denote
$Y_\omega = p_d(q_d^{-1}(\omega))$.
In the following, $Z_d$ will denote a fiber of $p_d$.
The relevance of $Y_\omega$ for the study of degree $d$ rational
curves comes from the fact:

\begin{prop}
\label{existe}
  For any degree $d$ rational curve $C$, there exists at least one
  element $\omega \in F_d$ such that 
  $C \subset Y_\omega \simeq Y_d$. For a general curve, 
  the point $\omega$ is unique.
\end{prop}

\begin{proo}
Assume first that $C$ is irreducible. Then,
this was observed in \cite{BKT} for (isotropic) Grassmannians. Its is 
obvious when $F_d$ is a point. For the Cayley plane and $d=2$ this 
is Lemma \ref{e6p1coniques}. For the Freudenthal variety and $d=2$ again, this 
is Lemma \ref{e7p7coniques}.

We now extend this result
to the case of a reducible curve. We consider the moduli space
$M_d(X)$ parametrizing rational curves of degree $d$ in $X$, and the
subset $M_d^Y(X)$ of curves included in $Y_\omega$ for some
$\omega \in F_d$. Consider the relative moduli space
${\cal M}_d \rightarrow F_d$
whose fiber over $\omega \in F_d$ parametrizes the curves in
$Y_\omega$. We have a natural map 
${\cal M}_d \rightarrow M_d(X)$, whose image is by definition $M^Y_d(X)$.
Since ${\cal M}_d$ is proper, $M^Y_d(X)$ is closed in $M_d(X)$. We have
seen above that it contains an open subset of $M_d(X)$, which is irreducible
by \cite{Thomsen,irred}; therefore $M_d^Y(X) = M_d(X)$.
\end{proo}

\smallskip

 \begin{fact}
\label{unique}
   There exists a unique degree $d$ morphism $f:\p^1\to Y_d$ passing
   through three general points of $Y_d$.
 \end{fact}
\begin{proo}
This was proved in \cite{BKT} for (isotropic) Grassmannians.  
For the other cases and $d=2$, this is simply the fact that through 
three general points on a quadric, there exists a unique conic 
(the intersection of the quadric with the plane generated by the 
three points).
Finally, for $E_7/P_7$ and  $d=3$, this is Lemma \ref{e7p7cubiques}. 
\end{proo}

We now prove that $d(x,y)$ classifies the $G$-orbits in $X \times X$.

Let $(x,y) \in X \times X$ and assume $y \not = x$.
Up to the action of $G$, we may assume that
$x$ is the base point and that $y$ is the class of an element $v$ in the
Weyl group. Let $d$ denote the number of occurences of the reflection
$s_\alpha$ in a reduced decomposition of $v$. 
Since the reflections $s_\beta$ for $\beta \not = \alpha$
belong to $P$, we may assume that the quiver $Q_v$ of $v$ has only peaks 
corresponding to the root $\alpha$. So $Q_v$ has only one maximal element,
$(\alpha,d_{max}+1-d)$ (see definition \ref{definition_carquois}); 
therefore $Q_v = \{q \in Q_X : q \preccurlyeq (\alpha,d_{max}+1-d) \}$.

Such subquivers are parametrized by $d \in [0,d_{max}]$, 
so there are at most $d_{max}+1$ orbits in $X \times X$.
For $e \in [0,d_{max}]$, we denote $v_e$ the corresponding element
in $W_X$ (so that $v = v_d$).
Note also that the Schubert variety corresponding to $v_e$ is 
isomorphic with $Y_e$.
Now we will prove that $d=d(x,y)$.
On the one hand, there is a $\p^1$ between the base point in $G/P$ and
the class of $s_\alpha$, so by induction we deduce that $d(x,y) \leq d$.
On the other hand, assume there is a degree $e$ rational curve through
$x$ and $y$; we will show that $e \geq d$. By proposition \ref{existe},
there is an element $\omega \in Y_e$ such that $x$ and $y$ belong to
$Y_\omega$. Since the condition $x,y \in Y_\omega$
is a closed condition on $\omega$, we may furthermore assume that $\omega$ is
$B$-stable. Therefore, $p_e(q_e^{-1}(\omega))$ is the Schubert cell
corresponding to the element $v_e\in W_X$. We thus have $v \leq v_e$ 
for the Bruhat order. This implies $d \leq e$.

\smallskip

We now conclude the proof of Proposition \ref{proposition_yd}. Using
our classification of the couples in $X \times X$, it is easy to check that 
for $x \not = y \in X$ and $d=d(x,y)$, there exists
a unique $\omega \in F_d$ such that $x,y \in Y_\omega$.
Moreover, $(x,y)$ is a generic couple in $Y_\omega \times Y_\omega$.
Therefore, if $z$ is a generic point in $Y_\omega$, by
fact \ref{unique}, $z$ belongs to $Y(x,y)$, so that
$Y(x,y) \supset Y_\omega$. 
On the other hand, if $C$ is a union of rational curves
of degree $d$, then by proposition \ref{existe}, there exists $\beta$
such that $Y_\beta \supset C$; since $x,y \in Y_\beta$ we
deduce $\beta = \omega$ and $Y(x,y) \subset Y_\omega$.
\end{proo}

\begin{rema}
(\i) The varieties $Y(x,y)$ are lines when $d(x,y)=1$, so $F_1$
parametrizes lines in $X$ (therefore $F_1=F$).

(\i\i) The variety $F_2$ parametrizes the maximal quadrics on $X$; the 
nodes defining $F_2$ are such that when we suppress them and keep 
the connected component of the remaining diagram containing the 
node that defines $X$, we get the weighted Dynkin diagram of a
quadric. The maximal quadrics on $X$ are then obtained as Tits shadows
(see \cite{landsberg}).

(\i\i\i) The bounds on the degrees are the easy bounds for the
vanishing of Gromov-Witten invariants. Namely, the degree $d$
Gromov-Witten invariants all vanish as soon as 
$c_1(X)d+\dim(X)> 3\dim(X)$. (Actually, in type A
these easy bounds are even more restrictive in general.)  

\end{rema}

We now deduce from the last two results some equalities of dimensions.
Proposition \ref{existe} implies that
the dimension of the scheme of degree $d$ rational
curves on $X$ equals the dimension of the scheme of degree $d$ rational
curves on $Y_d$ plus the dimension of $F_d$, that is, according to
(\ref{c1(yd)}),
\begin{equation}
\label{relation_dim}
\dim(X)+d \cdot c_1(X) = \dim(F_d)+3\dim(Y_d).
\end{equation}

\begin{rema}
  Conversely, this equality together with the irreducibility of the
  variety ${\bf 
  Mor}_d(\p^1,X)$ of degree $d$ morphisms from $\p^1$ to $X$ (see for
  example \cite{Thomsen} or \cite{irred}) implies that, for any
  morphism $f:\p^1\to X$, there exists $\omega \in F_d$ such that $f$ 
  factors through
  $Y_\omega$ and that for a general morphism $f$,
  there is a finite number of such points $\omega$. It is an easy
  verification that if there are more than one point $\omega$ then $f$ is
  not general. 
\end{rema}

A consequence of Fact \ref{unique}
is that the dimension of the scheme of degree $d$
rational curves on $Y_d$ is $3\dim(Y_d)-3$. Since this dimension can
also be computed from the index of $Y_d$, we get the relation 
\begin{equation}
\label{c1(yd)}
d \cdot c_1(Y_d)=2\dim(Y_d).
\end{equation}

\bigskip

We now derive some combinatorial properties which will be useful
in \cite{cmp2}.
Since the quiver of $Y_d$ is the set of vertices under
$(\alpha,d_{max}+1-d)$, the quiver of $Y_d^*$ is the set of vertices not above
$(\iota(\alpha),d)$, if $\iota$ denotes the Weyl involution of the
simple roots.
Let $\delta(u)$ denote the number of occurences of $s_\alpha$ in a reduced
expression of $u$. We therefore have $X(u) \subset Y_d^*$ \iff
$d \leq \delta(u)$.

From the proof of proposition \ref{couples}, it also follows that there
is a curve of degree $d$ through $1$ and $u$ \iff $d \geq \delta(u)$.

Finally, we relate our integer $\delta(u)$ to an integer defined in
\cite{FW}. Recall the definition \cite[lemma 4.1]{FW} that
two elements $u,v \in W/W_P$ are adjacent if there exists a reflection
$s$ such that $u=vs$. A chain between $u$ and
$v$ is a sequence $u_0,\ldots,u_r$ such that
$u \preccurlyeq u_0,u_r \preccurlyeq v^*$, and each $u_i$ is adjacent
to $u_{i+1}$. Such a chain has a natural degree.
We then consider the following definition, suggested by
\cite[theorem 9.1]{FW}:
\begin{defi}
Let $u,v \in W_X$. Let $\delta(u,v)$ denote the minimal degree of a
chain between $u$ and $v$.
\end{defi}
Note that $\delta$ is symmetric in $u,v$, that $\delta(u,v)=0$ \iff
$u \leq v^*$, and that it is a non-decreasing fonction in the two variables.

\begin{lemm}
For $u \in W_X$, we have $\delta(u,w_X)=\delta(u)$.
\end{lemm}
\begin{proo}
If $u \in W/W_P$, let $x(u)=uP/P \in X=G/P$ denote the corresponding $T$-fixed
point. By theorem 9.1 in \cite{FW}, $\delta(u,w_X)=\delta(w_X,u)$ is the minimal
degree of a curve meeting $\{x(w_X)\}$ and $X(u^*)$. 
Since we can assume that this
curve is $T$-invariant, it will pass through
$x(w_X)$ and $x(v)$ with $v \in W/W_P$ and $v \leq u^*$. Applying the involution
$x \mapsto w_0xw_0w_X$ of $X$, we deduce that this degree is the minimal
degree of a curve through the base point and $v^*$. This minimal degree is
$\delta(v^*)$, and is itself minimal when $v=u^*$, in which case it equals
$\delta(u)$.
\end{proo}

\bigskip

For $x \in X$, we now give a nice description of the subvariety 
$F_{d,x} := \{ w \in F_d : x \in F_\omega \} \subset F_d$.
Although this will not be used
in the sequel, it shows that $F_d$ generalizes very well $F_1$. In fact,
the lines through a fixed point are parametrized by the closed $P$-orbit
in the projectivization of the tangent space at this point, and we will
show that more generally, the set of $Y_\omega$'s through a fixed point are
parametrized by the closed $P$-orbit in the projectivization of the
$d$-th normal space, according to the following definition~:

\begin{defi}
Let $Z \subset \p V$ be a projective variety and let $z \in V- \{0\}$
such that $[z] \in Z$. Let $d$ be an integer. We recall~:
\begin{itemize}
\item
The $d$-th affine tangent space $\widehat{T^d_{[z]}Z} \subset V$
is generated by the $d$-th derivatives at $z$ of curves in the cone over $Z$.
\item
The $d$-th normal space $N^d_{[z]}Z$ is the quotient
$\widehat{T^d_{[z]}Z} / \widehat{T^{d-1}_{[z]}Z}$.
\end{itemize}
\end{defi}

\begin{rema}
(\i) $N^1$ is the tangent space twisted by $-1$.

(\i\i) The $d$-th normal spaces to (co)minuscule homogeneous spaces are given
in \cite[proposition 3.4]{landsberg}. If $L$ denotes a Levi factor of $P$,
they are irreducible $L$-modules.
\end{rema}

\smallskip

For $\omega \in F_d$, let $Y_\omega \subset X$ denote the corresponding
variety. For $x \in X$, let $F_{d,x} \subset F_d$ be the subvariety of 
$\omega$'s such that
$x \in F_\omega$. If $Z \subset \p V$ is any subset, let $\scal Z \subset V$
denote the linear span of its cone in $V$. We have the following~:

\begin{prop}
Let $x \in X$ and $\omega \in F_{d,x}$.
Then $\dim (\scal {Y_\omega} / \widehat{T^{d-1}_xX} ) = 1$.
Moreover, the morphism
$$
\begin{array}{rcl}
F_{d,x} & \rightarrow & \p N^d_xX\\
\omega  & \mapsto     & [ \scal {Y_\omega} ]
\end{array}
$$
is a closed immersion, with image the closed $L$-orbit in
$\p N^d_xX$.
\end{prop}
\begin{proo}
Let $x \in X$ be the base point, and let $L$ be a Levi factor of $P$.
By proposition \ref{proposition_yd}, the subvarieties 
$F_{d,x} \subset F_d$ parametrize the elements $\omega$
in $F_d$ which are incident to $x$ in the sense of Tits geometries
(namely the stabilizer of $\omega$ and that of $x$ intersect along
a parabolic subgroup), so $F_{d,x}$ is
homogeneous under $L$. For example, if $X=\G(p,n)$, then
$F_{d,x} \simeq \G(p-d,p) \times \G(d,n-p)$.

Therefore, to check that 
$\dim (\scal {Y_\omega} / \widehat{T^{d-1}_xX} ) = 1$, it is enough
to consider one particular example for $\omega$; we leave this to the
reader, as well as the fact that the corresponding class
$[ \scal {Y_\omega} ] \in \p N^d_xX$ belongs to the closed
$L$-orbit $\cal O$. For example, when $X=\G(p,n)$, if
$x$ denotes the linear space generated by
$e_1,\ldots,e_p$ and $\omega$ the flag
$(\scal{e_1,\ldots,e_{p-d}},\scal{e_1,\ldots,e_{p+d}})$, then
any element in $Y_\omega$ is the linear span of
$e_1,\ldots,e_{p-d},f_1,\ldots,f_d$, with $f_i \in \scal{e_1,\ldots,e_{p+d}}$.
Therefore, in the Pl\"ucker coordinates, this element is
equivalent to a multiple of
$e_1\wedge \ldots \wedge e_d \wedge e_{p+1} \wedge \ldots \wedge e_{p+d}$
modulo $\widehat{T^{d-1}_xX}$.

We therefore have
an $L$-equivariant morphism $F_{d,x} \rightarrow {\cal O}$, which must be
an isomorphism for example because the two varieties have the same
dimension.
\end{proo}

\subsectionplus{Higher Gromov-Witten invariants}

\label{subsection_higher}

In this section, we deduce from the preceeding geometric results
a way of computing higher degree 
Gromov-Witten invariants.

We keep the previous notations. As in the case of lines, 
the varieties $Z_d$ and $Y_d$ are 
(co)minuscule homogeneous variety. In particular they are endowed
with well defined quivers $Q_{Z_d}$ and $Q_{Y_d}$. The variety $Y_d$ 
is a Schubert subvariety of $X$ and can thus be written as
$X(w_{Y_d})$ for some $w_{Y_d}\in W_X$. We denote its Poincar{\'e} dual by
$X(w_{Y_d}^*)$ or $Y_d^*$.

We define as is the previous section
the quivers $Q_{I_d}$ and $Q_{F_d}$ of $I_d$ and $F_d$,  by
adding $Q_{Z_d}$ at the end of $Q_{X}$ (resp. by adding $Q_{Z_d}$ at
the end of $i_X(Q_{Y_d})$). These quivers
correspond to the particular reduced decomposition obtained 
through the formula
$$w_{I_d}=w_{X}w_{Z_d}=w_{Y_d}w_{F_d}.$$

\begin{example}
\label{exemple-grass1}
Let $X=\G(p,n)$. We describe the quivers $Q_{X}$, $Q_{Y_d}$,
$Q_{Z_d}$ and $Q_{F_d}$ and their different inclusions
for $d\leq\min(p,n-p)$. We have already seen that $Q_X$ is
a $p\times(n-p)$ rectangle. We draw it as on the left picture.

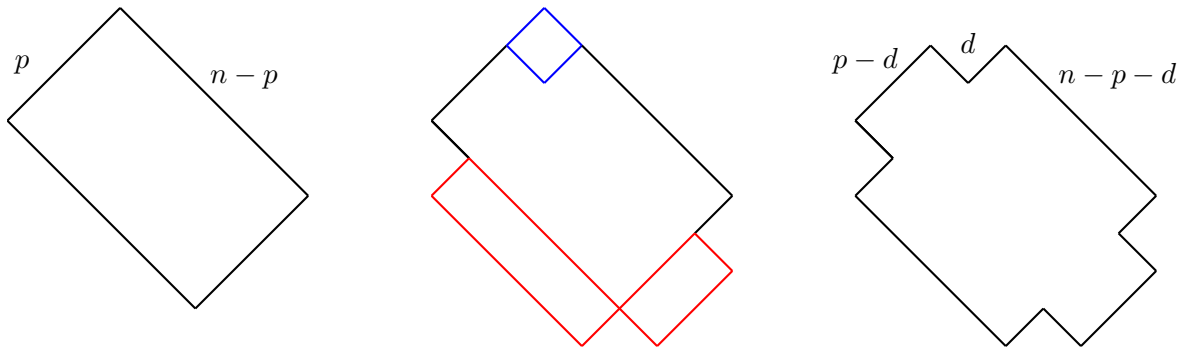
\psset{xunit=0.5cm}
\psset{yunit=0.5cm}
\centerline{\begin{pspicture*}(-3,-8)(6,3)
\psline(0,2)(-3,-1)
\psline(-3,-1)(2,-6)
\psline(2,-6)(5,-3)
\psline(5,-3)(0,2)
\put(-1.4,0.2){$p$}
\put(1.2,0){$n-p$}
\end{pspicture*}
\hskip 1 cm
\begin{pspicture*}(-3,-8)(6,3)
\psline(-1,1)(-3,-1)
\psline[linecolor=blue](0,2)(-1,1)
\psline[linecolor=blue](0,0)(-1,1)
\psline[linecolor=blue](0,0)(1,1)
\psline(-3,-1)(-2,-2)
\psline(4,-4)(5,-3)
\psline[linecolor=red](2,-6)(3,-7)
\psline(-2,-2)(-3,-1)
\psline[linecolor=blue](1,1)(0,2)
\psline(5,-3)(1,1)
\psline[linecolor=red](2,-6)(-2,-2)
\psline[linecolor=red](-2,-2)(-3,-3)
\psline[linecolor=red](-3,-3)(1,-7)
\psline[linecolor=red](1,-7)(2,-6)
\psline[linecolor=red](4,-4)(2,-6)
\psline[linecolor=red](5,-5)(4,-4)
\psline[linecolor=red](5,-5)(3,-7)
\end{pspicture*}
\hskip 0.5 cm
\begin{pspicture*}(-4,-8)(6,3)
\psline(-1,1)(-3,-1)
\psline(0,0)(-1,1)
\psline(0,0)(1,1)
\psline(-3,-1)(-2,-2)
\psline(4,-4)(5,-3)
\psline(2,-6)(3,-7)
\psline(-2,-2)(-3,-1)
\psline(5,-3)(1,1)
\psline(-2,-2)(-3,-3)
\psline(-3,-3)(1,-7)
\psline(1,-7)(2,-6)
\psline(5,-5)(4,-4)
\psline(5,-5)(3,-7)
\put(-1.8,0.2){$p-d$}
\put(-.1,0.4){$d$}
\put(1.2,0){$n-p-d$}
\end{pspicture*}}
On the middle picture we have drawn $Q_{I_d}$ and
inside it,  $Q_{Z_d}$ in red  and the complement of $i_X(Q_{w_{Y_d}})$
in $Q_X$ in blue. 
The quiver on the right is $Q_{F_d}$.

The similar pictures in the isotropic cases are (we denote
$\diamond = 2d, N=n-1$ in the quadratic case and $\diamond = d,N=n$ in the
symplectic case):

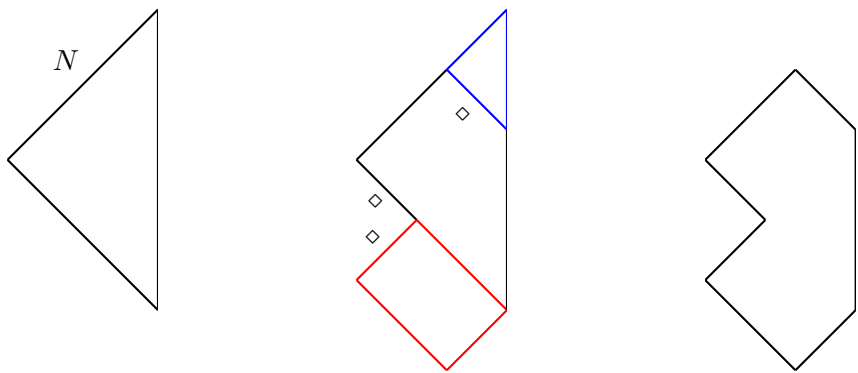
\psset{unit=.4cm}

\centerline{
\begin{pspicture*}(0,-3)(5,13)
\put(1.5,8){$N$}
\psline(5,10)(0,5)
\psline(0,5)(5,0)
\psline(5,10)(5,0)
\end{pspicture*}
\hskip 2.5cm
\begin{pspicture*}(0,-3)(5,13)
\psline[linecolor=blue](5,10)(3,8)
\psline(3,8)(0,5)
\psline(0,5)(2,3)
\put(0.4,3.4){$\diamond$}
\psline[linecolor=red](2,3)(5,0)
\psline[linecolor=blue](5,10)(5,6)
\psline(5,6)(5,0)
\psline[linecolor=blue](3,8)(5,6)
\put(3.3,6.3){$\diamond$}
\psline[linecolor=red](5,0)(3,-2)
\psline[linecolor=red](3,-2)(0,1)
\psline[linecolor=red](0,1)(2,3)
\put(0.3,2.2){$\diamond$}
\end{pspicture*}
\hskip 2.5cm
\begin{pspicture*}(0,-3)(5,13)
\psline(3,8)(0,5)
\psline(0,5)(2,3)
\psline(5,6)(5,0)
\psline(3,8)(5,6)
\psline(5,0)(3,-2)
\psline(3,-2)(0,1)
\psline(0,1)(2,3)
\end{pspicture*}
}

\end{example}

\psset{unit=1cm}

Let $X(w)$ be a Schubert subvariety in $X$. Consider the Schubert
subvariety $F_d(\wh)=q_d(p_d^{-1}(X(w))$ of $F_d$. The quiver of
$F_d(\wh)$ in $Q_{F_d}$ is obtained as follows: attach $Q_{Z_d}$ to
the bottom end of $Q_w\cap i_X(Q_{w_{Y_d}^*})$. In
particular, we have the inequality
\begin{equation}
\label{codimFd(wh)}
\codim_{F_d}(F_d(\wh))\geq\codim_X(X(w))-\dim(Y_d),
\end{equation}
with equality if and only if $X(w)\subset X(w_{Y_d}^*)$ .
(This generalizes inequalities obtained in \cite{BKT} for (isotropic)
Grassmannians.)

Now let $f:\p^1\to X$ be a degree $d$ morphism such that
$f(\p^1)$ meets $X(w)$. Then there exists $y\in F_d$ such that
$p_d(q_d^{-1}(y))$ meets $X(w)$, or equivalently, such that $y\in
q(p^{-1}(X(w)))=F_d(\wh)$. This is the key point to compute 
degree $d$ Gromov-Witten invariants on $X$ in terms of classical
invariants on $F_d$. 

\begin{lemm}
Let $X(u)$, $X(v)$ and $X(w)$ be three Schubert subvarieties of $X$. 
Suppose that
$$\codim(X(u))+\codim(X(v))+\codim(X(w)) = \dim(X)+d\cdot c_1(X).$$
Then for $g$, $g'$ and $g''$ three general elements in $G$, the
intersection $g\cdot F(\uh)\cap g'\cdot F(\vh)\cap g''\cdot F(\wh)$ is
a finite set of reduced points. Moreover this finite set is empty 
if one of the Schubert varieties is not contained in $X(w_{Y_d}^*)$. 

Let $y$ be a point in this intersection. Then the variety 
$p_d(q_d^{-1}(y))$ meets each of $g\cdot X(u)$, $g'\cdot X(v)$ and
$g''\cdot X(w)$ in a unique point and these points are in general
position in $p_d(q_d^{-1}(y))$.
\end{lemm}

\begin{proo}
Remark that the codimension condition, inequality (\ref{codimFd(wh)})
and equality (\ref{relation_dim}) imply that
$$
\begin{array}{l}
 \codim(F_d(\uh))+\codim(F_d(\vh))+\codim(F_d(\wh))\\
\hspace{3cm} \geq \; \codim(X_d(u)) + \codim(X_d(v)) + \codim(X_d(w))
- 3 \dim(Y_d)\\
\hspace{3cm} = \;  \dim(F_d),
\end{array}
$$
with equality if and only if the three Schubert varieties are
contained in $X(w_{Y_d}^*)$. Actually, this is true except for
Grassmannians and for $d>\min(p,n-p)$, in which case the previous 
inequality is always strict. The first part of the
lemma is thus implied by Bertini's theorem (see \cite{kleiman}). 

Furthermore, by Bertini again, we may assume that any $y$ in the 
intersection is general in $g\cdot F_d(\uh)$, $g'\cdot F_d(\vh)$ and 
$g''\cdot F_d(\wh)$. In particular, by lemma \ref{intersection} applied to
$v=w_{Y_d}$,  the variety $p_d(q_d^{-1}(y))$ meets each of $g\cdot
X(u)$, $g'\cdot X(v)$ and $g''\cdot X(w)$ transversely in a unique
point. Finally,
the stabiliser of $y$ acts transitively on $p_d(q_d^{-1}(y))$ and by
modifying $g$, $g'$ and $g''$ by elements in this stabiliser me may
assume that the points are in general position in $p_d(q_d^{-1}(y))$. 
\end{proo}

\begin{coro}
\label{gw-d}
Let $X(u)$, $X(v)$ and $X(w)$ be three Schubert subvarieties of $X$. 
Suppose that the sum of
their codimensions is $\dim(X)+d\cdot c_1(X)$. Then
$$I_d([X(u)],[X(v)],[X(w)])=I_0([F_d(\uh)],[F_d(\vh)],[F_d(\wh)]).$$
In particular, this invariant vanishes as soon as one of the three
Schubert varieties is not contained in $X(w_{Y_d}^*)$.
\end{coro}

\begin{proo}
The image of any morphism $f$ counting in the invariant
$I_d([X(u)],[X(v)],[X(w)])$ is contained in a variety
$p_d(q_d^{-1}(y))$ with $y\in g\cdot F_d(\uh)\cap g'\cdot
F_d(\vh)\cap g''\cdot F_d(\wh)$ for general elements $g$, $g'$ and
$g''$ in $G$. The preceding lemma implies that this intersection is
either empty or a finite number of reduced points. Given such a
$y$, the morphism $f$ has to pass through three fixed general points in
$p_d(q_d^{-1}(y))$, and by fact \ref{unique}, there exists a
unique such morphism. 

In particular, for $X=\G(p,n)$ and  $d>\min(p,n-p)$, all degree $d$
Gromov-Witten invariants vanish.
\end{proo}

\begin{rema}
  This result has been proved in \cite{BKT} for (isotropic)
  Grassmannians through  a case by case analysis. The vanishing
  condition generalizes the conditions of \cite{yong} for ordinary
  Grassmannians and of \cite{BKT} for isotropic ones. 
\end{rema}

\begin{example}
\label{exemple-grass2}
Let $X=\G(p,n)$. The quiver $Q_w$ of a Schubert subvariety $X(w)$ 
has the following form (recall that it is the
complement of the partition associated to $w$ inside  
the rectangle $p\times (n-p)$):

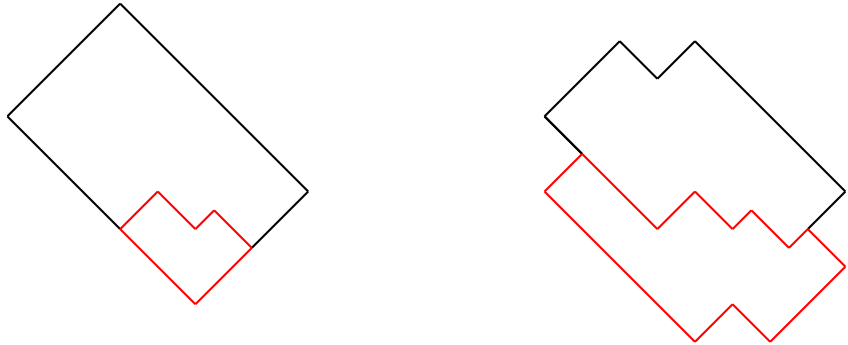
\psset{xunit=0.5cm}
\psset{yunit=0.5cm}
\centerline{\begin{pspicture*}(-3,-8)(6,3)
\psline(0,2)(-3,-1)
\psline(-3,-1)(0,-4)
\psline(3.5,-4.5)(5,-3)
\psline(5,-3)(0,2)
\psline[linecolor=red](0,-4)(2,-6)
\psline[linecolor=red](2,-6)(3.5,-4.5)
\psline[linecolor=red](2.5,-3.5)(3.5,-4.5)
\psline[linecolor=red](2,-4)(2.5,-3.5)
\psline[linecolor=red](1,-3)(2,-4)
\psline[linecolor=red](1,-3)(0,-4)
\end{pspicture*}
\hskip 2 cm
\begin{pspicture*}(-4,-8)(6,3)
\psline(-1,1)(-3,-1)
\psline(0,0)(-1,1)
\psline(0,0)(1,1)
\psline(-3,-1)(-2,-2)
\psline(4,-4)(5,-3)
\psline[linecolor=red](3.5,-4.5)(4,-4)
\psline(-2,-2)(-3,-1)
\psline(5,-3)(1,1)
\psline[linecolor=red](0,-4)(-2,-2)
\psline[linecolor=red](2.5,-3.5)(3.5,-4.5)
\psline[linecolor=red](2,-4)(2.5,-3.5)
\psline[linecolor=red](2,-6)(3,-7)
\psline[linecolor=red](1,-3)(2,-4)
\psline[linecolor=red](1,-3)(0,-4)
\psline[linecolor=red](-2,-2)(-3,-3)
\psline[linecolor=red](-3,-3)(1,-7)
\psline[linecolor=red](1,-7)(2,-6)
\psline[linecolor=red](5,-5)(4,-4)
\psline[linecolor=red](5,-5)(3,-7)
\end{pspicture*}}
On the left picture we have drawn the quiver $Q_{X}$ and
inside it in red the Schubert subquiver $Q_{w}$. On the right 
we have the quiver $Q_{F_d}$ and inside in red the subquiver of 
$F_d(\wh)$.
\end{example}



\sectionplus{The quantum Chevalley formula 
and a higher Poincar{\'e} duality}

In this section we give, thanks again to the
combinatorics of quivers, a simple
combinatorial version of the quantum Chevalley formula (proposition
\ref{chevalley}). We also describe what we call a higher Poincar{\'e}
duality (proposition \ref{q-poincare}): a duality on Schubert classes
defined in terms
of  degree $d$ Gromov-Witten invariants.

\subsectionplus{Quivers and the quantum Chevalley formula}

A general quantum Chevalley formula has been obtained by W. Fulton and
C. Woodward in \cite{FW}, following ideas of D. Peterson. In this
subsection we recover this formula for any (co)minuscule homogeneous
variety $X$, with a very simple combinatorial description in terms of
quivers.

Indeed, since the codimension one Schubert subvariety $H$ is certainly 
not contained in $X(w_{Y_d}^*)$ for $d\ge 2$, the 
vanishing criterion of corollary \ref{gw-d} ensures that the quantum
product with $H$ only involves Gromov-Witten invariants degree zero and
one.

Let us first give a few more notations to describe our quantum Chevalley
formula. If $i$ is a peak of $Q_w$, we denote by $Q_{w(i)}$ the full
subquiver of $Q_w$ obtained by removing the vertex $i$ and by $w(i)$
the corresponding element in $W$. 
Embed the quivers $Q_{X}$, $Q_{F}$ and $Q_{Z}$ in $Q_{I}$ as
explained in subsection 
\ref{subsection_fano}. If $X(w)$ is a Schubert subvariety
of $X$, consider the Schubert variety $F(i_F(\wh))$ in
$F$. If it exists, denote by $w_q$ the element in $W$ such that
$F(i_F(\wh))=F(\wh_q)$.

\begin{prop}\label{chevalley}
For any Schubert subvariety $X(w)$ of $X$, we have
$$[X(w)]*[H]=\sum_{i\in p(Q_w)}[X(w(i))]+q [X(i_X(w_q))].$$
\end{prop}

\begin{proo}
  The degree zero part of the right hand side is the classical
product $[X(w)]\cdot [H]$, for which we have just reformulated the 
classical Chevalley formula in terms of quivers. 

For the degree one part, we need to compute the Gromov-Witten 
invariants $I_1([H],[X(w)],[X(v)])$ for all $v\in W_X$. 
Since $q(p^{-1}(H))=F$, this amounts by corollary
\ref{XF_degre1} to compute
  $I_0([F],[F(\wh)],F(\vh)])$. By Poincar{\'e} duality on $F$ 
this invariant is  zero unless $\vh=i_F(\wh)$, in which case it is
equal to one. But this exactly means that $w_q$ exists and
is equal to $v$. 
\end{proo}

\begin{defi}
\label{1-poincare}
We will say that $X(w)$ and $X(w_q)$ are {\it 1-Poincar{\'e} dual}.
\end{defi}

\begin{example}
Let $X$, $X(w)$ and $X(v)$ as in example \ref{exemplebis} (\i\i). Then
the quivers of $F(i_F(\wh))$ and $F(i_F(\vh))$ are the following:

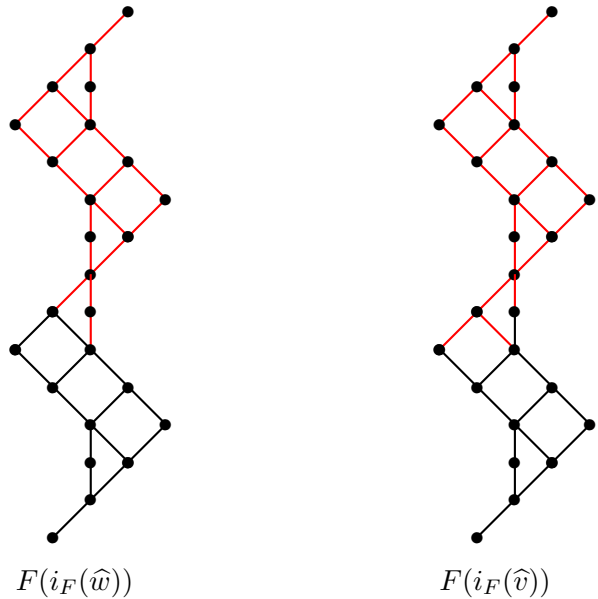
\psset{xunit=0.5cm}
\psset{yunit=0.5cm}
\centerline{\begin{pspicture*}(-5,-11)(2,6.3)
\psline[linecolor=red](-0,0)(-1,1)
\psline[linecolor=red](-2,4)(-2,3)
\psline[linecolor=red](-2,2)(-1,1)
\psline[linecolor=red](-2,2)(-3,3)
\psline[linecolor=red](-1,-1)(-4,2)
\psline[linecolor=red](-2,2)(-2,3)
\psline[linecolor=red](-3,3)(-4,2)
\psline[linecolor=red](-2,2)(-3,1)
\psline[linecolor=red](-1,1)(-2,0)
\psline[linecolor=red](-0,0)(-3,-3)
\psline(-3,-3)(-4,-4)
\psline[linecolor=red](-2,-1)(-2,0)
\psline[linecolor=red](-2,-1)(-2,-2)
\psline[linecolor=red](-1,5)(-3,3)
\put(-0.1,-0.1){$\bullet$}
\put(-0.6,0.4){$\bullet$}
\put(-1.1,0.9){$\bullet$}
\put(-1.6,1.4){$\bullet$}
\put(-1.1,-0.1){$\bullet$}
\put(-1.6,0.4){$\bullet$}
\put(-2.1,0.9){$\bullet$}
\put(-0.6,-0.6){$\bullet$}
\put(-1.1,1.4){$\bullet$}
\put(-1.1,-0.6){$\bullet$}
\put(-2.1,-2.1){$\bullet$}
\put(-1.6,-1.6){$\bullet$}
\put(-1.1,-1.1){$\bullet$}
\put(-0.6,-0.6){$\bullet$}
\put(-0.6,2.4){$\bullet$}
\put(-1.1,1.9){$\bullet$}
\psline[linecolor=red](-2,-2)(-2,-4)
\psline(-3,-3)(0,-6)
\psline(-4,-4)(-1,-7)
\psline(0,-6)(-3,-9)
\psline(-2,-6)(-2,-8)
\psline(-2,-4)(-3,-5)
\psline(-1,-5)(-2,-6)
\put(-0.1,-3.1){$\bullet$}
\put(-0.6,-2.6){$\bullet$}
\put(-1.1,-2.1){$\bullet$}
\put(-1.6,-1.6){$\bullet$}
\put(-1.1,-3.1){$\bullet$}
\put(-1.6,-2.6){$\bullet$}
\put(-2.1,-2.1){$\bullet$}
\put(-0.6,-3.6){$\bullet$}
\put(-1.1,-1.6){$\bullet$}
\put(-1.1,-3.6){$\bullet$}
\put(-1.6,-4.6){$\bullet$}
\put(-1.1,-4.1){$\bullet$}
\put(-0.6,-3.6){$\bullet$}
\put(-2,-5.2){$F(i_F(\wh))$}
\end{pspicture*}
\hskip 2 cm
\begin{pspicture*}(-5,-11)(2,6.3)
\psline[linecolor=red](-0,0)(-1,1)
\psline[linecolor=red](-2,4)(-2,3)
\psline[linecolor=red](-2,2)(-1,1)
\psline[linecolor=red](-2,2)(-3,3)
\psline[linecolor=red](-1,-1)(-4,2)
\psline[linecolor=red](-2,2)(-2,3)
\psline[linecolor=red](-3,3)(-4,2)
\psline[linecolor=red](-2,2)(-3,1)
\psline[linecolor=red](-1,1)(-2,0)
\psline[linecolor=red](-0,0)(-1,-1)
\psline[linecolor=red](-1,-1)(-4,-4)
\psline[linecolor=red](-2,0)(-2,-2)
\psline[linecolor=red](-1,5)(-3,3)
\put(-0.1,-0.1){$\bullet$}
\put(-0.6,0.4){$\bullet$}
\put(-1.1,0.9){$\bullet$}
\put(-1.6,1.4){$\bullet$}
\put(-1.1,-0.1){$\bullet$}
\put(-1.6,0.4){$\bullet$}
\put(-2.1,0.9){$\bullet$}
\put(-0.6,-0.6){$\bullet$}
\put(-1.1,1.4){$\bullet$}
\put(-1.1,-0.6){$\bullet$}
\put(-2.1,-2.1){$\bullet$}
\put(-1.6,-1.6){$\bullet$}
\put(-1.1,-1.1){$\bullet$}
\put(-0.6,-0.6){$\bullet$}
\put(-0.6,2.4){$\bullet$}
\put(-1.1,1.9){$\bullet$}
\psline(-2,-6)(-2,-8)
\psline(-2,-4)(0,-6)
\psline[linecolor=red](-2,-4)(-3,-3)
\psline(-4,-4)(-1,-7)
\psline(0,-6)(-3,-9)
\psline[linecolor=red](-2,-3)(-2,-2)
\psline(-2,-3)(-2,-4)
\psline(-2,-4)(-3,-5)
\psline(-1,-5)(-2,-6)
\put(-0.1,-3.1){$\bullet$}
\put(-0.6,-2.6){$\bullet$}
\put(-1.1,-2.1){$\bullet$}
\put(-1.6,-1.6){$\bullet$}
\put(-1.1,-3.1){$\bullet$}
\put(-1.6,-2.6){$\bullet$}
\put(-2.1,-2.1){$\bullet$}
\put(-0.6,-3.6){$\bullet$}
\put(-1.1,-1.6){$\bullet$}
\put(-1.1,-3.6){$\bullet$}
\put(-1.6,-4.6){$\bullet$}
\put(-1.1,-4.1){$\bullet$}
\put(-0.6,-3.6){$\bullet$}
\put(-2,-5.2){$F(i_F(\vh))$}
\end{pspicture*}}

\noindent
and in particular $F(i_F(\wh))$ cannot be of the form $F(\uh)$ for
some $u$. On the contrary for $u$ with the following quiver, 

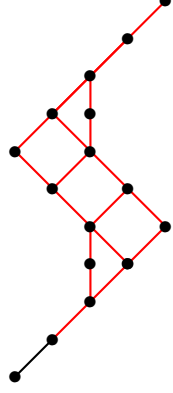
\psset{xunit=0.5cm}
\psset{yunit=0.5cm}
\centerline{\begin{pspicture*}(-5,-4.3)(2,6.3)
\psline[linecolor=red](-0,0)(-1,1)
\psline[linecolor=red](-0,6)(-3,3)
\psline[linecolor=red](-2,4)(-2,3)
\psline[linecolor=red](-2,2)(-1,1)
\psline[linecolor=red](-2,2)(-3,3)
\psline[linecolor=red](-1,-1)(-4,2)
\psline[linecolor=red](-2,2)(-2,3)
\psline[linecolor=red](-3,3)(-4,2)
\psline[linecolor=red](-2,2)(-3,1)
\psline[linecolor=red](-1,1)(-2,0)
\psline[linecolor=red](-0,0)(-1,-1)
\psline[linecolor=red](-1,-1)(-3,-3)
\psline(-3,-3)(-4,-4)
\psline[linecolor=red](-2,0)(-2,-2)
\psline[linecolor=red](-1,5)(-3,3)
\put(-0.1,-0.1){$\bullet$}
\put(-0.6,0.4){$\bullet$}
\put(-1.1,0.9){$\bullet$}
\put(-1.6,1.4){$\bullet$}
\put(-1.1,-0.1){$\bullet$}
\put(-1.6,0.4){$\bullet$}
\put(-2.1,0.9){$\bullet$}
\put(-0.6,-0.6){$\bullet$}
\put(-1.1,1.4){$\bullet$}
\put(-1.1,-0.6){$\bullet$}
\put(-2.1,-2.1){$\bullet$}
\put(-1.6,-1.6){$\bullet$}
\put(-1.1,-1.1){$\bullet$}
\put(-0.6,-0.6){$\bullet$}
\put(-0.1,2.9){$\bullet$}
\put(-0.6,2.4){$\bullet$}
\put(-1.1,1.9){$\bullet$}
\end{pspicture*}}

\noindent
we have $F(i_F(\vh))=F(\uh)$. In particular we obtain with the
notations of subsection \ref{subsubsection_cayley}:
$$H*\sigma'_{12}=\sigma_{13}\qquad and\qquad
H*\sigma''_{12}=\sigma_{13}+q\sigma_1.$$
\end{example}

\subsectionplus{Higher quantum Poincar{\'e} duality}

\label{subsection_higher_quantum}

Poincar{\'e} duality can be reformutated as follows: there exists an
involution of $W_X$, 
given by $v\mapsto v^* = w_0vw_0w_X$, such that
$$I_0([X],[X(v)],[X(w)])= \delta_{w,v^*}.$$
We
have seen that for degree one invariants, the hyperplane class $[H]$ 
plays the role of $[X]$: $X(v)$ and $X(w)$ are 1-Poincar{\'e} dual
(see definition \ref{1-poincare}) if and 
only if $$I_1([H],[X(v)],[X(w)])=1.$$

More generally, the class $[Y_d^*]$ will play the role of $[X]$ for
degree $d$ Gromov-Witten invariants. We will define an involution
$v \mapsto v_{q^d}$ of a subset of $W_X$, with a simple combinatorial
interpretation, and such that
$$I_d([Y_d^*],[X(v)],[X(w)])=\delta_{w,v_{q^d}}$$
(with the understanding that if $v_{q^d}$ is not defined, then the invariant
is zero).

Before giving a precise definition of $v_{q^d}$, let us describe
Poincar{\'e} duality on $F_d$. As for $X=F_0$ or $F=F_1$, compiling
reduced expressions  
$w_X=s_{\b_1}\cdots s_{\b_N}$ and $w_{Z_d}=s_{\b'_1}\cdots
s_{\b'_{M_d}}$, where $M_d=\dim(Z_d)$,  we obtain the reduced expression
$$w_{F_d} = s_{\b_2}\cdots s_{\b_{N}} s_{\b'_1}\cdots s_{\b'_{M_d}}.$$
Modulo commutation relations, this expression is symmetric, that is, 
of the form $s_{\ga_1}\cdots s_{\ga_{R_d}}$ with $R_d=\dim(F_d)$
and $i(\ga_k)=\ga_{R_d+1-k}$. The associated quiver $Q_{w_{F_d}}$ is
symmetric and we denote by $i_{F_d}$ the induced involution on
subquivers. The same proof as for proposition \ref{dualite} gives the
following result:

\begin{prop}
\label{poinc-f_d}
   Let $F_d(w)$ be a Schubert subvariety of $F_d$ such that
   $w=s_{\ga_{k+1}}\cdots s_{\ga_{R_d}}$. Then the classes $[F_d(w)]$
   and $[F_d(i_{F_d}(w))]$ are Poincar{\'e} dual.
\end{prop} 

\begin{rema}
(\i) Beware  that not all Schubert varieties $F_d(w), w\in W_{F_d}$
satisfy the hypothesis of the proposition. This is
because $F_d$ is not minuscule and  in consequence there may be 
braid relations. However, all Schubert varieties $F_d(w)$ associated to a
Schubert subquiver $Q_{F_d}(w)$ of $Q_{F_d}$ satisfy the property.

(\i\i) We will denote by $F_d(u^\star)$ the Poincar{\'e} dual of $F_d(u)$.
\end{rema}

The quiver $Q_{F_d}$ contains $Q_{w_{Y_d}^*}=Q_{Y_d^*}$ and is
symmetric. We denote by $i_{F_d}$ the associated involution. The
subquiver 
$$Q_{w_{T_d}}:=Q_{Y_d^*}\cap i_{F_d}(Q_{Y_d^*})\subset Q_{Y_d^*}\subset Q_X$$
is a symmetric Schubert subquiver. We let $T_d=X(w_{T_d})$ and 
denote by $i_{T_d}$ the involution on $Q_{T_d}$.

\smallskip
The varieties $T_d$ are given by the following table. Observe that 
they are always smooth, and that the vertices of $Q_{T_{d-1}}$ are those
under the vertex $(\iota(\alpha),d)$, where $\iota$ is the Weyl
involution of the simple roots and $\alpha$ is the root defining $X$.

$$\begin{array}{ccc}
 X & d & T_d  \\
 & &  \\
\G(p,n) & d\leq\min(p,n-p) & \G(p-d,n-2d)\\
\G_{\omega}(n,2n) & d\leq n & \G_{\omega}(n-d,2n-2d)\\
\G_{Q}(n,2n) & d\leq \frac{n}{2} & \G_{Q}(n-2d,2n-4d)\\
\QQ^n & d=1 &  \PP^1\\
 & d=2 & \{{\rm pt}\} \\
E_6/P_1 & d=1 & \PP^5 \\
 & d=2 & \{{\rm pt}\} \\
E_7/P_7 & d=1 & \QQ^{10}\\
 & d=2 & \PP^1\\
 & d=3 & \{{\rm pt}\}
\end{array}$$

\smallskip

\begin{defi}
  The application $v\mapsto v_{q^d}$ is defined for all $v\in W_X$
  such that $Q_v$ is contained in $Q_{w_{T_d}}$, or
  equivalently $X(v)\subset T_d$, by 
$$Q_{v_{q^d}}=i_{F_d}(Q_v)\cap
  Q_X=i_{T_d}(Q_v\cap Q_{T_d}).$$ 
\end{defi}

Otherwise said, the map $v\mapsto v_{q^d}$ is given by Poincar{\'e}
duality inside $T_d$.

\begin{prop}
\label{q-poincare}
  The Gromov-Witten invariant
  $I_d([X(w_{Y_d}^*)],[X(v)],[X(w)])$ vanishes unless
$w=v_{q^d}$.  
In that  case the invariant is equal to one.
\end{prop}

\begin{proo}
  The proof is similar to that of proposition \ref{chevalley}. 
From corollary  \ref{gw-d}, we know that 
$$I_d([X(w_{Y_d}^*)],[X(v)],[X(w)]) = 
  I_0([F_d(\widehat{w_{Y_d}^*})],[F_d(\vh)],[F_d(\wh)]).$$ 
But $F_d(\widehat{w_{Y_d}^*})=F_d$, so this invariant is trivial unless 
$[F_d(\vh)]$ and  $[F_d(\wh)]$ are Poincar{\'e} dual in $F_d$. But
  proposition \ref{poinc-f_d} applies to $F_d(\vh)$ and the invariant
  vanishes unless the quivers of $F_d(\vh)$ and $F_d(\wh)$ are
  symetric under $i_{F_d}$ and in that case the invariant equals
  one. This is equivalent to $w=v_{q^d}$.
\end{proo}

Before dealing with examples, let us state the following lemma
generalizing Lemma \ref{cas-gen1} for $d\geq2$; the proof is the
same. It will be useful
together with Lemma \ref{annul} to prove the vanishing of some
Gromov-Witten invariants.

\begin{lemm}
  \label{cas-gen2}
Let $u$ and $v$ in $W_{F_d}$ such that $F_d(u)$ and $F_d(v)$ are
represented by Schubert subquivers $Q_{F_d}(u)$ and $Q_{F_d}(v)$ of
$Q_{F_d}$.

(\i) If $Q_{F_d}(u)\subset Q_{F_d}(v)$ then ${F_d}(u)\subset {F_d}(v)$
(see Fact \ref{cas-gen} (\i)).

(\i\i) Conversely, if ${F_d}(u)\subset {F_d}(v)$, then we have the
inclusion $Q_{F_d}(u)\cap i_{F_d}(Q_{Z_d})\subset Q_{F_d}(v)\cap
i_{F_d}(Q_{Z_d})$.
\end{lemm}

\begin{example}
(\i) Suppose again that $X=\G(p,n)$ is a Grassmannian. 
We give on the left picture the
quiver $Q_{F_d}$ with inside it in blue the quiver 
$i_{F_d}(Q_{Z_d})$ and in red the quiver $Q_{Y_d}$ one can add to
$i_{F_d}(Q_{Z_d})$ to get $Q_{X}$.

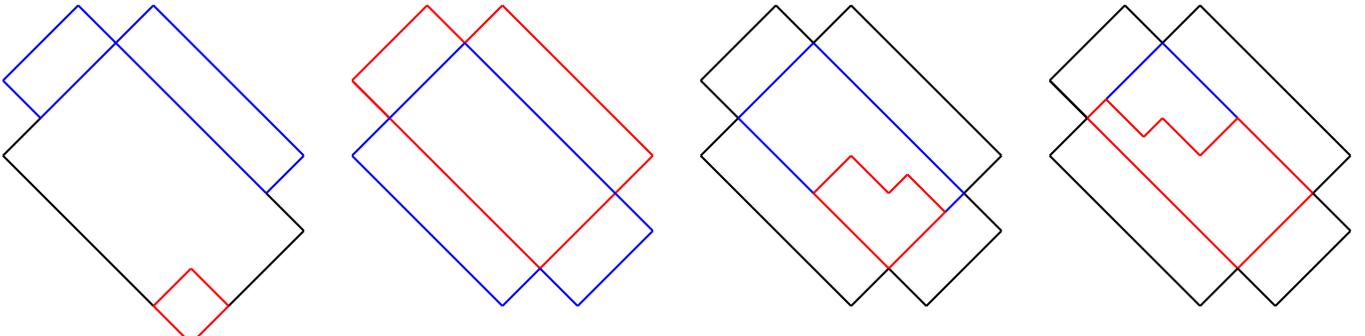
\psset{xunit=0.5cm}
\psset{yunit=0.5cm}
\centerline{\begin{pspicture*}(-3,-8)(6,3)
\psline[linecolor=blue](-1,1)(-3,-1)
\psline[linecolor=blue](0,0)(-2,-2)
\psline[linecolor=blue](0,0)(-1,1)
\psline[linecolor=blue](-3,-1)(-2,-2)
\psline[linecolor=blue](4,-4)(5,-3)
\psline[linecolor=blue](1,1)(0,0)
\psline[linecolor=blue](5,-3)(1,1)
\psline[linecolor=blue](4,-4)(0,0)
\psline(-2,-2)(-3,-3)
\psline(-3,-3)(1,-7)
\psline[linecolor=red](1,-7)(2,-6)
\psline[linecolor=red](1,-7)(2,-8)
\psline[linecolor=red](3,-7)(2,-8)
\psline(5,-5)(4,-4)
\psline(5,-5)(3,-7)
\psline[linecolor=red](2,-6)(3,-7)
\end{pspicture*}
\begin{pspicture*}(-3,-8)(6,3)
\psline[linecolor=red](-1,1)(-3,-1)
\psline[linecolor=red](0,0)(-1,1)
\psline[linecolor=red](-3,-1)(2,-6)
\psline[linecolor=red](2,-6)(5,-3)
\psline[linecolor=red](-2,-2)(-3,-1)
\psline[linecolor=red](1,1)(0,0)
\psline[linecolor=red](5,-3)(1,1)
\psline[linecolor=blue](0,0)(-3,-3)
\psline[linecolor=blue](-3,-3)(1,-7)
\psline[linecolor=blue](1,-7)(2,-6)
\psline[linecolor=blue](3,-7)(2,-6)
\psline[linecolor=blue](5,-5)(0,0)
\psline[linecolor=blue](5,-5)(3,-7)
\end{pspicture*}
\begin{pspicture*}(-3,-8)(6,3)
\psline(-1,1)(-3,-1)
\psline(0,0)(-1,1)
\psline(-3,-1)(-2,-2)
\psline(4,-4)(5,-3)
\psline(1,1)(0,0)
\psline(5,-3)(1,1)
\psline[linecolor=blue](0,0)(-2,-2)
\psline(-2,-2)(-3,-3)
\psline(-3,-3)(1,-7)
\psline(1,-7)(2,-6)
\psline(3,-7)(2,-6)
\psline[linecolor=blue](4,-4)(0,0)
\psline(5,-5)(4,-4)
\psline(5,-5)(3,-7)
\psline[linecolor=blue](-2,-2)(0,-4)
\psline[linecolor=blue](3.5,-4.5)(4,-4)
\psline[linecolor=red](0,-4)(2,-6)
\psline[linecolor=red](2,-6)(3.5,-4.5)
\psline[linecolor=red](2.5,-3.5)(3.5,-4.5)
\psline[linecolor=red](2,-4)(2.5,-3.5)
\psline[linecolor=red](1,-3)(2,-4)
\psline[linecolor=red](1,-3)(0,-4)
\end{pspicture*}
\begin{pspicture*}(-3,-8)(6,3)
\psline(-1,1)(-3,-1)
\psline(0,0)(-1,1)
\psline(0,0)(1,1)
\psline(-3,-1)(-2,-2)
\psline[linecolor=red](-2,-2)(2,-6)
\psline[linecolor=red](2,-6)(4,-4)
\psline(4,-4)(5,-3)
\psline(-2,-2)(-3,-1)
\psline(5,-3)(1,1)
\psline(-3,-3)(-2,-2)
\psline[linecolor=blue](-1.5,-1.5)(0,0)
\psline[linecolor=blue](2,-2)(0,0)
\psline[linecolor=red](-1.5,-1.5)(-2,-2)
\psline[linecolor=red](-1.5,-1.5)(-0.5,-2.5)
\psline[linecolor=red](0,-2)(-0.5,-2.5)
\psline[linecolor=red](0,-2)(1,-3)
\psline[linecolor=red](2,-2)(1,-3)
\psline(-3,-3)(1,-7)
\psline(1,-7)(2,-6)
\psline(3,-7)(2,-6)
\psline[linecolor=red](4,-4)(2,-2)
\psline(5,-5)(4,-4)
\psline(5,-5)(3,-7)
\end{pspicture*}}
In the second picture we have drawn inside $Q_{F_d}$ the quiver $Q_{Y_d^*}$
  in red and the quiver $i_{F_d}(Q_{Y_d^*})$ in blue. In the third one, we
  have $Q_{T_d}$ in blue and inside it $Q_{w}$
  in red. Finally on the right picture we have the quiver
  $Q_{w_{q^d}}$ in red.

In terms of partitions, consider a Schubert class $[X(\lambda)]$,
where $\lambda$ is a partition whose Ferrers diagram is contained
in the rectangle $p\times (n-p)$. For $d\le \min(p,n-p)$, this 
class admits a $d$-Poincar{\'e} dual if and only if $\lambda_d=n-p$
and $\lambda_p\ge d$. Then $\lambda$ is uniquely defined by the 
partition $\mu$, whose diagram is contained in the rectangle 
$(p-d)\times (n-p-d)$ (the quiver $Q_{T_d}$), such that
$\mu_i=\lambda_{d+i}-d$. Let 
$\mu^*$ be the partition complementary to $\mu$ in that smaller
rectangle. Then the $d$-Poincar{\'e} dual $[X(\lambda)]$ is the
class $[X(\lambda_{q^d})]$, where the partition $\lambda_{q^d}$
is defined by 
$$\lambda_{q^d,d}=n-p, \qquad \lambda_{q^d,i}=d+\mu^*_{i-d}
\quad{\mathrm for}\; i>d.$$

(\i\i) For an isotropic Grassmannian $\G_{\omega}(n,2n)$
(resp. $\G_Q(n,2n)$), the picture similar to that for ordinary
Grassmannians is:

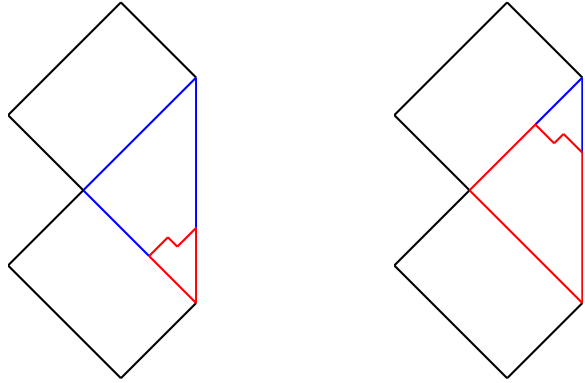
\psset{unit=.5cm}

\centerline{
\begin{pspicture*}(0,-3)(5,10)
\psline(3,8)(0,5)
\psline(0,5)(2,3)
\psline(3,8)(5,6)
\psline[linecolor=blue](5,6)(5,2)
\psline[linecolor=blue](3.75,1.25)(2,3)
\psline[linecolor=blue](2,3)(5,6)
\psline(5,0)(3,-2)
\psline(3,-2)(0,1)
\psline(0,1)(2,3)
\psline[linecolor=red](5,0)(3.75,1.25)
\psline[linecolor=red](3.75,1.25)(4.25,1.75)
\psline[linecolor=red](4.25,1.75)(4.5,1.5)
\psline[linecolor=red](4.5,1.5)(5,2)
\psline[linecolor=red](5,2)(5,0)
\end{pspicture*}
\hskip 2.5cm
\begin{pspicture*}(0,-3)(5,10)
\psline(3,8)(0,5)
\psline(0,5)(2,3)
\psline(3,8)(5,6)
\psline[linecolor=blue](5,6)(5,4)
\psline[linecolor=blue](5,6)(3.75,4.75)
\psline(5,0)(3,-2)
\psline(3,-2)(0,1)
\psline(0,1)(2,3)
\psline[linecolor=red](5,0)(2,3)
\psline[linecolor=red](2,3)(3.75,4.75)
\psline[linecolor=red](3.75,4.75)(4.25,4.25)
\psline[linecolor=red](4.25,4.25)(4.5,4.5)
\psline[linecolor=red](4.5,4.5)(5,4)
\psline[linecolor=red](5,4)(5,0)
\end{pspicture*}
}
\psset{unit=1cm}

Recall that the Schubert classes are
indexed by strict partitions $\lambda$ made of integers smaller 
or equal to $N$ (recall the notation of example \ref{exemple-grass1})-- 
this is often denoted 
$\lambda\subset\rho_N$,
where $\rho_N=(N,N-1,\ldots,2,1)$. 

The  Poincar{\'e}
dual of the Schubert class $[X(\lambda)]$ is the class
$[X(\lambda^*)]$, where $\lambda^*$ is the partition whose parts
complement the parts of $\lambda$ in the set $\{1,\ldots,N\}$. 
More generally, the Schubert class 
$[X(\lambda)]$ has a 
$d$-Poincar{\'e} dual $[X(\lambda_{q^d})]$ if and only if 
$\lambda_{\diamond} = N-\diamond$; in this case, denote by
$\mu\subset\rho_{N-\diamond}$ the partition defined by
$\mu_i = \lambda_{N-\diamond+i}$. Then we have
$$
(\lambda_{q^d})_i=
\left \{
\begin{array}{ccc}
N-i & \mbox{ if } & i \leq  \diamond\\
\mu^*_{i-\diamond} & \mbox{ if } & i > \diamond,
\end{array}
\right .
$$
where $\mu^*$ is the complement of $\mu$ inside $\rho_{N-\diamond}$.

(\i\i\i) If $X=E_6/P_1$ and $d=2$, then $[Y_d]=\sigma_8$ 
and the only Schubert class $[X(w)]$ such
that $\sigma_8*[X(w)]$ has a non trivial degree two term is
$[X(w)]=\sigma_{16}$. In this case we have
$$\sigma_8*\sigma_{16}=q^2\sigma_0.$$
Indeed, the $q^2$ term comes from Poincar{\'e} duality, we proved in
example \ref{sigma8} that all degree one invariants
$I_1(\sigma_8,\sigma_{16},\sigma_u)$ vanish and for dimension reasons,
there is no $q^0$ term.
\end{example}

The previous observation can be generalized as follows:

\begin{prop}
Let $d_{\max}$ be the maximal power of $q$ 
in the quantum product of two Schubert classes. Then 
we have the following formulae:
\begin{eqnarray}
\nonumber [\{{\rm pt}\}]*[\{{\rm pt}\}] &= &q^{d_{\max}}[Y_{d_{\max}}], \\
\nonumber [Y_{d_{\max}}^*]*[\{{\rm pt}\}]&= &q^{d_{\max}}[X].
\end{eqnarray}
  \end{prop}

 The values of $d_{\max}$ are the following. For a
    Grasmannian $\G(p,n)$, $d_{\max}=\min(p,n-p)$. For 
$\G_Q(n,2n)$, $d_{\max}=[n/2]$ and for $\G_{\omega}(n,2n)$,
$d_{\max}=n$. Finally, $d_{\max}=2$ for quadrics or the Cayley plane, 
while $d_{\max}=3$ for the Freudenthal variety.
    
\begin{proo}
 If $q^d[X(w)]$ appears in the product $[\{{\rm pt}\}]*[\{{\rm pt}\}]$,
then $X(w)$ must contain $Y_d$ by corollary \ref{gw-d}. In
particular we must have 
$$\dim(X)+d\cdot c_1(X)=2\dim(X)+\dim(X(w))\geq 2\dim(X)+\dim(Y_d).$$
But since the $d$-Poincar{\'e} dual to the class of a point is our
variety $T_d$, we have the relation
$$\dim(X)+d\cdot c_1(X)=\dim(X)+\codim(T_d)+\dim(Y_d).$$
Comparing with the previous inequality, we get
$\codim(T_d)\ge\dim(X)$, hence $T_d=X$ and  $d=d_{\max}$, and 
$\dim(X(w))=\dim(Y_d)$, thus $X(w)=Y_d$. 
In particular we only have the term
$q^{d_{\max}}[Y^*_{d_{\max}}]$ in $[\{{\rm pt}\}]*[\{{\rm pt}\}]$. 
The higher Poincar{\'e} duality 
implies that the coefficient is one and the first identity follows.

To prove the second one, using the first identity and the
associativity of the quantum product, we see that we only need
to prove that if 
$q^d$ appears in the product $[Y_{d_{\max}}^*]*[\{{\rm pt}\}]$, 
then $d=d_{\max}$. Let us set $v=w_{Y_{d_{\rm max}}}^*$
and $w=1$ so that $X(w)=\{{\rm pt}\}$. 
The Schubert variety $F_d(\wh)$ contains $F_d(\vh^\star)$ if and only if
there is an inclusion $Q_{F_d}(\vh^\star)\cap i_{F_d}(Q_{Z_d})\subset
Q_{F_d}(\wh)\cap i_{F_d}(Q_{Z_d})$ of quivers (cf Lemma
\ref{cas-gen2}). But both quivers $Q_{F_d}(\vh^\star)$ and
$Q_{F_d}(\wh)$ are contained in $i_{F_d}(Q_{Z_d})$ and if $d< d_{\rm
  max}$, the quiver of $F_d(\wh)$ does not contain the quiver of
$F_d(\vh^\star)$. The Lemma \ref{annul} gives the vanishing of $I_0([Y_{d_{\rm 
    max}}^*],[\{{\rm pt}\}],[F_d(\uh)])$ for any $u\in W_X$ and in
particular $q^d$ appears in $[Y_{d_{\max}}^*]*[\{{\rm pt}\}]$ if and
only if $d=d_{\max}$.
  \end{proo}

\subsectionplus{The smallest power in a quantum product of Schubert classes}

In \cite{FW}, Fulton and Woodward described the minimal power of
$q$ that can appear in the quantum product $[X(u)]*[X(v)]$ of two
Schubert classes. In this subsection
we give a new combinatorial description of this minimal power for
(co)minuscule homogeneous varieties. This generalizes the
reinterpretation of Fulton and Woodward's result by A. Buch in
\cite[Theorem 3]{Buch}.

Let $u\in W_X$, we define the element $\overline{u}\in W_X$ by
$Q_{\overline{u}}=i_X(i_{T_d}(Q_u\cap Q_{T_d})).$

\begin{prop}
  Let $X(u)$ and $X(v)$ be two Schubert subvarieties in $X$. If $q^d$
  appears in $[X(u)]*[X(v)]$ then the quiver $X(\overline{u}^*)$ is
  a subquiver of $X(v)$.
\end{prop}

\begin{proo}
The hypothesis that $q^d$ appears in the quantum product
$[X(u)]*[X(v)]$ is equivalent to the 
  existence of an element $w\in W_X$ such that $I_d([X(u)],
  [X(v)],[X(w)])\neq0$. By corollary \ref{gw-d}, this is equivalent to
  the non vanishing of the 
  classical invariant $I_0([F_d(\uh)],[F_d(\vh)],[F_d(\wh)])$. This in
  particular implies that the product $[F_d(\uh)]\cdot[F_d(\vh)]$ is
  non zero and, by Lemma \ref{annul}, that $F_d(\vh)\supset
  F_d(\uh^\star)$. Thanks to Lemma \ref{cas-gen2}, we have the inclusions
  $Q_{F_d}(\uh^\star)\cap i_{F_d}(Q_{Z_d})\subset Q_{F_d}(\vh)\cap
i_{F_d}(Q_{Z_d})$ of quivers. Because $Q_{F_d}(\uh^\star)$ is contained in
$Q_{F_d}(Z_d)$ and $Q_{F_d}(\vh)\cap i_{F_d}(Q_{Z_d})$ contains
$Q_{Z_d}$, this inclusion is equivalent to $Q_{F_d}(\uh^\star)\cap
Q_{T_d})\subset Q_{F_d}(\vh)\cap Q_{T_d}$.
This is equivalent to the
  inclusion $i_X(Q_{\overline{u}})\subset Q_v$ or $X(\overline{u}^*)\subset
  X(v)$. 
\end{proo}

The following corollary gives the smallest power of $q$ in a quantum
product. It is a generalization of lemma \ref{annul} giving a
condition for $q^0$ to appear in that product. It is also a
generalisation of theorem 3 in \cite{Buch}:

\begin{coro}\label{smallest}
  Let $X(u)$ and $X(v)$ be two Schubert subvarieties in $X$. The
  smallest power $q^d$ that appears in $[X(u)]*[X(v)]$ is the smallest
  $d$ such that 
$$Q_u\subset Q_{w_{Y_d}^*},\qquad Q_v\subset Q_{w_{Y_d}^*}, \qquad
i_X(Q_{\overline{u}})\subset Q_v.$$
\end{coro}

\begin{proo}
We already know by the vanishing criterion of corollary \ref{gw-d} and
the previous proposition that the conditions are necessary for the
product $[X(u)]*[X(v)]$  to have a $q^d$ term.
 
Conversely, let us denote by 
$\widetilde{u}$ the element in $W_X$ such that
$Q_{\widetilde{u}}=Q_u\cap Q_{T_d}$. 
Since $Q_{\widetilde{u}}$ is contained in $Q_u$, the product $H^{*a}*\sigma_u$
with $a=\dim X(u)-\dim X(\widetilde{u})$ contains
$\sigma_{\widetilde{u}}=\sigma_{w_{T_d}}\cdot\sigma_{\overline{u}}$.
Multiplying by 
$\sigma_{w_{Y_d}^*}$ gives by higher Poincar{\'e} duality that the product
$\sigma_{w_{Y_d}^*}*H^{*a}*\sigma_u$ contains
$\sigma_{w_{Y_d}^*}*\sigma_{w_{T_d}}\cdot\sigma_{\overline{u}} =
q^d\sigma_{\overline{u}}$. Finally, the product 
$\sigma_{w_{Y_d}^*}*H^{*a}*\sigma_u*\sigma_v$
contains $q^d\sigma_{\overline{u}}*\sigma_v$ which has a non zero
$q^d$ term by Lemma \ref{annul} and the hypothesis
$i_X(Q_{\overline{u}})\subset Q_v$.
By non negativity
of the invariants, the product $\sigma_u*\sigma_v$ must contain a
power in $q$ smaller or equal to $d$.
\end{proo}

\sectionplus{The quantum cohomology of the exceptional 
minuscule varieties}

\label{section_quantum_exceptional}

In this section we apply our quantum Chevalley formula to the 
exceptional minuscule varieties. In turns out that, together with our
computation in example \ref{sigma8}, this suffices to 
deduce the whole quantum Chow ring from the classical one. 

\subsectionplus{The quantum Chow ring of the Cayley plane}
\label{subsection_quantum_cayley}

The quantum Chevalley formula for the Cayley plane $\OO \PP^2$ 
can be conveniently visualized on the Hasse diagram. The
following follows from subsection \ref{subsection_higher_quantum} and
example \ref{exemple_incidence_e6}. The Schubert 
classes having a 1-Poincar{\'e} dual class are those in the Bruhat 
interval $[\s_{16},\s''_{11}]$, represented in blue on the picture 
below. (Note that $\s''_{11}$ is the class of the Schubert variety 
$T_1$, with the notations of the previous
subsection.)
For such a class $\sigma$, the 1-Poincar{\'e}
dual is obtained by first applying the obvious symmetry in this
diagram, and then the usual Poincar{\'e} duality. This means that 
the $q$-term in $\sigma*H$ is the Schubert class  corresponding to
$\sigma$ in the isomorphic interval $[\s''_5,[X]]$, in red on the 
picture. 

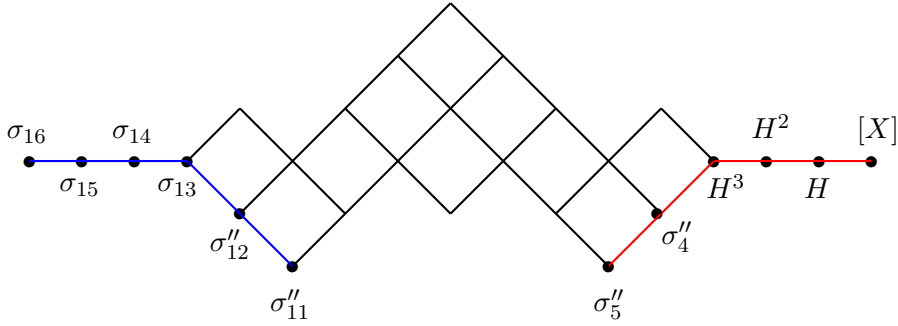
\psset{xunit=3.5mm}
\psset{yunit=3.5mm}
\centerline{
\begin{pspicture*}(-8,-2)(30,20)
\put(-2.2,2.7){$\bullet$}
\put(-1.5,2.7){$\bullet$}
\put(-.8,2.7){$\bullet$}
\put(-.1,2.7){$\bullet$}
\put(.6,2){$\bullet$}
\put(1.3,1.3){$\bullet$}
\put(6.9,2.7){$\bullet$}
\put(7.6,2.7){$\bullet$}
\put(8.3,2.7){$\bullet$}
\put(9,2.7){$\bullet$}
\put(6.15,2){$\bullet$}
\put(5.5,1.3){$\bullet$}
\put(.3,1.6){$\sigma_{12}''$}
\put(-.4,2.4){$\sigma_{13}$}
\put(-1,3.1){$\sigma_{14}$}
\put(-1.7,2.4){$\sigma_{15}$}
\put(-2.4,3.1){$\sigma_{16}$}
\put(8.2,2.3){$H$}
\put(8.9,3.1){$[X]$}
\put(7.5,3.1){$H^2$}
\put(6.9,2.3){$H^3$}
\put(6.3,1.7){$\s_4''$}
\put(5.4,.8){$\s_5''$}
\put(1.1,.8){$\s_{11}''$}
\psline[linecolor=red](20,8)(26,8)
\psline[linecolor=blue](-6,8)(0,8)
\psline(18,10)(20,8)
\psline(14,6)(18,10)
\psline[linecolor=red](16,4)(20,8)
\psline(12,8)(14,10)
\psline(0,8)(2,10)
\psline(2,6)(10,14)
\psline(10,14)(18,6)
\psline(4,4)(12,12)
\psline(10,6)(12,8)
\psline[linecolor=blue](0,8)(4,4)
\psline(2,10)(6,6)
\psline(6,10)(10,6)
\psline(8,12)(16,4)
\end{pspicture*}}

Here is the quantum version of Proposition \ref{e6_classique}.

\begin{theo}
\label{e6_quantique}
Let ${\cal H}_q := \Z[h,s,q]/(3hs^2-6h^5s+2h^9,s^3-12h^8s+5h^{12}-q)$.
Mapping $h$ to $H$, $s$ to $\sigma'_4$ and preserving $q$ yields
an isomorphism of graded algebras 
$${\cal H}_q \simeq QA^*(\OO\PP^2).$$
\end{theo}

\begin{proo}
By proposition \ref{e6_classique} and \cite[proposition 11]{fulton},
it is enough to show that the displayed relations hold in
$QA^*(\OO\PP^2)$. Recall that the index of the Cayley plane is twelve. 
The relation
$3H*\sigma^{*2} - 6 H^{*5}*\sigma +2 H^{*9} = 0$ holds because 
its degree is strictly less than 12 and it holds 
in the classical Chow ring. 

The
Fano variety of lines through a given point can be identified 
with the spinor variety $G_Q(5,10)$ in its minimal embedding 
(the projectivization of a half-spin representation). In particular 
its degree equals $12$ (see \cite[3.1]{landsberg}). Applying proposition
\ref{h},
we get the relation  
\begin{eqnarray}\label{12}
H^{*12}=H^{12}+12q,
\end{eqnarray} 
which we could also deduce from the quantum Chevalley formula. 
Now we use the result of example \ref{sigma8}, according to which the
multiplication of $\sigma_8$ by any class of degree four does not 
require any quantum correction. Since
$\sigma_8 = \sigma^2 + 2 H^4\sigma - H^8$, we get:
$$
\begin{array}{rcccl}
\sigma^{*3} + 2 H^{*4}*\sigma^{*2}-H^{*8}*\sigma & = &
\sigma'_{12} & = & 15 H^8\sigma - 19/3 H^{*12} + 76q,\\
H^{*4}*\sigma^{*2} + 2 H^{*8} * \sigma - H^{*12} & = &
\sigma'_{12}+\sigma''_{12} & = & 4 H^8\sigma -5/3 H^{*12} + 20q.
\end{array}
$$
Indeed, for the first line we have used that 
$\sigma_8*\sigma'_{4} =\sigma_8\cdot\sigma'_{4} = \sigma'_{12}$, 
the expression for $\sigma'_{12}$ obtained in the proof
of Proposition \ref{e6_classique}, and the identity (\ref{12}). For
the second line we have used that $\sigma_8H^4=\sigma'_{12}+\sigma''_{12}$.

Eliminating $H^8\sigma$, we get
$4\sigma^{*3} - 7 H^{*4} * \sigma^{*2} - 
34 H^{*8}*\sigma + 46/3 H^{*12} - 4q = 0$.
Taking into account the  relation already proved, this yields
$4\sigma^{*3} - 48 H^{*8}*\sigma + 20 H^{*12} - 4 q =0$, as claimed.
\end{proo}

We can be more specific about the quantum multiplication of Schubert 
classes. The first quantum corrections appear in degree twelve. In 
this degree, the only cases which do not  follow directly from the
quantum Chevalley formula are the products of a degree eight class by a
degree four class. We have 
$$\sigma_8'*\sigma_4''=\sigma_8'\sigma_4''+q, \qquad 
\sigma_8''*\sigma_4'=\sigma_8''\sigma_4'+q,$$
while the other products have no quantum correction.

In fact, to prove this, let us first compare the quantum monomials
$H^{*i} * \sigma^{*j}$ of degree 12 with the corresponding classical products.
The classical Chevalley formula gives
$$
H^3 \cdot \sigma^2 = 2 \sigma''_{11} + 6 \sigma'_{11}
\mbox{ and }
H^7 \cdot \sigma = 5 \sigma''_{11} + 14 \sigma'_{11}.
$$
Therefore, our quantum Chevalley formula yields
$$
H^4 * \sigma^2 = H^4 \cdot \sigma^2 + 2q
\mbox{ and }
H^8 * \sigma = H^8 \cdot \sigma + 5q.
$$
Since $\sigma \cdot \sigma_8 = \sigma * \sigma_8$, it follows
that $\sigma^{*3} = \sigma^3 + q$.
The claims about quantum products of Schubert classes of degree 4 and
8 follow directly.

\lpara

It is then easy, inductively, to obtain the following formulae:
\begin{eqnarray} \nonumber
\s'_{12} &=& \s^3+2H^4\s^2-H^8\s, \\ \nonumber
\s''_{12} &=& -\s^3-H^4\s^2+3H^8\s-H^{12}, \\ \nonumber
\s_{13} &=& H\s^3+2H^5\s^2-H^9\s, \\ \nonumber
\s_{14} &=& 2H^6\s^2+11H^{10}\s-5H^{14}, \\ \nonumber
\s_{15} &=& -H^3\s^3+2H^7\s^2+23H^{11}\s-10H^{15}, \\ \nonumber
\s_{16} &=& \s^4-2H^4\s^3-10H^8\s^2+H^{16}.
\end{eqnarray}
We have deliberately omitted the signs for the quantum product. 
In fact we have:

\begin{prop}
The previous Giambelli type formulas for the Schubert classes of 
the Cayley plane hold in the classical as well as in the quantum 
Chow ring.
\end{prop}

Of course in degree smaller than twelve, this is also the case of the 
Giambelli type formulas given in the proof of Proposition
\ref{e6_classique}. The coincidence of classical and quantum Giambelli
formulas already appeared for minuscule homogeneous varieties of
classical type (see \cite{bertram,KT2}).

\subsectionplus{The quantum Chow ring of the Freudenthal variety}

\label{subsection_quantum_freudenthal}

The quantum Chevalley formula for the Freudenthal variety $E_7/P_7$
can again easily be visualized on the Hasse diagram. The quiver of
the Fano variety of lines in $E_7/P_7$ looks as follows
(the quiver of $E_7/P_7$ is in blue and the quiver of $Z$ is red):

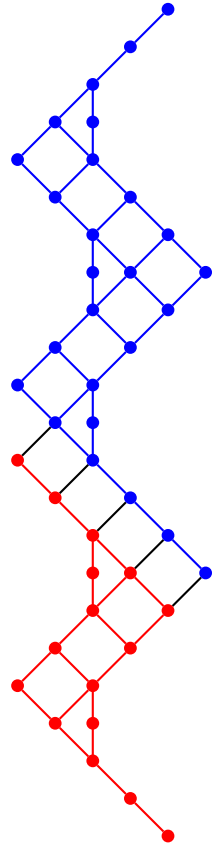
\psset{unit=5mm}
\begin{pspicture*}(-14,.5)(5.50,24.50)

\psline[linecolor=blue](3.88,23.88)(3.12,23.12)
\psline[linecolor=blue](2.88,22.88)(2.12,22.12)
\psline[linecolor=blue](1.88,21.88)(1.12,21.12)
\psline[linecolor=blue](2.0,21.88)(2.0,21.12)
\psline[linecolor=blue](0.88,20.88)(0.12,20.12)
\psline[linecolor=blue](1.12,20.88)(1.88,20.12)
\psline[linecolor=blue](2.0,20.88)(2.0,20.12)
\psline[linecolor=blue](0.12,19.88)(0.88,19.12)
\psline[linecolor=blue](1.88,19.88)(1.12,19.12)
\psline[linecolor=blue](2.12,19.88)(2.88,19.12)
\psline[linecolor=blue](1.12,18.88)(1.88,18.12)
\psline[linecolor=blue](2.88,18.88)(2.12,18.12)
\psline[linecolor=blue](3.12,18.88)(3.88,18.12)
\psline[linecolor=blue](2.0,17.88)(2.0,17.12)
\psline[linecolor=blue](2.12,17.88)(2.88,17.12)
\psline[linecolor=blue](3.88,17.88)(3.12,17.12)
\psline[linecolor=blue](4.12,17.88)(4.88,17.12)
\psline[linecolor=blue](2.0,16.88)(2.0,16.12)
\psline[linecolor=blue](2.88,16.88)(2.12,16.12)
\psline[linecolor=blue](3.12,16.88)(3.88,16.12)
\psline[linecolor=blue](4.88,16.88)(4.12,16.12)
\psline[linecolor=blue](1.88,15.88)(1.12,15.12)
\psline[linecolor=blue](2.12,15.88)(2.88,15.12)
\psline[linecolor=blue](3.88,15.88)(3.12,15.12)
\psline[linecolor=blue](0.88,14.88)(0.12,14.12)
\psline[linecolor=blue](1.12,14.88)(1.88,14.12)
\psline[linecolor=blue](2.88,14.88)(2.12,14.12)
\psline[linecolor=blue](0.12,13.88)(0.88,13.12)
\psline[linecolor=blue](1.88,13.88)(1.12,13.12)
\psline[linecolor=blue](2.0,13.88)(2.0,13.12)
\psline(0.88,12.88)(0.12,12.12)
\psline[linecolor=blue](1.12,12.88)(1.88,12.12)
\psline[linecolor=blue](2.0,12.88)(2.0,12.12)
\psline[linecolor=red](0.12,11.88)(0.88,11.12)
\psline(1.88,11.88)(1.12,11.12)
\psline[linecolor=blue](2.12,11.88)(2.88,11.12)
\psline[linecolor=red](1.12,10.88)(1.88,10.12)
\psline(2.88,10.88)(2.12,10.12)
\psline[linecolor=blue](3.12,10.88)(3.88,10.12)
\psline[linecolor=red](2.0,9.88)(2.0,9.12)
\psline[linecolor=red](2.12,9.88)(2.88,9.12)
\psline(3.88,9.88)(3.12,9.12)
\psline[linecolor=blue](4.12,9.88)(4.88,9.12)
\psline[linecolor=red](2.0,8.88)(2.0,8.12)
\psline[linecolor=red](2.88,8.88)(2.12,8.12)
\psline[linecolor=red](3.12,8.88)(3.88,8.12)
\psline(4.88,8.88)(4.12,8.12)
\psline[linecolor=red](1.88,7.88)(1.12,7.12)
\psline[linecolor=red](2.12,7.88)(2.88,7.12)
\psline[linecolor=red](3.88,7.88)(3.12,7.12)
\psline[linecolor=red](0.88,6.88)(0.12,6.12)
\psline[linecolor=red](1.12,6.88)(1.88,6.12)
\psline[linecolor=red](2.88,6.88)(2.12,6.12)
\psline[linecolor=red](0.12,5.88)(0.88,5.12)
\psline[linecolor=red](1.88,5.88)(1.12,5.12)
\psline[linecolor=red](2.0,5.88)(2.0,5.12)
\psline[linecolor=red](1.12,4.88)(1.88,4.12)
\psline[linecolor=red](2.0,4.88)(2.0,4.12)
\psline[linecolor=red](2.12,3.88)(2.88,3.12)
\psline[linecolor=red](3.12,2.88)(3.88,2.12)
\pscircle*[linecolor=blue](0.0,20.0){.17}
\pscircle*[linecolor=blue](0.0,14.0){.17}
\pscircle*[linecolor=red](0.0,12.0){.17}
\pscircle*[linecolor=red](0.0,6.0){.17}
\pscircle*[linecolor=blue](1.0,21.0){.17}
\pscircle*[linecolor=blue](1.0,19.0){.17}
\pscircle*[linecolor=blue](1.0,15.0){.17}
\pscircle*[linecolor=blue](1.0,13.0){.17}
\pscircle*[linecolor=red](1.0,11.0){.17}
\pscircle*[linecolor=red](1.0,7.0){.17}
\pscircle*[linecolor=red](1.0,5.0){.17}
\pscircle*[linecolor=blue](2.0,22.0){.17}
\pscircle*[linecolor=blue](2.0,21.0){.17}
\pscircle*[linecolor=blue](2.0,20.0){.17}
\pscircle*[linecolor=blue](2.0,18.0){.17}
\pscircle*[linecolor=blue](2.0,17.0){.17}
\pscircle*[linecolor=blue](2.0,16.0){.17}
\pscircle*[linecolor=blue](2.0,14.0){.17}
\pscircle*[linecolor=blue](2.0,13.0){.17}
\pscircle*[linecolor=blue](2.0,12.0){.17}
\pscircle*[linecolor=red](2.0,10.0){.17}
\pscircle*[linecolor=red](2.0,9.0){.17}
\pscircle*[linecolor=red](2.0,8.0){.17}
\pscircle*[linecolor=red](2.0,6.0){.17}
\pscircle*[linecolor=red](2.0,5.0){.17}
\pscircle*[linecolor=red](2.0,4.0){.17}
\pscircle*[linecolor=blue](3.0,23.0){.17}
\pscircle*[linecolor=blue](3.0,19.0){.17}
\pscircle*[linecolor=blue](3.0,17.0){.17}
\pscircle*[linecolor=blue](3.0,15.0){.17}
\pscircle*[linecolor=blue](3.0,11.0){.17}
\pscircle*[linecolor=red](3.0,9.0){.17}
\pscircle*[linecolor=red](3.0,7.0){.17}
\pscircle*[linecolor=red](3.0,3.0){.17}
\pscircle*[linecolor=blue](4.0,24.0){.17}
\pscircle*[linecolor=blue](4.0,18.0){.17}
\pscircle*[linecolor=blue](4.0,16.0){.17}
\pscircle*[linecolor=blue](4.0,10.0){.17}
\pscircle*[linecolor=red](4.0,8.0){.17}
\pscircle*[linecolor=red](4.0,2.0){.17}
\pscircle*[linecolor=blue](5.0,17.0){.17}
\pscircle*[linecolor=blue](5.0,9.0){.17}
\end{pspicture*}

The Schubert 
classes having a 1-Poincar{\'e} dual class are those in the Bruhat 
interval $[\s_{27},\s_{17}]$, represented in blue on the picture 
below, $\sigma_{17}$ being the class of $T_1$. For such a class
$\sigma$, the 1-Poincar{\'e}
dual is obtained by first applying the obvious symmetry of this
interval, and then the usual Poincar{\'e} duality. So
the $q$-term in $\sigma*H$ is the Schubert class  corresponding to
$\sigma$ in the isomorphic interval $[\s_{10},[X]]$, in red on the 
picture. 

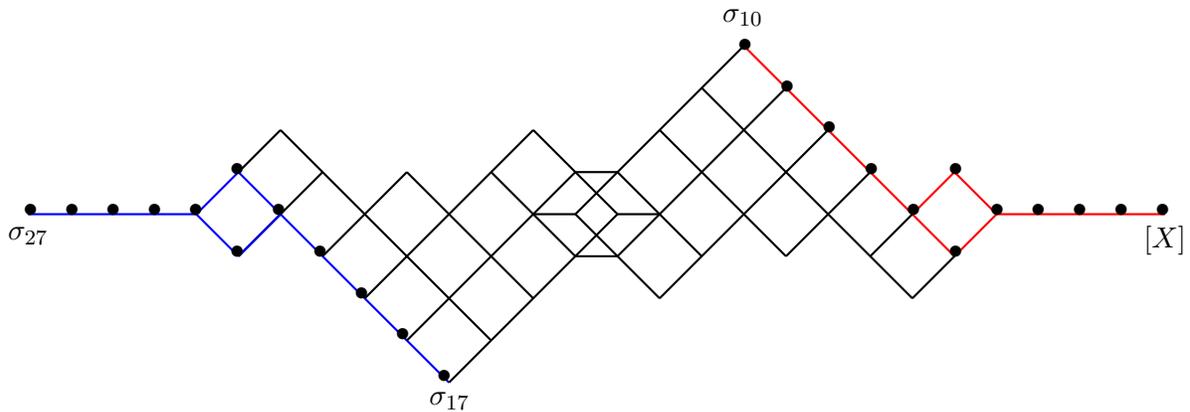
\psset{unit=1cm}
\psset{xunit=2.8mm}
\psset{yunit=2.8mm}
\centerline{
\begin{pspicture*}(-10,-3)(51,20)
\psline[linecolor=blue](-8,8)(0,8)
\psline[linecolor=blue](0,8)(2,10)
\psline[linecolor=blue](0,8)(2,6)
\psline[linecolor=blue](2,10)(12,0)
\psline(2,6)(4,8)
\psline(2,10)(4,12)
\psline(4,8)(6,10)
\psline[linecolor=blue](2,6)(4,8)
\psline(6,6)(10,10)
\psline(8,4)(16,12)
\psline(10,2)(18,10)
\psline(12,0)(20,8)
\psline(4,12)(14,2)
\psline(10,10)(16,4)
\psline(14,10)(18,6)
\psline(16,12)(20,8)
\psline(16,8)(18,8)
\psline(18,6)(20,6)
\psline(18,10)(20,10)
\psline(20,8)(22,8)
\psline(18,8)(26,16)
\psline(18,8)(22,4)
\psline(20,6)(28,14)
\psline(22,4)(30,12)
\psline(20,10)(24,6)
\psline(22,12)(28,6)
\psline(24,14)(34,4)
\psline[linecolor=red](26,16)(36,6)
\psline(28,6)(32,10)
\psline(32,6)(34,8)
\psline[linecolor=red](34,8)(36,10)
\psline[linecolor=red](36,10)(38,8)
\psline[linecolor=red](38,8)(46,8)
\psline[linecolor=red](36,6)(38,8)
\psline(34,4)(36,6)
\multiput(-2.3,2.2)(.55,0){5}{$\bullet$}
\multiput(.45,2.75)(.55,-.55){6}{$\bullet$}
\put(.45,1.65){$\bullet$}
\multiput(10.55,2.2)(.55,0){5}{$\bullet$}
\multiput(10,1.65)(-.56,.55){6}{$\bullet$}
\put(10,2.75){$\bullet$}
\put(7,4.8){$\sigma_{10}$}
\put(12.6,1.8){$[X]$}
\put(-2.5,1.9){$\sigma_{27}$}
\put(3.1,-.3){$\sigma_{17}$}
\end{pspicture*}}

As for the Cayley plane, 
we can deduce all the quantum products of Schubert
classes in degree
$18$. For example, let $c,d$ be classes of degrees $5$ and $13$. Since
the product by the hyperplane class defines an isomorphism between
$A^{12}$ and $A^{13}$, we can write $d=He$ for some class $e$ of 
degree $12$, either in the classical or the quantum Chow ring. 
Using the associativity of the quantum product, we get 
$$c*d=c*(He)=c*H*e=(c*e)*H=(ce)*H.$$
This can be computed from the classical intersection product and 
the quantum Chevalley formula. 

A priori, this method does not work if we want to compute the quantum
product of two classes of degree $9$. Indeed the product with the 
hyperplane class does not define an isomorphism between
$A^{8}$ and $A^{9}$. Nevertheless, the fact that the quantum
correction of a product like $\s_9*\s'_9$ is a {\it non negative} 
integer multiple of $q$ allows us to carry the computation over. 
For example, we have 
$$(\s'_9+\s''_9)*\s_9=H\s''_8*\s_9=(\s''_8\s_9)*H=(2\s'_{17}+5\s''_{17})*H.$$
Since the class $\s_{17}$ does not appear, the quantum Chevalley formula 
shows that there is no quantum correction in this product, and therefore
there is none either in the two products $\s'_9*\s_9$ and $\s''_9*\s_9$. 

Our conclusion is the following:

\begin{prop}\label{e7q}
Among the quantum products of two Schubert classes of total degree $18$, 
the only ones that require a quantum correction are 
\begin{eqnarray}
\nonumber \s_{17}*H &=& \s_{17}H+q, \\
\nonumber \s'_{13}*\s'_5 &=& \s'_{13}\s'_5+q, \\
\nonumber \s''_{13}*\s''_5 &=& \s''_{13}\s''_5+q, \\
\nonumber \s_9*\s_9 &=& (\s_9)^2+q, \\
\nonumber \s'_9*\s'_9 &=& (\s'_9)^2+q, \\
\nonumber \s''_9*\s''_9 &=& (\s''_9)^2+q.
\end{eqnarray}
\end{prop}

This completely determines the quantum Chow ring. The quantum version
of Theorem \ref{e7_classique} is the following. 

\begin{theo}
\label{e7_quantique}
Let ${\cal H} = \ZZ[h,s,t,q]/(s^2 - 10sh^5 + 2th + 4h^{10},
2st - 12sh^9 + 2th^5 + 5h^{14}, t^2 + 922 sh^{13} - 198 th^9 - 385 h^{18}-q)$.
Mapping $h$ to $H$, $s$ to $\sigma'_5$, $t$ to $\sigma_9$ and preserving $q$
yields an isomorphism of graded algebras 
$${\cal H}_q \simeq QA^*(E_7/P_7).$$
\end{theo}

\begin{proo} As for the proof of Theorem \ref{e6_quantique}, we just
  need to prove that the displayed relations hold in $QA^*(E_7/P_7)$.
This is clear for the first two, since they are of degree smaller that
18 and they hold in the classical Chow ring by Theorem
\ref{e7_classique}. 
That the third equation also holds is a direct consequence of 
Proposition \ref{e7q}.
\end{proo}

As we did for the Cayley plane, we could also derive Giambelli 
type formulas for the Schubert classes, holding both in the classical
and the quantum Chow ring.

\bigskip\noindent
Pierre-Emmanuel {\sc Chaput}, Laboratoire de Math\'ematiques Jean Leray, 
UMR 6629 du CNRS, UFR Sciences et Techniques,  2 rue de la Houssini\`ere, BP
92208, 44322 Nantes cedex 03, France. 

\noindent {\it email}: pierre-emmanuel.chaput@math.univ-nantes.fr

\medskip\noindent
Laurent {\sc Manivel}, 
Institut Fourier, UMR 5582 du CNRS,  Universit\'e de Grenoble I, 
BP 74, 38402 Saint-Martin d'H\'eres, France.
 
\noindent {\it email}: Laurent.Manivel@ujf-grenoble.fr

\medskip\noindent
Nicolas {\sc Perrin},  Institut de Math\'ematiques,  
Universit\'e Pierre et Marie Curie, Case 247, 4 place Jussieu,  
75252 PARIS Cedex 05, France.

\noindent {\it email}: nperrin@math.jussieu.fr
\end{document}